%% file: main.tex
\PassOptionsToPackage{dvipsnames}{xcolor} %%xcoloroptioncrashsolver

\documentclass[11pt,letterpaper]{amsart}
\usepackage[a4paper, left=1in, right=1in, top=1.1in, bottom=1.1in]{geometry}

\usepackage[english]{babel}
\usepackage[utf8]{inputenc}
\usepackage{amsmath}
\usepackage{graphicx}
\usepackage{amssymb}
\usepackage{stmaryrd}
\usepackage{amsthm}
\usepackage{tikz-cd}
\usepackage{mathrsfs}
\usepackage{csquotes}
\usepackage[colorinlistoftodos]{todonotes}
\usepackage{yfonts}
\usepackage{ dsfont }
\usepackage{xcolor}
\usepackage{multicol}
\usepackage{indentfirst}
\usepackage[hyphens]{url}
\usepackage{hyperref}
\usepackage[hyphenbreaks]{breakurl}
\usepackage{subfiles}
\usepackage{longtable}
\usepackage{float}

 \usepackage[
    backend=biber,
    style=alphabetic, maxbibnames=99, url=false
  ]{biblatex}
\usepackage{xurl}
\hypersetup{
    colorlinks,
    linkcolor={RoyalBlue},
    citecolor={ForestGreen},
    urlcolor={RoyalPurple}
}
\hypersetup{breaklinks=true}

\usepackage[inline]{enumitem}   
\makeatletter

\newcommand{\inlineitem}[1][]{%
\ifnum\enit@type=\tw@
    {\descriptionlabel{#1}}
  \hspace{\labelsep}%
\else
  \ifnum\enit@type=\z@
       \refstepcounter{\@listctr}\fi
    \quad\@itemlabel\hspace{\labelsep}%
\fi}
\makeatother
\allowdisplaybreaks

 \usepackage[
    backend=biber,
    style=alphabetic, maxbibnames=99, url=false
  ]{biblatex}
\usepackage[capitalise]{cleveref}

\title[Extremal effective curves and non-semiample line bundles on $\M{g}{n}$]{Extremal effective curves and\\ non-semiample line bundles on $\M{g}{n}$}

\author{Daebeom Choi}
\address{Department of Mathematics\\
    University of Pennsylvania\\
    Philadelphia, PA 19104-6395}
\email{dbchoi@sas.upenn.edu}

\date{\today}
\theoremstyle{definition}
\newtheorem{thm}{Theorem}[section]

\newtheorem{lem}[thm]{Lemma}

\newtheorem{prop}[thm]{Proposition}
\newtheorem{qes}[thm]{Question}

\newtheorem{defn}[thm]{Definition}
\newtheorem{eg}[thm]{Example}

\newtheorem{conj}[thm]{Conjecture}

\newtheorem{cor}[thm]{Corollary}

\newtheorem{rmk}[thm]{Remark}

\newcommand{\R}{\mathbb{R}}

\newcommand{\Z}{\mathbb{Z}}
\newcommand{\N}{\mathbb{N}}
\newcommand{\Q}{\mathbb{Q}}

\newcommand{\K}{\text{Knu}}
\newcommand{\Pic}{\text{Pic}}
\newcommand{\HH}{\text{H}}
\newcommand{\A}{\text{A}}

\newcommand{\M}[2]{\overline{\rm{M}}_{#1, #2}}
\newcommand{\MM}[2]{\overline{\mathcal{M}}_{#1, #2}}

\newcommand{\Mg}[1]{\overline{\rm{M}}_{#1}}

\newcommand{\NE}[1]{\overline{\text{NE}}_1(#1)}

\DeclareSymbolFont{yhlargesymbols}{OMX}{yhex}{m}{n}
\DeclareMathAccent{\widetriangle}{\mathord}{yhlargesymbols}{"E6} %% widetriangle

\begin{document}
    
\begin{abstract}
We develop a new method for establishing the extremality in the closed cone of effective curves on the moduli space of curves and determine the extremality of many boundary $1$-strata. As a consequence, by using a general criterion for non-semiampleness which extends Keel's argument, we demonstrate that a substantial portion of the cone of nef divisors of $\M{g}{n}$ is not semiample. As an application, we construct the first explicit example of a non-contractible extremal ray of the closed cone of effective curves on $\M{3}{n}$. Our method relies on two main ingredients: (1) the construction of a new collection of nef divisors on $\M{g}{n}$, and (2) the identification of a tractable inductive structure on the Picard group, arising from Knudsen’s construction of $\M{g}{n}$.
\end{abstract}

\maketitle

\section{Introduction}\label{sec:intro}

\subfile{sections/introduction}

\section{Preliminaries}\label{sec:prelim}

\subfile{sections/prelim}

\section{Non-semiampleness Criterion}\label{sec:nonsemi}

\subfile{sections/nonsemi}

\section{Semiampleness of \texorpdfstring{$\psi$}{TEXT} classes on \texorpdfstring{$\M{1}{n}$}{TEXT}}\label{sec:psisemi}

\subfile{sections/psisemi}

\section{Examples}\label{sec:examples}

\subfile{sections/examples}

\section{The case of \texorpdfstring{$\M{0}{n}$}{TEXT} and \texorpdfstring{$\M{1}{n}$}{TEXT}}\label{sec:case}

\subfile{sections/case01}

\section{New families of nef divisors on \texorpdfstring{$\M{g}{n}$}{TEXT}}\label{sec:nef}

\subfile{sections/newnef}

\section{Divisors on Knudsen's construction}\label{sec:Knudsen}

\subfile{sections/Knudsen}

\section{Method for proving regular extremality}\label{subsec:method}

\subfile{sections/strategy}

\section{Extremality of F-curves}\label{sec:extremality}

\subfile{sections/extremality}

\section{Proof of \texorpdfstring{\cref{thm:main1}}{\texttwoinferior} and \texorpdfstring{\cref{cor:main1cor}}{\texttwoinferior}} \label{subsec:proof}

\subfile{sections/mainproof}

\section{Further Questions}\label{sec:discussion}

\subfile{sections/discussion}

\printbibliography

\end{document}

%% file: sections/introduction.tex
In this paper, we investigate the birational geometry of $\M{g}{n}$, the moduli space of stable, $n$-pointed genus $g$ curves, by introducing a new method to establish extremality in the cone of effective curves. This approach provides new insights into the dual cone of nef divisors and, in particular, illuminates the subtle nature of semiample divisors on \(\M{g}{n} \).

By Stein factorization, the investigation of morphisms from $\M{g}{n}$ to projective varieties can be reduced to the study of contractions of $\M{g}{n}$. Mori theory offers a framework for such analysis by associating to each contraction $f: X \to Y$ of projective varieties the relative cone of curves $\NE{f}$, which corresponds to a face of the closed cone of curves $\NE{X}$. In certain favorable cases — such as the case where $X$ is a log Fano variety defined over a field of characteristic zero — this theory yields a well-behaved bijection between contractions and faces.

Two fundamental issues can obstruct a correspondence between faces of the effective cone of curves and projective contractions. First, $\NE{X}$ may be non-polyhedral, having infinitely many faces. Second, it may happen that not every face of $\NE{X}$ corresponds to a contraction, since the cone of nef divisors may contain divisors that are not semiample. While conjecturally the first issue will not arise on \( \M{g}{n} \) (cf. \cref{conj:F}), our first main theorem shows that the second phenomenon is especially pronounced for $\M{g}{n}$: a large portion of the nef cone consists of non-semiample divisors.

\begin{thm}\label{thm:main1}
Assume that the characteristic of the base field is $0$, and that either $g=2$, $n\geq 2$, or $g\geq 3$, $n\geq 1$. Then there is a codimension $\lfloor \frac{g}{2} \rfloor$ face $F$ of the nef cone of $\M{g}{n}$ such that its general element is not semiample. More precisely, there is a codimension one linear subspace $E$ of $F$ such that any nef divisor in $F\setminus E$ is not semiample.
\end{thm}

If $g=2$ or $3$, then $F$ has codimension one, and for all $g$, the codimension is constant as $n\to \infty$. This indicates that the non-semiample region of the nef cone occupies a non-negligible portion of it. \cref{thm:main1} implies the following, more directly related to the second issue:

\begin{cor}\label{cor:main1cor}
Assume that the characteristic of the base field is $0$ and $g\ge 3$, $n\ge 1$. Then there exists a dimension $\lfloor \frac{g}{2} \rfloor$ face $C$ of $\NE{\M{g}{n}}$ which is not contractible, i.e. there is no projective contraction of $\M{g}{n}$ whose relative cone of curves is $C$.
\end{cor}

Hence, in particular, there exists a noncontractible extremal ray of $\NE{\M{3}{n}}$. However, by \cref{thm:M3ncont}, this extremal ray of $\NE{\M{3}{n}}$ is contractible in positive characteristic, and the contraction admits a very simple description. This highlights the dependence of the geometry of $\M{g}{n}$ on the characteristic of the base field. We will return to this point at the end of the introduction and in \cref{sec:discussion}.

The construction of nef but non-semiample line bundles on $\M{g}{n}$ dates back at least to Keel \cite{Ke99}, who proved that in characteristic $0$, $\psi$-classes are not semiample in general. In addition, we prove in \cref{thm:psisemi} that $\psi$-classes are semiample on $\M{1}{n}$, and thus we completely determine the cases in which $\psi$-classes are semiample (\cref{cor:psisem}). Moreover, in \cref{sec:case}, we explain that they play a significant role in contracting F-curves on $\M{1}{n}$ (cf.~\cref{rmk:cont}).

To prove \cref{thm:main1}, we utilize a non-semiampleness criterion, \cref{thm:semin}. The proof of \cref{thm:semin} follows Keel’s method, employing arguments similar to those of \cite{Cho23}. One major difficulty in the proof of \cref{thm:main1} lies in computing the dimension of the space of nef divisors satisfying the criterion of \cref{thm:semin}. Surprisingly, this difficulty is intimately related to the first issue of Mori theory discussed above.

The irreducible components of $1$-dimensional boundary strata on $\M{g}{n}$, so-called \textbf{F-curves}, are classified into six types in \cite{GKM02}. The \textbf{F-conjecture} asserts that $\NE{\M{g}{n}}$ is generated by F-curves, and hence is polyhedral. There has been extensive work on the F-conjecture (e.g., \cite{KM13, GKM02, GF03, Gib09, Lar11, Fe15, MS19, Fe20, FM25}), but still it remains wide open. We will recall the relevant background in \cref{sec:prelim}, including the notation for F-curves.

We define an extremal ray of a cone to be a \textbf{regular extremal} if the corresponding face of the dual cone has codimension $1$. Note that if the cone is polyhedral, then every extremal ray is regular, but this is not true in general (cf. \cref{rmk:regext}). This notion was implicitly used in \cite{Mul21}, where the author proved that the cone of moving curves of $\M{g}{n}$, for $g,n\ge 2$, is non-polyhedral (and hence the cone of effective divisors is also non-polyhedral). This was done by constructing a non-regular extremal ray \cite[Theorem 1.1]{Mul21}.

The following result shows that many extremal F-curves are indeed regular, providing further evidence for the F-conjecture. Since the exact statements and use notations from \cref{sec:prelim}, we provide here only abbreviated versions of the theorems. For the full statements, we refer to \cref{thm:main2det} and \cref{thm:main3det}.

\begin{thm}\label{thm:main2}Assume that the characteristic of the base field is not equal to $2$. Then:
    \begin{enumerate}
        \item Type 1 and type 4 F-curves on $\M{g}{n}$ span regular extremal rays of $\NE{\M{g}{n}}$, whereas type 2 F-curves do not.
        \item Apart from three exceptional families, each of which spans a regular extremal ray, no type 3 F-curve spans an extremal ray of $\NE{\M{g}{n}}$.
        \item Knudsen-type F-curves (cf.~\cref{defn:Knudsen}) are regular extremal. In particular, every F-curve on $\M{0}{n}$ and $\M{1}{n}$ is regular extremal. 
        \item There exist three additional families of type 6 F-curves, each of which spans a regular extremal ray.
    \end{enumerate}
    In characteristic $2$, the statements concerning Type~3 curves remain valid. Moreover, if \cite[Theorem~0.3]{GKM02} holds in characteristic~2, then all of the statements follow.
    \end{thm}

In the same vein, our method can be applied to study the regular extremality of F-curves for small genus.

\begin{thm}\label{thm:main3}
    Assume that the characteristic of the base field is not equal to $2$. Then:
    \begin{enumerate}
        \item Every type 5 or type 6 F-curve on $\M{2}{n}$ is regular extremal.
        \item Apart from one (resp. two) exceptional family, every type 5 and type 6 F-curve on $\M{3}{n}$ (resp. $\M{4}{n}$) is regular extremal.
    \end{enumerate}
    Moreover, if \cite[Theorem~0.3]{GKM02} holds in characteristic~2, then all of the statements follow.
\end{thm}

The dependence on the assumption that the characteristic is not equal to $2$ arises from the fact that \cite{GKM02}, and consequently many results on the F-conjecture, rely on this assumption. As stated above, if one can establish the relevant results on the F-conjecture in characteristic $2$, then our results can be extended to characteristic $2$ as well. Moreover, by \cref{cor:ext01}, every F-curve on $\M{0}{n}$ and $\M{1}{n}$ is regular extremal, regardless of the characteristic.

As suggested by the statements, F-curves of type 5 and 6 are more subtle than the others. For type 5 curves, we have an explicit conjecture (\cref{conj:exttype5}) describing which of them are regular extremal. In contrast, for type 6 curves, no conjectural description is currently available. We refer to the beginning of \cref{subsec:56} for a detailed discussion.

These theorems can be understood in the context of the general problem of understanding extremal effective cycles. The study of extremal effective cycles on $\M{g}{n}$ is a well-developed area, with many important results as in \cite{Ver02, CC14, CC15, Sch15, Sco20, Bla22}. We refer to \cite{Che18} for a survey of these results. However, many findings focus on relatively low-codimension cases. In contrast, this paper investigates the situation for extremal curves. In this setting, since many boundary strata are known to be non-extremal, distinguishing extremal from non-extremal curves is already a nontrivial task.

\cref{thm:main2} and \cref{thm:main3} require a construction of sufficiently many nef divisors on $\M{g}{n}$ that contract a fixed curve and span a codimension one subspace. In some simpler cases, such as $g=0$ or $g=1$ (see \cref{sec:case}), $\psi$-classes and their pullbacks do the job. However, this approach fails for higher genus. One needs other nef divisors, and even with enough nef divisors, computing the dimension of their span becomes challenging.

The two main advances of this paper are the construction of new nef divisors on $\M{g}{n}$ (see \cref{sec:nef}), and the development of a new induction scheme to verify the dimension of the span of these divisors (see \cref{subsec:method}).

The new nef divisors introduced in \cref{sec:nef}, referred to here as \textbf{semigroup kappa divisors}, are certain sums of the $\kappa$ class and boundary classes (cf. \cref{defn:ensemble}). By examining their intersection with F-curves \cite[Theorem 2.1]{GKM02} and boundary restriction, we prove that they are nef in \cref{thm:kapnef}. These divisors are particularly useful for contracting certain F-curves, as the intersection number of the kappa class with any F-curve is $1$.

To verify that the set of nef divisors we construct spans a codimension-\(1\) subspace, we employ a new two-step induction argument based on Knudsen's construction of \( \M{g}{n} \), which we describe in detail in \cref{subsec:method}. In \cite{Cho24, Cho25}, the author observed a particularly tractable inductive structure on the Picard group of $\M{g}{n}$ for $g\le 1$, using Knudsen’s construction. In this paper, we extend this observation to arbitrary genus in \cref{thm:Knuses}. A crucial point is that the relative cone arising from Knudsen’s construction is a simplicial cone generated by F-curves. This enables explicit computation of the dimension of the span of certain nef divisors.

Finally, we note that \cref{thm:main2} and \cref{thm:main3} hold in almost arbitrary characteristic, whereas \cref{thm:main1} and \cref{cor:main1cor} hold only over fields of characteristic zero. Moreover, as shown in \cite{Ke99} and also in \cref{cor:main1cor}, there exist divisors that are semiample, as well as morphisms that occur only over fields of positive characteristic. This demonstrates that the geometry of $\M{g}{n}$ depends strongly on the characteristic of the base field, and raises the question of providing a modular description of such morphisms in positive characteristic (cf.~\cref{subsec:morph}). 

One may compare this with the complete subvariety problem for the moduli space of abelian varieties $A_g$ \cite{KS03, GMMT25}, where the maximal-dimensional complete subvariety of $A_g$ is described in terms of the $p$-rank, a structure that exists only in positive characteristic. Just as this problem highlights the role of the $p$-rank, we expect that the answer to the question raised in \cref{subsec:morph} will likewise shed further light on the essential differences between $\M{g}{n}$ in positive and in zero characteristic. In particular, since \cref{thm:M3ncont} provides an example of a contraction of \( \M{3}{n} \) with a remarkably simple description that exists only in positive characteristic, it is an interesting problem to find a modular interpretation of this morphism. Such a description could potentially reveal new structural features of curves in positive characteristic.

\subsection{Structure of the paper}

In \cref{sec:prelim}, we review the basics of divisors and curves on $\M{g}{n}$. In \cref{sec:nonsemi}, we present the non-semiampleness criterion (\cref{thm:semin}). In \cref{sec:psisemi}, we explain when $\psi$-classes are semiample, in particular proving that they are semiample on $\M{1}{n}$. Several examples of semiample and non-semiample nef line bundles are given in \cref{sec:examples}. The remainder of the paper is devoted to the extremality of F-curves. In \cref{sec:case} we provide an explicit solution to this problem for $\M{0}{n}$ and $\M{1}{n}$. \cref{sec:nef}  and \cref{sec:Knudsen} are devoted to describe two main tools for proving extremality: a new family of nef line bundles (\cref{thm:kapnef}) and an exact sequence (\cref{thm:Knuses}). In \cref{sec:extremality}, we prove the regular extremality of certain F-curves using these tools. Further discussion and open problems are presented in \cref{sec:discussion}.

\subsection{Notations and Conventions}

Throughout, any Picard group $\text{Pic}(X)$, closed cone of curves $\NE{X}$, the cone of nef divisors $\text{Nef}(X)$, Chow group $A_d(X)$, curve class $[C]$, and divisors/line bundles will be considered over $\Q$. Thus, unless otherwise stated, these terms will refer to their $\Q$-versions, such as the $\Q$-Picard group, $\Q$-divisors, and so on. For a variety $X$, $\rho(X)$ denotes the Picard number of $X$. We denote by $[n]$ the set $\{1, 2, \dotsc, n\}$, and by $\pi_S$ the projection map $\M{g}{n} \to \M{g}{S^c}$ that forgets the marked points indexed by $S \subseteq [n]$. For any cone \( N \subseteq \mathbb{R}^n \), we define \( N \otimes \mathbb{R} \) to be the subspace of \( \mathbb{R}^n \) spanned by \( N \).

\section*{Acknowledgement}

The author would like to thank Angela Gibney for her continued support and for many helpful discussions. We are also grateful to Andreas Fernandez Herrero for valuable conversations related to contracting curves on the moduli space of curves.

%% file: sections/prelim.tex
In this section, we review background on curves and divisors on $\M{g}{n}$. Let $\partial^r\M{g}{n}$ denote the locus of stable curves with at least $r$ nodes. Then $\partial^r\M{g}{n}$ is a pure codimension $r$ subvariety of $\M{g}{n}$, and its irreducible components are called \textbf{codimension }$r$\textbf{ boundary strata}. In particular, if $r = 3g - 4 + n$, then $\partial^r\M{g}{n}$ is pure of dimension $1$, and its irreducible components are called \textbf{F-curves}. A divisor on $\M{g}{n}$ is \textbf{F-nef} if it intersects all F-curves non-negatively.

\begin{conj}\label{conj:F}[F-conjecture, \cite{KM13, GKM02}]
F-curves generate $\NE{\M{g}{n}}$. Equivalently, a divisor on $\M{g}{n}$ is nef if and only if it is F-nef.
\end{conj}

By a sequence of papers, including \cite{KM13, GKM02, GF03, Gib09, Lar11, Fe15, MS19, Fe20, FM25}, the F-conjecture is known to hold for $g + n \le 8$ or for $g \le 44$ with $n = 0$, when the characteristic of the base field is not equal to $2$. In characteristic $2$, the methods of \cite{Lar11, FM25} (resp. \cite{Fab90}) prove the conjecture for $g = 0$, $n \le 8$ (resp. $n = 0$, $g = 2, 3$). For later use, in \cref{prop:charpsemi}, we verify the F-conjecture in arbitrary characteristic for $\M{1}{2}$, $\M{2}{1}$, $\M{2}{2}$, and $\M{3}{1}$.

One significant advantage of \cref{conj:F} is that, as we will explain in \cref{thm:explicit}, the cone of F-nef divisors admits a very explicit description. To this end, we will provide a more detailed description of F-curves and divisors on $\M{g}{n}$, following \cite{GKM02}.

From now on, we will identify F-curves with their classes in $\text{A}_1(\M{g}{n})$, i.e., up to numerical equivalence. There are six types of F-curves (cf. \cite[Theorem 2.2]{GKM02}), which are described as follows.

\begin{description}
    \item[Type 1] Let $i:\M{1}{1} \to \M{g}{n}$ be the map that attaches a fixed semistable curve $C$ of genus $g - 1$ with $n + 1$ marked points and then stabilizes. The image $\text{F}_1$ of $i$ is the \textbf{F-curve of type 1}.
    
    \item[Type 2] Let $i:\M{0}{4} \to \M{g}{n}$ be the map that attaches a fixed semistable curve $C$ of genus $g - 3$ with $n + 4$ marked points to the four marked points on curves parametrized by $\M{0}{4}$, and then stabilizes. The image $F_2$ of $i$ is the \textbf{F-curve of type 2}.
    
    \item[Type 3] Choose natural numbers $g_1 + g_2 = g - 2$ and a decomposition $I_1 \sqcup I_2 = [n]$. Fix semistable curves $C_1$ (resp. $C_2$) of genus $g_1$ (resp. $g_2$) with $|I_1| + 1$ (resp. $|I_2| + 3$) marked points. Let $i:\M{0}{4} \to \M{g}{n}$ be the map that attaches $C_1$ to the first point and $C_2$ to the remaining three points, and then stabilizes. The image $F_3^{g_1}(I_1)$ of $i$ is an \textbf{F-curve of type 3}. Note that $g_2$ and $I_2$ are determined by $g_1$ and $I_1$.
    
    \item[Type 4] Choose natural numbers $g_1 + g_2 = g - 2$ and a decomposition $I_1 \sqcup I_2 = [n]$. Fix semistable curves $C_i$ of genus $g_i$ with $|I_i| + 2$ marked points for $i = 1, 2$. Let $i:\M{0}{4} \to \M{g}{n}$ be the map that attaches $C_1$ and $C_2$ to two of the four points and then stabilizes. The image $F_4^{g_1}(I_1)$ of $i$ is an \textbf{F-curve of type 4}.
    
    \item[Type 5] Choose natural numbers $g_1 + g_2 + g_3 = g - 1$ and a decomposition $I_1 \sqcup I_2 \sqcup I_3 = [n]$. Fix semistable curves $C_i$ of genus $g_i$ with $|I_i| + 1$ marked points for $i = 1, 2$, and $C_3$ of genus $g_3$ with $|I_3| + 2$ marked points. Let $i:\M{0}{4} \to \M{g}{n}$ be the map that attaches $C_1$ and $C_2$ to one point each, and $C_3$ to the remaining two points, and then stabilizes. The image $F_5^{g_1, g_2}(I_1, I_2)$ of $i$ is an \textbf{F-curve of type 5}.
    
    \item[Type 6] Choose natural numbers $g_1 + g_2 + g_3 + g_4 = g - 1$ and a decomposition $I_1 \sqcup I_2 \sqcup I_3 \sqcup I_4 = [n]$. Fix semistable curves $C_i$ of genus $g_i$ with $|I_i| + 1$ marked points. Let $i:\M{0}{4} \to \M{g}{n}$ be the map that attaches each $C_i$ to one of the four points and then stabilizes. The image $F_6^{g_1, g_2, g_3, g_4}(I_1, I_2, I_3, I_4)$ of $i$ is an \textbf{F-curve of type 6}.
\end{description}

While the actual curve on $\M{g}{n}$ depends on the choice of semistable curves, its numerical class remains independent. For a pictorial description of these curves, we refer to \cite[Figure (2.3)]{GKM02}. Note that if $g = 0$, then only type 6 curves exist; if $g = 1$, then types 1, 5, and 6 exist; if $g = 2$, then all types except type 2 exist; and if $g \ge 3$, then F-curves of all types exist.

Note that the $\lambda$-class $\lambda$, the $\psi$-classes $\psi_i$, and the boundary classes $\delta_{i,I}$ generate the $\Q$-Picard group of $\M{g}{n}$. Furthermore, if $g \ge 3$, they form a basis of the $\Q$-Picard group. For relations among these classes, we refer to \cite[Theorem 4]{AC09}. Following \cite{GKM02}, we express a divisor $D$ as a linear combination:
\[
D = a\lambda - b_{\text{irr}} \delta_{\text{irr}} + \sum_{i=1}^n b_{0,i} \psi_i - \sum b_{i, I} \delta_{i, I}.
\]

\begin{thm}\label{thm:explicit}\cite[Theorem 2.1]{GKM02}    
The intersections of F-curves with the divisor $D$ are as follows:
\begin{enumerate}
    \item $D \cdot F_1 = a - 12b_{\text{irr}} + b_{1, \emptyset}$.
    \item $D \cdot F_2 = b_{\text{irr}}$.
    \item $D \cdot F_3^{i}(I) = b_{i, I}$.
    \item $D \cdot F_4^{i}(I) = 2b_{\text{irr}} - b_{i, I}$.
    \item $D \cdot F_5^{i,j}(I,J) = b_{i, I} + b_{j, J} - b_{i+j, I \cup J}$.
    \item $D \cdot F_6^{i,j,k,l}(I,J,K,L) = b_{i, I} + b_{j, J} + b_{k, K} + b_{l, L} - b_{i+j, I \cup J} - b_{i+k, I \cup K} - b_{i+l, I \cup L}$.
\end{enumerate}

\end{thm}

%% file: sections/nonsemi.tex
\begin{thm}\label{thm:semin} (\text{char }$k=0$)
    Let $g, n\ge 2$ and $L$ be a line bundle on $\M{g}{n}$ such that
    \[ L=a\lambda -b_{\text{irr}}\delta_{\text{irr}}-\sum_{i,I}b_{i,I}\delta_{i,I}. \]
    Assume that there exist $i,j\in [n]$ such that 
    \begin{enumerate}
        \item $b_{0,\left\{i,j\right\}}=0$,
        \item $b_{0,i}\ne b_{0,j}$, and
        \item $b_{0,k}=0$ for every $k\ne i,j$.
    \end{enumerate}
    Then $L$ is not semiample.
\end{thm}

The proof is almost identical to that of \cite[Theorem 4.6]{Cho23}. However, we will reproduce the proof here since a more specific circumstance is considered and a different language is used there. The proof is based on Keel's counterexample \cite[Section 3]{Ke92}.

\begin{proof}
    First, we will consider the case of $n=2$. In this case, condition (3) is vacuous. Since $\text{Pic}(\overline{\mathcal{M}}_{g,n})\simeq \text{Pic}(\M{g}{n})$, we will consider $L$ as a line bundle on $\overline{\mathcal{M}}_{g,n}$. Let $C$ be a projective smooth curve of genus $g$. Then, the projection $\pi_1:C\times C\to C$ with the diagonal embedding $\Delta:C\to C\times C$ can be considered as a family of genus $g$ curves with a marked point. Hence, there exists a corresponding morphism
    \[
    \begin{tikzcd}[ampersand replacement=\&]
    C\times C \arrow[r, "u"] \arrow[d, "\pi_1"] \& \overline{\mathcal{M}}_{g,2} \arrow[d, "\pi_1"]\\
    C \arrow[r] \& \overline{\mathcal{M}}_{g,1}
    \end{tikzcd}
    \]
    Consider $u^{\ast}L$. Note that the image of $u$ does not intersect any boundary divisors except $\delta_{0,\left\{1,2\right\}}$. Since $b_{0,\left\{1,2\right\}}=0$, the boundary divisors do not contribute to $u^\ast L$. Moreover, since the image of the composition map $C\times C\to \overline{\mathcal{M}}_{g,2}\to \overline{\mathcal{M}}_{g}$ is a point corresponding to $C$, we have $u^\ast \lambda=0$. Hence, $u^\ast L = b_{0,1} u^\ast \psi_1 + b_{0,2} u^\ast \psi_2$. Since the image of $u \circ \Delta$ is contained in $\Delta_{0,\left\{1,2\right\}}$, we have $\Delta^\ast u^\ast L = 0$.

    Now assume that $L$ is a semiample line bundle. Then $u^\ast L$ is also a semiample line bundle. Let $f: C\times C \to \mathbb{P}^n$ be a morphism such that $f^\ast \mathcal{O}(1) = u^\ast L^{\otimes m}$. Since $\Delta^\ast u^\ast L = 0$, the map $f$ contracts the image of $\Delta$. Hence, there exists an effective divisor $T$ on $C\times C$ such that $\Delta \cap T = \emptyset$ and $u^\ast L^{\otimes m} = [T]$. In particular, the restriction of $u^\ast L^{\otimes m}$ to any infinitesimal neighborhood of $\Delta$ is trivial. Let $\Delta_2$ be the second-order infinitesimal neighborhood of $\Delta$.

    By the definition of $\psi$ classes (cf. \cite[Proposition 3.9]{Cho23}), we have $u^\ast \psi_i = \pi_i^\ast\Omega^1_{C}(\Delta)$ for $i = 1,2$. By \cite[Lemma 3.5]{Ke99}, we have $u^\ast \psi_1|_{\Delta_2} = -u^\ast \psi_2|_{\Delta_2}$. Hence,
    \[
    u^\ast L^{\otimes m}|_{\Delta_2} = m(b_{0,1} - b_{0,2})u^\ast \psi_1|_{\Delta_2} = m(b_{0,1} - b_{0,2})\pi_1^\ast\Omega^1_{C}(\Delta)|_{\Delta_2}.
    \]
    Since $b_{0,1} \ne b_{0,2}$, this is not a trivial line bundle by \cite[Lemma 3.4]{Ke99}. This leads to a contradiction. Hence, $L$ is not a semiample line bundle.

    Now consider the general case. We may assume that $i=1$ and $j=n$. Then there exist a projective smooth curve $C$ of genus $g$, a projective smooth curve $D$, and morphisms $s_1,\cdots, s_{n-1}:D\to C$ such that the trivial family of curves $\pi_1:D\times C\to D$ with sections $s_1,\cdots, s_{n-1}$ forms a nontrivial family of genus $g$ curves with $n-1$ marked points. The construction uses elliptic curves. See \cite{Za05} or \cite[proof of Corollary 4.8]{Cho23}.

    As in the case of $n=2$, we have an induced morphism $u:D\times C\to \overline{\mathcal{M}}_{g,n}$,
    \begin{align*}
         u^\ast L&= b_{0,1} u^\ast \psi_1 + b_{0,n} u^\ast \psi_n-\sum_{i=2}^{n-1}b_{0,\left\{i,n\right\}}u^\ast\delta_{0,\left\{i,n\right\}}\\& = b_{0,1} \pi_1^\ast s_1^\ast\Omega^1_{C}(s_1) + b_{0,n} \pi_2^\ast\Omega^1_{C}(s_1)-\sum_{i=2}^{n-1}b_{0,\left\{i,n\right\}}[s_i]
    \end{align*}
    where we identify $s_i:D\to D\times C$ with its image, and $s_{1}^\ast u^\ast L=0$. Assume that $L$ is a semiample line bundle. Then, by the same argument as above, there exists an effective divisor $T$ on $C\times D$ such that $u^\ast L^{\otimes m}=[T]$ for some $m>0$ and $T\cap s_1=\emptyset$.

    Consider $f=(s_1, \text{id}):D\times C\to C\times C$. Let $R$ be the effective divisor $T+\sum_{i=2}^n [s_i]$, and define $R':=f_\ast R$, $U:=C\times C\setminus R'$, and $V:=f^{-1}(U)$. Since $f^{-1}(\Delta)=s_1$ and $s_1\cap R=\emptyset$, we have $\Delta \subseteq U$. Since $R\cap V=\emptyset$, we obtain
    \[
    u^\ast L|_V = b_{0,1} \pi_1^\ast s_1^\ast\Omega^1_{C}(s_1)|_V + b_{0,n} \pi_2^\ast\Omega^1_{C}(s_1)|_V=(f|_V)^\ast(b_{0,1}\pi_1^\ast\Omega^1_{C}(\Delta)+b_{0,n}\pi_2^\ast\Omega^1_{C}(\Delta)).
    \]
    On the other hand, since $u^\ast L^{\otimes m}=[T]$ and $T\cap V=\emptyset$, we conclude that $u^\ast L^{\otimes m}|_V$ is trivial. By a standard push-pull argument (e.g., \cite[Lemma 4.7]{Cho23}), there exists $d>0$ such that
    \[
    db_{0,1}\pi_1^\ast\Omega^1_{C}(\Delta)+db_{0,n}\pi_2^\ast\Omega^1_{C}(\Delta)=0 \quad \text{on } U.
    \]
    Since $\Delta\subseteq U$, its second-order infinitesimal neighborhood $\Delta_2$ is also contained in $U$, so its restriction to $\Delta_2$ is also zero. However, by \cite[Lemma 3.5]{Ke92}, we have
    \[
    (db_{0,1}\pi_1^\ast\Omega^1_{C}(\Delta)+db_{0,n}\pi_2^\ast\Omega^1_{C}(\Delta))|_{\Delta_2}=d(b_{0,1}-b_{0,n})\pi_1^\ast\Omega^1_{C}(\Delta)|_{\Delta_2},
    \]
    which is nonzero by \cite[Lemma 3.4]{Ke92} and condition (2). This leads to a contradiction. Hence, $L$ is not semiample.

\end{proof}

The condition $b_{0,i}\ne b_{0,j}$ is essential: see \cref{prop:lambdapsi}.

%% file: sections/psisemi.tex
\begin{thm}\label{thm:psisemi}
    $\psi_i$'s are semiample on $\M{1}{n}$.
\end{thm}

\begin{proof}
    We will use induction on $n$. If $n=1$, then this is trivial. If $n=2$, there are several ways to prove this. There exists a surjective map $f:\M{0}{5}\to \M{1}{2}$, which sends $(C, \left\{P_i\right\}_{i=1}^{5})$ to $(D, Q_1, Q_2)$, where $D$ is a double cover of $C$ ramified at $P_2, P_3, P_4, P_5$, $Q_2$ is the inverse image of $P_2$, and $Q_1$ is an inverse image of $P_1$. See, for example, \cite[Section 4.1]{Rul01}. Since $\M{1}{2}$ is normal, a line bundle $L$ on $\M{1}{2}$ is semiample if and only if $f^\ast L$ is. Since $\M{0}{5}$ is a Fano variety, any nef line bundle is semiample. Therefore, any nef line bundle on $\M{1}{2}$ is semiample. Alternatively, by \cite[Proposition 4.1.5]{Rul01}, for the natural map $\alpha:\M{1}{2}\to \overline{\rm{M}}_{2}$ corresponding to $\Delta_{\text{irr}}$, $\psi_1=\alpha^\ast (12\lambda-\delta_{\text{irr}})$. By \cite[Corollary 4.3]{Gib12}, $12\lambda-\delta_{\text{irr}}$ is semiample, hence $\psi_1$ is also semiample.
    
    Now assume $\psi_i$'s are semiample on $\M{1}{n-1}$ and \(n\ge 3\) . It is enough to prove that $\psi_1$ is semiample on $\M{1}{n}$. Let $\pi_{n-1}$ and $\pi_n:\M{1}{n}\to\M{1}{n-1}$ be the projection maps forgetting the $(n-1)$th and $n$th points, respectively. Then
    \[ \psi_1=\pi_{n-1}^\ast \psi_1+\delta_{0,\left\{1, n-1\right\}}=\pi_{n}^\ast \psi_1+\delta_{0,\left\{1, n\right\}}. \]
    By the induction hypothesis, $\pi_{n-1}^\ast \psi_1$ and $\pi_{n}^\ast \psi_1$ are semiample. Hence, the stable base locus $\textbf{B}(\psi_1)$ is contained in $\Delta_{0,\left\{1, n-1\right\}}$ and $\Delta_{0,\left\{1, n\right\}}$. However, their intersection is trivial, hence $\textbf{B}(\psi_1)=\emptyset$. Therefore, $\psi_1$ is semiample.

\end{proof}

\begin{cor}\label{cor:psisem}
    The $\psi$-classes are semiample on $\M{g}{n}$ if and only if one of the following holds: 
    \begin{enumerate} 
    \item The characteristic of the base field is positive. 
    \item $g \leq 1$. 
    \item $(g,n) = (2,1)$. 
    \end{enumerate}
\end{cor}

\begin{proof}
    Case (1) is proved in \cite{Ke99}. The $g=0$ case of (2) follows from the existence of Kapranov's construction \cite{Ka92, Ka93}. We have just proved the $g=1$ case. (3) follows from \cite{Rul01}, where Rulla proved that any nef line bundle on $\M{2}{1}$ is semiample. Now, consider the characteristic zero case. If $g \geq 2$ and $n \geq 2$, then \cref{thm:semin} (or the original argument in \cite{Ke99}) implies that the $\psi_i$'s are not semiample. For $g \geq 3$, let $f:\M{g-1}{3}\to \M{g}{1}$ be the map corresponding to $\Delta_{\text{irr}}$. Then $f^\ast \psi_1 = \psi_1$, and this is not semiample. Hence, $\psi_1$ is not semiample if $g \geq 3$.
\end{proof}

\begin{qes}\label{qes:psisemi}
    What is the contraction of $\M{1}{n}$ corresponds to $\psi_1$?
\end{qes}

Regarding the proof of \cref{thm:psisemi}, we should figure out the contraction of $\overline{\rm{M}}_{2}$ corresponds to $12\lambda-\delta_{\text{irr}}$.

%% file: sections/examples.tex
%\begin{eg}\label{eg:eg22}[char $k=0$]
%    If $L$ is a line bundle on $\M{2}{2}$ of the form
%    \[ \psi_1+\left(\text{linear combination of } \lambda,  \delta_{\text{irr}},\delta_{1, \emptyset}, \delta_{1, \left\{1\right\}}  \right) \]
%    then $L$ is not semiample. More precisely,  
%    \[ L=\psi_1-b_{\text{irr}}\delta_{\text{irr}}-b_{1, \emptyset}\delta_{1, \emptyset}-b_{1, \left\{1\right\}}\delta_{1, \left\{1\right\}} \]
%    is nef and not semiample if and only if
%    \begin{enumerate}
%        \item $b_{1, \emptyset}\ge 12 b_{\text{irr}}$,
%        \item $2 b_{\text{irr}}\ge b_{1, \emptyset}, b_{1, \left\{1\right\}}$,
%        \item $2 b_{1, \emptyset}+1\ge 2 b_{1, \left\{1\right\}} \ge 2 b_{1, \emptyset}$.
%    \end{enumerate}
%\end{eg}

%\begin{eg}
%    $\psi_1+\delta_{\text{irr}}+2\delta_1$ is a nef and non-semiample line bundle on $\M{2}{2}$.
%\end{eg}

%We can bootstrap this to construct many other semiample line bundles.

%\begin{prop}
%     If $L$ is a line bundle on $\M{3}{2}$ of the form
%    \[ \psi_1+\left(\text{linear combination of } \lambda, \psi_2, \delta_{\text{irr}}, \delta_{0, \left\{1,2\right\}}, \delta_{1, \emptyset}, \delta_{1, \left\{1\right\}}  \right) \]
%    then $L$ is not semiample. 
%\end{prop}

%\begin{proof}
%    Let $f:\M{2}{2}\to \M{3}{2}$ be the map attaching a fixed stable genus $1$ curve with $2$ marked points. Then $f^\ast L$ is of the form \cref{eg:eg22}.
%\end{proof}

%If $n=2$, then there is a more general theorem:

In this section, we present examples of nef divisors that are not semiample, as given by \cref{thm:semin}. Throughout this section, we assume that the characteristic of the base field is zero.

\begin{eg}\label{eg:M22}
    On $\M{2}{2}$, the following divisors are the set of extremal rays of the nef cone:
    \[\begin{gathered}
        \lambda,\ 12\lambda - \delta_{\text{irr}},\ \psi_1,\ \psi_2,\ \pi_1^\ast \psi,\ \pi_2^\ast \psi,\\
        \psi_1+\psi_2+\delta_{1,\left\{1\right\}},  \psi_1+\psi_2-2\delta_{0,\left\{1,2\right\}}+\delta_{1,\left\{1\right\}}\\
        \delta_{\text{irr}}+10\psi_1+10\psi_2+12\delta_{1, \emptyset}+2\delta_{1,\left\{1\right\}}=10(\lambda+\psi_1+\psi_2+\delta_{1, \emptyset})
    \end{gathered}\]
    Note that on $\M{2}{n}$, $10\lambda=\delta_{\text{irr}}+2\delta_1$.
\end{eg}

\begin{thm}\label{thm:cone22}
    There exists a $5$-dimensional subcone of the nef cone of $\M{2}{2}$ whose generic element is not semiample.
\end{thm}

\begin{proof}
    Consider the following divisors:
    \[ \lambda, 12\lambda - \delta_{\text{irr}}, \psi_1, \psi_2,  \psi_1+\psi_2+\delta_{1,\left\{1\right\}}, \lambda+\psi_1+\psi_2+\delta_{1, \emptyset}.  \]
    A generic element of the cone generated by these elements satisfies the condition of \cref{thm:semin}, and hence not semiample. It is straightforward to see that these elements generate a $5$-dimensional space, spanned by $\psi_1, \psi_2, \delta_{\text{irr}}, \delta_{1, \emptyset}, \delta_{1,\left\{1\right\}}$.
\end{proof}

Since the Picard number of $\M{2}{2}$ is $6$, this represents the largest possible such cone. Hence, this provides an instance of the cone described in \cref{thm:main1} within $\M{2}{2}$.

\begin{prop}\label{prop:semi31}
     If $D$ is a divisor on $\M{3}{1}$ of the form
     \[ a\lambda - b_{\text{irr}}\delta_{\text{irr}} + b_{0,1}\psi_1 - b_{1,\emptyset}\delta_{1,\emptyset} \]
     for some $b_{0,1} \ne b_{1,\emptyset}$, then $D$ is not semiample. 
\end{prop}

\begin{proof}
    Let $f: \M{2}{2} \to \M{3}{1}$ be the morphism obtained by attaching a genus one curve at a marked point. By assumption, the pullback $f^\ast D$ satisfies the condition of \cref{thm:semin}. Therefore, $f^\ast D$ is not semiample, and it follows that $D$ is also not semiample.
\end{proof}

Using \cref{prop:semi31}, we can directly verify \cref{thm:main1} and \cref{cor:main1cor} for $\M{3}{1}$.

\begin{eg}\label{eg:counter}
    There exists no contraction of $\M{3}{1}$ whose relative cone of curves is the extremal ray spanned by $F_3^1(1)$. Suppose, for contradiction, that such a contraction exists and let $D$ be the corresponding semiample divisor on $\M{3}{1}$. Since $D$ contracts $F_3^1(1)$, it must be of the form described in \cref{prop:semi31}. Then, by \cref{thm:semin}, it must satisfy $b_{0,1} = b_{1,\emptyset}$. Under the additional assumption that $D \cdot F_3^1(1) = 0$, this is precisely equivalent to $D \cdot F_5^{2,0}(\emptyset, 1) = 0$. Thus, $D$ must also contract $F_5^{2,0}(\emptyset, 1)$.
\end{eg}

\begin{eg}\label{eg:M31}
On $\M{3}{1}$, the following divisors are the set of extremal rays of the F-nef cone:
\[ \psi_1,\ \lambda,\ 12\lambda-\delta_{\text{irr}},\  10\lambda-\delta_{\text{irr}}-2\delta_{1,\emptyset}-2\delta_{1, \left\{1\right\}},\ 10\lambda-\delta_{\text{irr}}+2\psi_1-2\delta_{1,\emptyset}, 11\lambda-\delta_{\text{irr}}+3\psi_1-\delta_{1,\emptyset}-2\delta_{1, \left\{1\right\}}.  \]
They generates the nef cone by \cref{prop:charpsemi}.
\end{eg}

\begin{thm}\label{thm:cone31}
    There exists a $4$-dimensional subcone of the nef cone of $\M{3}{1}$ whose generic element is not semiample.
\end{thm}

\begin{proof}
    Consider the following divisors:
    \[\psi_1, \lambda, 12\lambda-\delta_{\text{irr}},  10\lambda-\delta_{\text{irr}}+2\psi_1-2\delta_{1,\emptyset}.  \]
    A generic element of the cone generated by these elements satisfies the condition of \cref{prop:semi31}, and hence not semiample. It is straightforward to see that these elements generate a $4$-dimensional space, spanned by $\lambda, \psi_1, \delta_{\text{irr}}, \delta_{1, \emptyset}$.
\end{proof}

Again, this is the largest possible such cone, since the Picard number of $\M{3}{1}$ is $5$. This naturally leads to the following question, which is partially answered by \cref{thm:main1}.

\begin{qes} \label{qes:nonsemi}
Let $\rho$ denote the Picard number of $\M{g}{n}$, with $g \ge 3$ and $n \ge 1$, or $g = 2$ and $n \ge 2$. Does there exist a $(\rho - 1)$-dimensional subcone of the nef cone, whose general element is not semiample?
\end{qes}

By the following observation, this is related to the extremality of F-curves.

\begin{prop}\label{prop:extcone} 
Let $C$ be a pointed $n$-dimensional polyhedral cone in $\mathbb{R}^n$, and let $C^\ast$ be its dual. Then $v \in C$ spans an extremal ray if and only if the dimension of $C^\ast \cap v^\perp$ is exactly $n-1$. 
\end{prop}

One important caveat is that, since it is not known whether the cone of curves on $\M{g}{n}$ is polyhedral, extremality alone is generally insufficient to draw definitive conclusions. To address this, we introduce the notion of a regular extremal ray in \cref{defn:convex}.

From the perspective of \cref{prop:extcone}, \cref{thm:cone22} corresponds to the fact that $F_3^0([2])$ spans an extremal ray in $\M{2}{2}$, while \cref{thm:cone31} reflects the extremality of $F_3^1(1)$ in $\M{3}{1}$. The connection between non-semiampleness and F-curves arises from \cref{thm:semin}: assumptions (1) and (3) in that theorem are equivalent to the condition that the divisor contracts certain F-curves. In particular, assumption (1) corresponds to the condition $L \cdot F_3^{0}(\{i, j\}) = 0$, and (3) is an empty condition when $n = 2$. This naturally raises the question of whether these F-curves are extremal. Even when $n > 2$, we follow a similar strategy as in \cref{prop:semi31}, by pulling back the divisor to the case $n = 2$. Therefore, in what follows, we focus on analyzing the extremality of F-curves. Indeed, in \cite{Bla22}, the author raised the following conjecture:

\begin{conj}\cite[Conjecture 1.3]{Bla22}
    All boundary strata in $\M{g}{n}$ are extremal.
\end{conj}

However, this conjecture is not true for F-curves, i.e. one-dimensional boundary strata. The easiest example is the type $2$ F-curve. If $g\ge 3$ and $L$ is a nef divisor such that $L\cdot F_2=0$, i.e., $b_{\text{irr}}=0$, then by taking intersections with F-curves of type 3 and type 4, we obtain $b_{i,I}=0$ for every $1\le i\le g-2$. Therefore, by \cref{prop:extcone}, $F_2$ cannot be extremal. Explicitly, 
\[ [F_2]=\frac{1}{2}\left([F_3(i,I)]+[F_4(i,I)] \right). \]

Moreover, by \cref{thm:type3next}, all type 3 curves except those listed in \cref{thm:main2} are not extremal. This provides further examples of non-extremal F-curves.

As an illustrative example, we examine the extremality of F-curves in $\overline{\mathrm{M}}_4$, $\M{3}{1}$, and $\M{2}{2}$ using the computer program \cite[extray]{Choigit25}. This program determines which F-curves are extremal in the cone of F-curves. By the discussion in \cref{sec:prelim}, we know that for these spaces, this cone coincides with the cone of curves, so extremality in the cone of F-curves is equivalent to extremality in the cone of curves. Note that, for the sake of time complexity, we exclude type 2 and type 3 F-curves not listed in \cref{thm:main2} from the program, since they are already known to be non-extremal.

\begin{eg}\label{eg:mext}
    On $\overline{\rm{M}}_4$, the following divisors are the set of extremal rays of the nef cone:
\[ \lambda,\ 12\lambda-\delta_{\text{irr}},\  10\lambda-\delta_{\text{irr}}-2\delta_{1},\ 10\lambda-\delta_{\text{irr}}-2\delta_{1}-2\delta_2,\ 21\lambda-2\delta_{\text{irr}}-3\delta_{1}-4\delta_{2}   \]
Note that $21\lambda-2\delta_{\text{irr}}-3\delta_{1}-4\delta_{2}$ contracts the type 6 F-curve.

Here is a summary of the extremality of F-curves:
\begin{table}[htbp]
\centering
\begin{tabular}{ |c|c|c|c| } 
 \hline
  & extremal? & Relation \\
  \hline
 $F_1$ & Y & - \\ 
 \hline
 $F_2$ & N & $\frac{1}{2}\left([F_3^1(\emptyset)]+[F_4^1(\emptyset)] \right)$ \\ 
 \hline
 $F_3^1(\emptyset)$ & N & $\frac{1}{2}\left([F_5^{1,1}(\emptyset,\emptyset)]+[F_3^2(\emptyset)]\right)$ \\ 
 \hline
 $F_3^2(\emptyset)=F_5^{1,2}(\emptyset, \emptyset)$ & Y & - \\ 
 \hline
 $F_4^1(\emptyset)$ & Y & - \\ 
 \hline
 $F_4^2(\emptyset)$ & Y & - \\ 
 \hline
 $F_5^{1,1}(\emptyset, \emptyset)$ & N & $\frac{1}{2}\left( [F_6^{1,1,1,1}(\emptyset, \emptyset, \emptyset, \emptyset)]+[F_3^2( \emptyset)] \right)$ \\ 
 \hline
 $F_6^{1,1,1,1}(\emptyset, \emptyset, \emptyset, \emptyset)$ & Y & - \\ 
 \hline
\end{tabular}
 \caption{F-curves on $\overline{\rm{M}}_4$.}
 \end{table}
\end{eg}
\vspace{-1.5em}

We have a similar table for $\M{3}{1}$ and $\M{2}{2}$.
\begin{longtable}{ |c|c|c|c| }
 \hline
  & extremal? & Relation \\
  \hline
 $F_1$ & Y & - \\ 
 \hline
 $F_2$ & N & $\frac{1}{2}\left([F_3^1(\emptyset)]+[F_4^1(\emptyset)] \right)$ \\ 
 \hline
 $F_3^0(\left\{1\right\})$ & N & $\frac{1}{2}\left([F_5^{0,1}(\left\{1\right\},\emptyset)]+[F_5^{0,2}(\left\{1\right\},\emptyset)] \right)$ \\ 
 \hline
 $F_3^1(\emptyset)$ & N & $\frac{1}{2}\left([F_3^1(\left\{1\right\})]+[F_5^{1,1}(\emptyset,\emptyset)] \right)$ \\ 
 \hline
 $F_3^1(\left\{1\right\})=F_5^{1,1}(\left\{1\right\},\emptyset)$ & Y & - \\ 
 \hline
 $F_4^1(\emptyset)=F_4^2(\left\{1\right\})$ & Y & - \\ 
 \hline
 $F_4^1(\left\{1\right\})=F_4^2(\emptyset)$ & Y & - \\ 
 \hline
 $F_5^{0,1}(\left\{1\right\},\emptyset)$ & N & $\frac{1}{3}\left( [F_6^{1,1,1,0}(\emptyset,\emptyset,\emptyset,\left\{1\right\})]+2[F_3^0( \left\{1\right\})] \right)$ \\ 
 \hline
 $F_5^{0,2}(\left\{1\right\},\emptyset)$ & Y & - \\ 
 \hline
 $F_5^{1,1}(\emptyset,\emptyset)$ & Y & - \\ 
 \hline
 $F_6^{1,1,1,0}(\emptyset,\emptyset,\emptyset,\left\{1\right\})$ & Y & - \\ 
 \hline
  \caption{F-curves on $\M{3}{1}$.}
\end{longtable}
\vspace{-1.5em}

\begin{longtable}{ |c|c|c|c| }
 \hline
  & extremal? & Relation \\
  \hline
 $F_1$ & Y & - \\ 
 \hline
 $F_3^0(\left\{1\right\})$ & N & $\frac{1}{2}\left([F_5^{0,1}(\left\{1\right\},\emptyset)]+[F_5^{0,1}(\left\{1\right\},\left\{2\right\})] \right)$ \\ 
 \hline
 $F_3^0( \left\{2\right\})$ & N & $\frac{1}{2}\left([F_5^{0,1}(\left\{2\right\},\emptyset)]+[F_5^{0,1}(\left\{2\right\},\left\{1\right\})] \right)$ \\ 
 \hline
 $F_3^0(\left\{1,2\right\})=F_5^{0,1}(\left\{1,2\right\},\emptyset)$ & Y & - \\ 
 \hline
 $F_4^1(\emptyset)$ & Y & - \\ 
 \hline
 $F_4^1(\left\{1\right\})$ & Y & - \\ 
 \hline
 $F_5^{0,0}(\left\{1\right\},\left\{2\right\})$ & Y & - \\ 
 \hline
 $F_5^{0,1}(\left\{1\right\},\emptyset)$ & Y & - \\ 
 \hline
 $F_5^{0,1}(\left\{1\right\},\left\{2\right\})$ & Y & - \\ 
 \hline
 $F_5^{0,1}(\left\{2\right\},\emptyset)$ & Y & - \\ 
 \hline
 $F_5^{0,1}(\left\{2\right\},\left\{1\right\})$ & Y & - \\ 
 \hline
$F_6^{1,1,0,0}(\emptyset,\emptyset,\left\{1\right\},\left\{2\right\})$ & Y & - \\ 
 \hline
  \caption{F-curves on $\M{2}{2}$.}
\end{longtable}

\vspace{-1.5em}

%Even after examining the extremality of the F-curve, it remains unclear whether the face of the nef cone defined by that F-curve contains many non-semiample line bundles. This issue corresponds to the problem of characterizing line bundles on contractions.

\begin{rmk}\label{rmk:type6}
    The previous tables, together with a result \cite[Theorem 1.1]{Bla22}, might suggest the conjecture that all type 6 F-curves are extremal. However, this is an illusion. On $\overline{\rm{M}}_7$,
    \[ [F_6^{1,1,2,3}(\emptyset,\emptyset,\emptyset,\emptyset)] = \frac{1}{2} \left( [F_6^{1,1,1,4}(\emptyset,\emptyset,\emptyset,\emptyset)] + [F_6^{1,2,2,2}(\emptyset,\emptyset,\emptyset,\emptyset)] \right). \]
    This counterexample arises from the observation in \cite[Table 1]{Mo17}, which is also mentioned in \cite[Remark 6.3]{CC15} and \cite[Remark 4.6]{Bla22}. The key point is that not every extremal ray of $\M{0}{g+n}$ remains extremal in the quotient $\M{0}{g+n}/S_g$.
\end{rmk}

From now on, we will prove the semiampleness of certain line bundles on $\M{2}{2}$. The following lemma will be useful for this purpose.

\begin{lem}\label{lem:descend}
    Let $f: X \to Y$ be a contraction of projective varieties. Let $L$ be a semiample $\Q$-line bundle on $X$. If $L$ intersects trivially with $\NE{f}$, then $L$ descends to $Y$, i.e., $L \in f^\ast \mathrm{Pic}(Y)$.
\end{lem}

\begin{proof}
    Let $f': X \to Y'$ be the contraction corresponding to $L$. Then, by \cite[Proposition 1.14(b)]{Deb01}, $f'$ factors through $f$. Since $L \in f'^\ast \mathrm{Pic}(Y')$, it follows that $L \in f^\ast \mathrm{Pic}(Y)$.
\end{proof}

\begin{prop}\label{prop:lambdapsi} 
For any $\epsilon>0$, $D_\epsilon=\epsilon(12\lambda-\delta_{\text{irr}})+\psi_1+\psi_2$ is semiample on $\M{2}{2}$. 
\end{prop}

\begin{proof} 
    First, we will prove $D_1$ is semiample. By \cite[Corollary 4.3]{Gib12}, $12\lambda-\delta_{\text{irr}}$ is semiample for $g\le 11$. Since $\psi_1+\psi_2-2\delta_{0,\left\{1,2\right\}}$ and $12\lambda-\delta_{\text{irr}}$ are semiample by \cref{cor:psisem}, the stable base locus of $D_1$ is contained in $\Delta_{0,\left\{1,2\right\}}$. Moreover, for $\xi:\M{2}{2}\to \overline{\rm{M}}_{3}$,  
    \[ D_1=\xi^\ast\left(12\lambda-\delta_{\text{irr}}\right)+\delta_{1,1} \]
    Hence, the stable base locus is contained in $\Delta_{1,1}$. However, this does not intersect with $\Delta_{0,\left\{1,2\right\}}$. Hence, $D_1$ is semiample. 

    Let $f:\M{2}{2} \to X$ be the contraction corresponding to $D_1$. Among the F-curves, $D_1$ contracts only $F_1$ and $F_3^0([2])$, by \cref{thm:explicit}. Hence, $D_1$ lies in the interior of the face of the nef cone that intersects $F_1$ and $F_3^0([2])$ trivially, i.e., the cone generated by
    \[
    12\lambda - \delta_{\text{irr}} ,\ \psi_1,\ \psi_2,\ \psi_1 + \psi_2 + \delta_{1,\{1\}},\ \delta_{\text{irr}} + 10\psi_1 + 10\psi_2 + 12\delta_{1, \emptyset} + 2\delta_{1,\{1\}}.
    \]
    Let $F = f^\ast \mathrm{Nef}(X)$ be the pullback of the nef cone of $X$, which is a subcone of the described cone. Since $D_1$ is semiample, $D_\epsilon$ is also semiample for $\epsilon \ge 1$. As $D_\epsilon$ contracts only $F_1$ and $F_3^0([2])$, it intersects trivially with $\NE{f}$. Thus, by \cref{lem:descend}, $D_\epsilon$ descends to $X$. Therefore, $D_\epsilon$ descends to $X$ for any $\epsilon > 0$. Moreover, for $\epsilon > 0$, $D_\epsilon$ is a nef line bundle that contracts only $F_1$ and $F_3^0([2])$, and is thus ample on $X$, hence semiample on $\M{2}{2}$.
\end{proof}

\begin{prop}
    For any $\epsilon>0$ and $0<a<1$, $D_{a,\epsilon}=\epsilon(12\lambda-\delta_{\text{irr}})+\psi_1+\psi_2+a\delta_{1, \left\{1\right\}}$ is semiample on $\M{2}{2}$. 
\end{prop}

\begin{proof}
    By \cref{prop:lambdapsi}, the stable base locus of $D_{a,\epsilon}$ is contained in $\Delta_{1, \{1\}}$. Let $\xi: \M{2}{2} \to \Mg{3}$ be the clutching map. Note that, by \cite[Proposition 3.3.6]{Rul01}, any nef line bundle on $\Mg{3}$ is semiample. In particular, the divisor $10\lambda - \delta_{\text{irr}} - 2\delta_1$ is semiample on $\Mg{3}$, and thus
    \begin{align*}
        \xi^\ast\left(10\lambda - \delta_{\text{irr}} - 2\delta_1\right)
        &= 10\lambda - \delta_{\text{irr}} + \psi_1 + \psi_2 - \delta_{1, \{1\}} - 2\delta_{1, \emptyset} - 2\delta_{0, \{1,2\}} \\
        &= \psi_1 + \psi_2 + \delta_{1, \{1\}} - 2\delta_{0, \{1,2\}}
    \end{align*}
    is semiample. Here, we used the relation $10\lambda - \delta_{\text{irr}} - 2\delta_1 = 0$ on $\M{2}{n}$ (cf.~\cite[Theorem 4]{AC09}). Then,
    \[
    D_{a, \epsilon} = \left(\epsilon(12\lambda - \delta_{\text{irr}}) + (1 - a)(\psi_1 + \psi_2)\right) + a\left(\psi_1 + \psi_2 + \delta_{1, \{1\}} - 2\delta_{0, \{1,2\}}\right) + 2a\delta_{0, \{1,2\}},
    \]
    so the stable base locus of $D_{a,\epsilon}$ is contained in $\Delta_{0,\{1,2\}}$. Since $\Delta_{0,\{1,2\}}$ and $\Delta_{1, \{1\}}$ are disjoint, the stable base locus of $D_{a, \epsilon}$ is empty, and thus $D_{a,\epsilon}$ is semiample.    
\end{proof}

\begin{prop}\label{prop:dim3semi}
    Any line bundle contained in the interior of the cone generated by
    \[ 12\lambda - \delta_{\text{irr}}, \psi_1+\psi_2,  \psi_1+\psi_2+\delta_{1,\left\{1\right\}}, \delta_{\text{irr}}+10\psi_1+10\psi_2+12\delta_{1, \emptyset}+2\delta_{1,\left\{1\right\}} \]
    is semiample.
\end{prop}

\begin{proof}
     Let $f:\M{2}{2}\to X$ be the contraction described in the proof of \cref{prop:lambdapsi}, i.e., the contraction corresponding to $F_1$ and $F_3^0([2])$. Since $\epsilon(12\lambda-\delta_{\text{irr}})+\psi_1+\psi_2$, $\epsilon(12\lambda-\delta_{\text{irr}})+\psi_1+\psi_2+a\delta_{1, \left\{1\right\}}$, and $12\lambda-\delta_{\text{irr}}$ are semiample line bundles that contract $F_1$ and $F_3^0([2])$, they descend to $X$ by \cref{lem:descend}. In particular, the Picard rank of $X$ is at least $3$.

    By \cref{thm:semin}, $f^\ast \text{Pic}(X)$ is strictly smaller than the face of the nef cone contracting $F_1$ and $F_3^0([2])$ (cf. proof of \cref{prop:lambdapsi}). Hence, the Picard rank of $X$ is exactly $3$. Therefore, $f^\ast \text{Pic}(X)$ is the intersection of the nef cone with the subspace spanned by $12\lambda-\delta_{\text{irr}}$, $\psi_1+\psi_2$, and $\delta_{1, \left\{1\right\}}$, which is precisely the cone described in the proposition. Its interior corresponds to pullback of ample line bundles on $X$, hence these are semiample.
\end{proof} 

\begin{prop}\label{prop:dim4semi}
    Any line bundle contained in the interior of the cone generated by
    \[ \lambda, 12\lambda - \delta_{\text{irr}}, \psi_1+\psi_2,  \psi_1+\psi_2+\delta_{1,\left\{1\right\}}, \delta_{\text{irr}}+10\psi_1+10\psi_2+12\delta_{1, \emptyset}+2\delta_{1,\left\{1\right\}} \]
    is semiample.
\end{prop}

\begin{proof}
    The proof is similar to the proof of \cref{prop:dim3semi}, but uses the contraction corresponding to $(12+\epsilon)\lambda-\delta_{\text{irr}}+\psi_1+\psi_2$, which only contracts $F_3^0([2])$.
\end{proof}

\begin{rmk}
\begin{enumerate}
    \item The contraction associated with $F_1$ is studied in \cite{CTV23, CTV21} and is known to have good properties, as $F_1$ is $K_X$-negative.
    
    \item Unlike the case of $g = 3$ stated in \cref{cor:main1cor}, there exists a contraction that contracts only the F-curve $F_3^0([2])$. The key difference is that, in the case $g \ge 3$, the codimension-one linear subspace spanning the face in \cref{thm:main1} is itself a face. However, in the case of $g = 2$, this is not true—the linear subspace lies in the interior of the cone.
    
    \item The contractions appearing in \cref{prop:dim3semi} and \cref{prop:dim4semi} have an interesting property: if we denote the contraction by $f: \M{2}{2} \to X$, then
    \[
    \NE{f}^\perp \subsetneq f^\ast \mathrm{Pic}(X).
    \]
    In particular, these contractions differ in nature from those considered in \cite{Cho24}, \cite[Section 4]{Cho25}, or \cref{sec:Knudsen}, where the equality $\NE{f}^\perp = f^\ast \mathrm{Pic}(X)$ holds. This serves as an explicit example illustrating why the condition $\NE{f}^\perp = f^\ast \mathrm{Pic}(X)$ is closely related to the semiampleness of divisors, as discussed in \cite{Cho24}.
\end{enumerate}

\end{rmk}

%\begin{thm}
%    $D=\psi_1+\psi_2+\delta_{1,\left\{1\right\}}$ is not semiample....?
%\end{thm}

%\begin{proof}
%    Assume that this is semiample, and let $f:\M{2}{2}\to \PP^n$ be the corresponding map to a projective space. Note that $D|_{\Delta_{0,\left\{1,2\right\}}}$ and $D|_{\Delta_{1,\left\{1\right\}}}$ is trivial. The second one follows from the fact that $\psi_1=\psi_2$ on $\M{1}{2}$. Hence, $f(\Delta_{0,\left\{1,2\right\}})$ and $f(\Delta_{1,\left\{1\right\}})$ are points. Moreover, $\Delta_{1,\left\{1,2\right\}}\simeq \M{1}{1}\times \M{1}{3}$ and
%    \[ D|_{\Delta_{1,\left\{1,2\right\}}}=\pi_2^\ast\left(\psi_1+\psi_2\right) \]
%    so $\dim f(\Delta_{0,\left\{1,2\right\}})=3$. TODO
%\end{proof}

%% file: sections/case01.tex
Here, we analyze the extremality of F-curves for genus $0$ and $1$, which is easier to prove and more explicit than the general genus case. We begin with some definitions and a lemma.

\begin{defn}\label{defn:convex}
    Let $C$ be a proper cone (i.e., a cone that does not contain a full straight line and has nonempty interior) in an $n$-dimensional Euclidean space. For any subset $S\subseteq C$, the \textbf{index of extremality} is defined as 
    \[
    I(S)=n-\dim S^\perp,
    \]
    where $S^\perp$ is the face of the dual cone $C^\ast$ orthogonal to $S$. A face $F\subseteq C$ is said to be \textbf{regular} if $\dim F=I(F)$. Moreover, if $F$ is a ray then it is said to be \textbf{regular extremal}.
\end{defn}

\begin{rmk}\label{rmk:regext}
    For any face $F$, $I(F)\ge \dim F$. If $C$ is a polyhedral cone, it is straightforward to verify that every face is regular extremal. However, this property does not hold in general. Consider the region 
    \[
    D=\left\{x^2+(|y|+1)^2\le 4 \right\}
    \]
    and let $C$ be the $3$-dimensional cone over $D$. In this setting, every ray on the boundary of $C$ is extremal, but only the rays corresponding to $(\pm\sqrt{3},0)$ are regular extremal.
\end{rmk}

Note that some previous papers have already implicitly utilized similar notions. For example, \cite{Mul21} proved that the closed cone of moving curves of $\M{g}{n}$ for $g,n\geq 2$ is not polyhedral (or equivalently, that the closed cone of pseudoeffective divisors is not polyhedral) by constructing an extremal ray $F$ of the moving cone with $2\le I(F)\le n$.

\begin{lem}\label{lem:convex}
    Let $X$ and $Y$ be projective, normal, $\Q$-factorial varieties, and let $f:X\to Y$ be a morphism. Assume the following conditions hold:
    \begin{enumerate}
        \item $\ker f^\ast$ is spanned by nef line bundles.
        \item Any nef line bundle on $X$ is the pullback of a nef line bundle on $Y$.
    \end{enumerate}
    Then $f_\ast$ maps any (regular) face of $\NE{X}$ to a (regular) face of $\NE{Y}$. In particular, for any face $F$ of $\NE{X}$, $I(F)=I(f_\ast F)$. 
\end{lem}

\begin{proof}
    Note that (2) implies that the induced map $f^\ast :\Pic(Y)\to \Pic(X)$ is surjective, or, equivalently, that $f_\ast:\NE{X}\to \NE{Y}$ is injective. This establishes that $I(F)\le I(f_\ast F)$ for any face $F$. Moreover, (1) implies that $f_\ast \A_1(X)_\R\cap \NE{Y}$ forms a face of $\NE{Y}$. By (2), we have 
    \[
    f_\ast (\A_1(X)_\R)\cap \NE{Y}=f_\ast \NE{X}.
    \]
    Thus, $f_\ast$ maps faces to faces. It remains to verify the second assertion.
    
    Condition (1) ensures that there exist $\rho(Y)-\rho(X)$ linearly independent nef line bundles on $Y$ that vanish on $X$. Furthermore, for any face $F$, there are $\rho(X)-I(F)$ linearly independent nef line bundles on $X$ that intersect $F$ trivially. By (2), we can extend these to nef line bundles on $Y$. Consequently, we obtain $\rho(Y)-I(F)$ independent line bundles that intersect $f_\ast F$ trivially. In particular, this yields the inequality $I(f_\ast F)\le I(F)$. Since we previously established that $I(F)\le I(f_\ast F)$, it follows that $I(F)=I(f_\ast F)$, completing the proof of the second assertion.

\end{proof}

\begin{rmk}
\cref{lem:convex} is a straightforward variant of \cite[Proposition 2.5]{CC15} and \cite[Lemma 2.7]{Bla22}. However, since it considers both the effective cone and its dual, \cref{lem:convex} has the advantage of establishing extremality within the closed cone of pseudoeffective cycles, not only the (possibly non-closed) cone of effective cycles.
\end{rmk}

Let $f_1:\M{0}{n}\to \M{0}{n+1}$ and $f_2:\M{1}{n}\to \M{1}{n+1}$ be the maps attaching a $3$-pointed, genus $0$ stable curve to a marked point, and let $f_3:\M{0}{n}\to \M{1}{n-1}$ be the map attaching an elliptic curve to a marked point. Note that $f_3$ is a special case of the flag map, $F:\M{0}{g+n}/S_g\to \M{g}{n}$, defined and studied in \cite{GKM02}.

\begin{thm}\label{thm:01push}
The maps $f_1$ and $f_2$ satisfy the conditions of \cref{lem:convex}. If the characteristic of the base field is not equal to $2$, then the same holds for $f_3$ and $F$.
\end{thm}

\begin{proof}
Note that $f_3$ is a special case of $F$, so it suffices to prove the claim for $f_1$, $f_2$, and $F$.

(1) First, consider $f_1$. It is well known that $\rho(\M{0}{n}) = 2^{n-1} - \binom{n}{2} - 1$ (e.g., \cite{AC09}). Hence, it is enough to produce
\[
\rho(\M{0}{n+1}) - \rho(\M{0}{n}) = 2^{n-1} - n
\]
linearly independent nef divisors on $\M{0}{n+1}$ which intersect $\Delta_{0, [n-1]}$ trivially. Let
\[
P = \left\{S\subseteq [n-1]\ \middle|\ |S^c|\ge 2 \right\},
\]
and for $S\in P$, define $\psi_S := \pi_S^\ast \psi_{n+1}$. These divisors are nef and intersect $\Delta_{0, [n-1]}$ trivially. It remains to prove that they are linearly independent. 

Assume
\[
D = \sum_{S\in P} a_S \psi_S = 0.
\]
We will prove that $a_S=0$ for all $S$ by induction on $|S^c|$. If $|S^c|=2$, i.e., $S^c=\left\{p,q\right\}$, then
\[
F_6^{0,0,0,0}(p,q,n,\left\{p,q,n\right\}^c)\cdot D = a_S = 0.
\]
Assume that $a_S=0$ for every $|S^c|\le m$. Then for any \(S\) such that $|S^c|=m+1$,
\[
F_6^{0,0,0,0}(S^c\cup \left\{n+1\right\}, S_1, S_2, S_3)\cdot D = a_S = 0
\]
for any nonempty $S_1$, $S_2$, and $S_3$. Therefore, the $\psi_S$ are linearly independent, and $f_1$ satisfies condition (1).

The proof for $f_2$ is similar. We have
\[
\rho(\M{1}{n+1}) - \rho(\M{1}{n}) = 2^n - 1.
\]
Define $P$ as above, and for $S\subsetneq [n-1]$, let $\psi_S^1 = \pi_S^\ast \psi_n$ and $\psi_S^2 = \pi_S^\ast \psi_{n+1}$. Then, it suffices to show that the divisors $\psi_S^1$, $\psi_S^2$ for $S\ne [n-1]$, together with $\pi_{[n-1]}^\ast \psi$, are linearly independent. Note that on $\M{1}{n}$, we have $\psi_1 = \psi_2$. The proof is similar to the case of $f_1$, so we omit it.

Now consider $F$. By \cite[Theorem 4.7]{GKM02}, there exists a nef divisor $D_{\text{GKM}}$ on $\M{g}{n}$ such that $F^\ast D_{\text{GKM}}=0$ and $D_{\text{GKM}}$ intersects F-curves of type 1–5 positively. Hence, $D_{\text{GKM}}\in \ker F^\ast$, and for any divisor $D\in \ker F^\ast$ and sufficiently large $r$, the divisor $D+rD_{\text{GKM}}$ is nef by \cite[Theorem 0.3]{GKM02}. Therefore, (1) holds.

(2) The maps $f_1$ and $f_2$ are sections of the projection map; hence, condition (2) is automatic for them. For $F$, this condition is proved in \cite[Theorem 0.7]{GKM02}.
\end{proof}

\begin{cor}\label{cor:ext01}
    Any F-curve on $\M{0}{n}$ or $\M{1}{n}$ spans a regular extremal ray of the closed cone of curves. Moreover, for each such F-curve, there exists a contraction of $\M{0}{n}$ or $\M{1}{n}$ whose relative cone of curves is precisely the ray spanned by the F-curve.
\end{cor}

\begin{proof}
    On $\M{0}{n}$, every F-curve arises as the image of the fundamental class of $\M{0}{4}$ under a composition of maps $f_1:\M{0}{n} \to \M{0}{n+1}$ for various values of $n$. By \cref{lem:convex} and \cref{thm:01push}, each such class spans an extremal ray. 
    
    We now consider $\M{1}{n}$ for general characteristic. Here, we can only use the maps $f_2$. Note that the type 1 F-curve is the pushforward of a curve on $\M{1}{1}$ via a sequence of $f_2$'s, and hence is regular extremal. Any type 5 or type 6 F-curve is the pushforward, along a sequence of $f_2$'s, of one of the following:

    \begin{enumerate}
        \item $F_5^{0,0}(\{1\}, \{2\})$ on $\M{1}{2}$,
        \item $F_5^{0,0}(\{1\}, \{2\})$ on $\M{1}{3}$,
        \item $F_6^{0,0,0,1}(\{1\}, \{2\}, \{3\}, \emptyset)$ on $\M{1}{3}$,
        \item $F_6^{0,0,0,1}(\{1\}, \{2\}, \{3\}, \{4\})$ on $\M{1}{4}$.
    \end{enumerate}

    Hence, it is enough to prove that these curves are regular extremal. For (1) and (3), note that they are contracted by all projections. Hence,
    \[
    \sum_{i=1}^n \pi_i^\ast \text{Pic}(\M{1}{n-1})
    \]
    contracts them. By \cite[Lemma 4.3]{Cho25}, this is a codimension $1$ subspace. Therefore, they are regular extremal.
    
    For (2) and (4), these curves are contracted by divisors in
    \[
    \sum_{i=1}^{n-1} \pi_i^\ast \text{Pic}(\M{1}{n-1}).
    \]
    This space has codimension $2$, which can be shown using \cite[Theorem 4.1]{Cho25} or via direct computation, at least for $n = 3, 4$. Moreover, they are also contracted by $\psi_n$, which is not contained in the above space by \cite[Lemma 4.3]{Cho25}. Therefore, they are also regular extremal.    

    For the second assertion, we proceed by induction on $n$. The base cases are well known. As shown in the proof of \cref{thm:01push}, the kernels of $f_1^\ast$ and $f_2^\ast$ are generated by $\psi$ classes and their pullbacks. Hence, by \cref{thm:psisemi}, they are generated by semiample divisors. Moreover, by the same reasoning as in \cref{thm:01push}, any semiample divisor on $\M{0}{n}$ or $\M{1}{n}$ is the pullback of a semiample divisor via $f_1^\ast$ or $f_2^\ast$. Consequently, by the same argument used in \cref{lem:convex}, the orthogonal complement $F^\perp$ is spanned by semiample divisors. Therefore, the product of the corresponding morphisms yields the desired contraction.
\end{proof}

We will provide another proof of \cref{cor:ext01} in \cref{subsec:56}.

\begin{rmk}\label{rmk:cont}
    By the proof of \cref{cor:ext01}, for any F-curve on $\M{0}{n}$ or $\M{1}{n}$, its orthogonal complement is spanned by either pullbacks of $\psi$-classes or ample line bundles along projection maps. Therefore, on $\M{0}{n}$, the corresponding morphism can be constructed as the product of all contractions of the form $f_{\text{Kap}} \circ \pi_S$ that contract $F$, where $f_{\text{Kap}}$ denotes Kapranov’s construction \cite{Ka92, Ka93} associated with the class $\psi_i$. An analogous argument reduces the construction of corresponding morphisms in genus one to \cref{qes:psisemi}.
\end{rmk}

%% file: sections/newnef.tex
Here, we introduce variants of the $\kappa$ class, which form a new family of nef divisors on $\M{g}{n}$. We refer to the beginning of \cref{subsec:method} for a motivation for constructing such divisors.

Note that $\kappa$ is an ample divisor on $\M{g}{n}$ and satisfies
\[
\kappa=12\lambda-\delta+\psi
\]
by Mumford's formula. Define
\[
B_{g,n}=\left\{ \delta_{i,I}\ |\ 0\le i\le g, I\subseteq [n], |I|\ge 1 \text{ if }i=0\text{ and }I\ne [n]\text{ if }i=g \right\},
\]
i.e., the set of all boundary classes except $\delta_{\text{irr}}$, including $\delta_{0,i}=-\psi_i$. Note that we identify $\delta_{i,I}$ and $\delta_{g-i, I^c}$.

\begin{defn}\label{defn:ensemble}
    A subset $D\subseteq B_{g,n}$ is called a \textbf{semigroup of boundary divisors} (semigroup for short, if there is no confusion) if it satisfies the condition:
    \[
    \delta_{i,I}, \delta_{j,J}\in D, \quad i+j\le g, \quad I\cap J=\emptyset\ \Rightarrow\ \delta_{i+j, I\cup J}\in D.
    \]
    when it makes sense, i.e. if $(i+j, I\cup J)\ne (0, \emptyset), (g, [n])$. 
    
    The corresponding \textbf{semigroup kappa divisor} is defined as
    \[
    \kappa_{g,n}(D):=\kappa+\sum_{\delta_{i,I}\in D}\delta_{i,I}.
    \]
\end{defn}

\begin{rmk}
    We identify $\delta_{i,I}$ with $\delta_{g-i, I^c}$, so the condition for a semigroup must also hold under this identification. For example, since $\delta_{3, \emptyset}=\delta_{g-3, [n]}$, the inclusion $\delta_{2, \emptyset},\delta_{3, \emptyset}\in D$ implies not only $\delta_{5, \emptyset}\in D$ but also that $\delta_{g-1, [n]}=\delta_{1, \emptyset}\in D$.
\end{rmk}

Due to this symmetry, it is sometimes cumbersome to check whether a set of boundary divisors forms a semigroup. The following lemma is helpful for this purpose.

\begin{lem}\label{lem:mult}
    Let $D \subseteq B_{g,n}$ and $p \in [n]$. Then $D$ is a semigroup if and only if the following holds: for any $\delta_{i,I}, \delta_{j,J}\in D$ such that $p \not\in I, J$,
    \begin{enumerate}
        \item If $I \cap J = \emptyset$ and $i + j \le g$, then $\delta_{i+j, I \cup J} \in D$.
        \item If $I \subseteq J$, $i \le j$, and $(i, I) \ne (j, J)$, then $\delta_{j - i, J \setminus I} \in D$.
    \end{enumerate}
\end{lem}

\begin{proof}
    First, assume that $D$ is a semigroup. Since $p \not\in I \cup J$, condition (1) follows directly from the semigroup property. For (2), note that under the given assumptions, we have $\delta_{g - j, J^c} = \delta_{j, J} \in D$, so by the semigroup condition, it follows that $\delta_{j - i, J \setminus I} \in D$.

    Conversely, assume that $\delta_{i,I}, \delta_{j,J} \in D$ with $i + j \le g$, $I \cap J = \emptyset$, and $(i + j, I \cup J) \ne (0, \emptyset), (g, [n])$. If $p \not\in I \sqcup J$, then $\delta_{i+j, I \cup J} \in D$ by (1). If $p \in I \sqcup J$, we may assume without loss of generality that $p \in J$ and $p \not\in I$. Then $\delta_{g - j, J^c} = \delta_{j, J} \in D$, and since $i \le g - j$ and $I \subsetneq J^c$, condition (2) implies $\delta_{g - j - i, J^c \setminus I} = \delta_{i + j, I \cup J} \in D$. Hence $D$ is a semigroup.
\end{proof}

Now we will prove that semigroup kappa divisors are nef. The following well-known lemma is useful in general.

\begin{lem}\label{lem:bdry}
    For the boundary divisor
    \[
        \Delta_{i, I}\simeq \M{i}{I+1}\times \M{g-i}{I^c+1},
    \]
    choose the attaching maps $\theta_1:\M{i}{I+1}\to \M{g}{n}$ and $\theta_2:\M{g-i}{I^c+1}\to \M{g}{n}$. A divisor $D$ on $\M{g}{n}$ is nef (resp.~semiample, ample) on $\Delta_{i, I}$ if and only if $\theta_1^\ast D$ and $\theta_2^\ast D$ are nef (resp.~semiample, ample).
\end{lem}

\begin{proof}
    First, we show that
    \[
        \Pic(\M{i}{I+1}) \times \Pic(\M{g-i}{I^c+1}) \simeq \Pic(\M{i}{I+1}\times \M{g-i}{I^c+1}).
    \]
    Choose a prime $l$ different from the characteristic of the base field. It suffices to prove this after tensoring with $\Q_l$. By \cite[Theorem 0.1]{Mor01}, we have 
    \[
    \HH^1_{\text{\emph{ét}}}(\M{g}{n}, \Q_l) = 0 
    \quad \text{and} \quad 
    \HH^2_{\text{\emph{ét}}}(\M{g}{n}, \Q_l) \simeq \Pic(\M{g}{n}) \otimes_\Q \Q_l
    \]
    for any $g,n$. Therefore, the statement follows from the Künneth formula. Hence, on $\Delta_{i,I}$ we have 
    \[
        D = \pi_1^\ast \theta_1^\ast D + \pi_2^\ast \theta_2^\ast D.
    \]
    The desired equivalence is then immediate.

\end{proof}

\begin{thm}\label{thm:kapnef}
    $\kappa_{g,n}(D)$'s are nef. Moreover, in positive characteristic, they are semiample.
\end{thm}

\begin{proof}
    \textbf{Step 1. }Let $D$ be a semigroup on $\M{g}{n}$ such that $\delta_{0.j}\in D$. Then there exists a semigroup $D'$ on $\M{g}{n-1}$ such that $\kappa_{g,n}(D)=\pi_j^\ast \kappa_{g,n-1}(D')$.
    
    We may assume that $j = n$. Define a set of boundary divisors on $\M{g}{n-1}$:
    \[
    D' := \left\{ \delta_{i, I}\ |\ I \subseteq [n-1],\ \delta_{i. I} \in D \right\}.
    \]
    First, we will prove that this is a well-defined semigroup. We need to show that if $\delta_{i. I} \in D'$, so $I \subseteq [n-1]$ and $\delta_{i. I} \in D$ as a divisor on $\M{g}{n}$, then $\delta_{g-i. [n-1] \setminus I} \in D'$, i.e., $\delta_{g-i. [n-1] \setminus I} \in D$ as a divisor on $\M{g}{n}$. Since $\delta_{i. I} \in D$, we have $\delta_{g-i. [n] \setminus I} \in D$. Also, by $\delta_{0,n} \in D$, it follows from the semigroup condition that $\delta_{g-i. [n-1] \setminus I} \in D$. Therefore, $D'$ is well-defined. The semigroup condition for $D'$ then follows directly from the condition for $D$.
    
    Now let
    \[
    \tilde{D}=\left\{(i, I)\ |\ 0\le i\le g,\ I\subseteq [n-1],\ (i, I)\ne (0, \emptyset),\ (g, [n-1])\text{ and }\delta_{i, I}\in D' \right\}.
    \]
    Note that $(i, I)\mapsto \delta_{i,I}$ is a two-to-one map from $\tilde{D}$ to $D'$, and 
    \[
    D=\left\{\delta_{i, I}\ |\ (i, I)\in \tilde{D} \right\}\sqcup  \left\{\delta_{0, n} \right\}= \left\{\delta_{i, I\cup \left\{n\right\}}\ |\ (i, I)\in \tilde{D} \right\}\sqcup \left\{\delta_{0, n} \right\}.
    \]
    This is straightforward, since if $\delta_{i, I}\in D$ and $(i, I)\ne (0, n),\ (g, [n-1])$, then exactly one of $(i, I)$ or $(g-i, I^c)$ lies in $\tilde{D}$, depending on whether $I$ contains $n$ or not. Therefore,
    \begin{align*}
        \pi_n^\ast \kappa_{g,n}(D') &= \pi_i^\ast \kappa + \frac{1}{2}\sum_{(i, I)\in \tilde{D}}\pi_n^\ast \delta_{i,I} = \kappa - \psi_n + \frac{1}{2}\sum_{(i, I)\in \tilde{D}} \left(\delta_{i,I} + \delta_{i, I\cup \left\{n\right\}}\right) \\
        &= \kappa + \frac{1}{2}\left(\sum_{(i, I)\in \tilde{D}} \delta_{i,I} - \psi_n\right) + \frac{1}{2}\left(\sum_{(i, I)\in \tilde{D}} \delta_{i,I\cup \left\{n\right\}} - \psi_n\right)= \kappa_{g,n-1}(D)
    \end{align*}
    so the conclusion follows. We have used \cite[Lemma 1(i)]{AC09} in the computation.

    \textbf{Step 2. }Let $D$ be a semigroup on $\M{g}{n}$ and let $\theta:\M{g'}{n'+1} \to \M{g}{n}$ be a map attaching a fixed stable curve. Then $\theta^\ast \kappa_{g,n}(D)$ is also a semigroup kappa divisor.
    
    The map $\theta$ is a composition of attaching maps that either increase the number of marked points (i.e., $n' < n$) or attach an elliptic tail. Hence, it suffices to prove the statement for such attaching maps.
    
    First, assume $n' < n$. Then, by \cite[Lemma 1]{AC09},
    \[
    \theta^\ast \kappa_{g,n}(D) = \kappa + \sum_{i \le g',\, I \subseteq [n']} \delta_{i,I}.
    \]
    Thus, it is enough to show that $D' = \left\{ \delta_{i,I} \ |\ i \le g',\ I \subseteq [n'] \right\}$ is a semigroup on $\M{g'}{n'+1}$. It is straightforward to verify that $D'$ satisfies the two conditions of \cref{lem:mult} with $p = n'+1$, using the corresponding conditions for $D$.
    
    Second, let $\theta:\M{g-1}{n+1} \to \M{g}{n}$ be the map attaching an elliptic curve to the $(n+1)$-th point. Then, by \cite[Lemma 1]{AC09},
    \[
    \theta^\ast \delta_{i,I} = 
    \begin{cases}
    \delta_{i,I} + \delta_{i-1, I \cup \{n+1\}} & \text{if } 1 \le i \le g - 1, \\
    \delta_{i,I} & \text{if } i = 0, \\
    \delta_{i-1, I \cup \{n+1\}} & \text{if } i = g.
    \end{cases}
    \]
    We choose the representative of $D$ as
    \[
    \tilde{D} := \left\{ (i,I)\ \mid\ n \not\in I, \delta_{i, I}\in D \right\},
    \]
    so for any element $x \in D$, there exists a unique $(i,I) \in \tilde{D}$ such that $x = \delta_{i,I}$. Then it suffices to show that
    \begin{align*}
    D' &= \left\{ \delta_{i, I} \ \mid\ i \le g - 1,\ (i,I) \in \tilde{D} \right\} 
    \sqcup \left\{ \delta_{i - 1, I \cup \{n+1\}} \ \mid\ 1 \le i,\ (i,I) \in \tilde{D} \right\} \\
    &= \left\{ \delta_{i, I} \ \mid\ i \le g - 1,\ (i,I) \in \tilde{D} \right\} 
    \sqcup \left\{ \delta_{g - i, [n] \setminus I} \ \mid\ 1 \le i,\ (i,I) \in \tilde{D} \right\}
    \end{align*}
    is a semigroup on $\M{g - 1}{n + 1}$. Note that, as a divisor on $\M{g}{n}$, we have $\delta_{g - i, [n] \setminus I} = \delta_{i, I} \in D$. Hence, the two conditions of \cref{lem:mult} for $D'$ with $p = n + 1$ follow from the corresponding conditions for $D$.
    
    \textbf{Step 3. }Proof of the theorem
    
    Consider the first assertion. We will use induction on $3g - 3 + n$, i.e., the dimension of $\M{g}{n}$. The base cases are $(g,n) = (0,4)$ and $(1,1)$. For $(g,n) = (1,1)$, any semigroup kappa divisor is either $\kappa$ or $\kappa - \psi_1$, both of which are nef. For $(g,n) = (0,4)$, we prove a stronger result: any semigroup kappa divisor intersects type 6 F-curves nonnegatively.
    
    By the formula \cref{thm:explicit}, the intersection of $\kappa$ with any F-curve is $1$. Choose a type 6 F-curve $F_6^{g_1, g_2, g_3, g_4}(I_1, I_2, I_3, I_4)$ on $\M{g}{n}$. If $0$ (resp. $1,2,3,4$) among the $\delta_{g_i,I_i}$ are contained in $D$, then at least $0$ (resp. $0,1,3,3$) among $\delta_{0, I_1 \cup I_2},\ \delta_{0, I_1 \cup I_3},\ \delta_{0, I_1 \cup I_4}$ are contained in $D$ by the definition of a semigroup. Hence, the intersection of $\sum_{\delta_{i,I} \in D} \delta_{i,I}$ with the $F$-curve is at least $-1$. Therefore, the intersection of $\kappa_{g,n}(D)$ with any type 6 F-curve is nonnegative.
    
    Now we proceed with the induction step. If $D$ contains $\delta_{0,i}$ for some $i$, then by Step 1, it is the pullback of a semigroup kappa divisor along the projection map $\pi_i$. Hence, $\kappa_{g,n}(D)$ is nef by the induction hypothesis. Now assume that $D$ contains no $\delta_{0.i}$'s. Then $\kappa$ is ample and $\sum_{\delta_{i,I} \in D} \delta_{i,I}$ is effective. Thus, it is enough to prove that $\kappa_{g,n}(D)|_{\Delta_{i,I}}$ is nef for every boundary divisor of form $\Delta_{i,I}$.
    
   To prove that $\kappa_{g,n}(D)$ is nef on this divisor, by \cref{lem:bdry}, it is enough to show that $\theta^\ast \kappa_{g,n}(D)$ is nef for any attaching map $\theta: \M{g'}{n'+1} \to \M{g}{n}$. By Step 2, this pullback is also a semigroup kappa divisor. By the induction hypothesis, it is nef. This completes the proof of the first assertion.
    
    Now consider the second assertion. Again, we use induction. Since any nef line bundle on $\M{0}{4}$ and $\M{1}{1}$ is semiample, the base cases follow from the first assertion. Now consider the induction step. If $D$ contains $\delta_{0,i}$ for some $i$, then the same argument applies and shows the semiampleness. If $D$ contains no $\delta_{0.i}$'s, then $\kappa_{g,n}(D)$ is the sum of an ample and an effective divisor, so it is big. Since it is big and nef, we can apply \cite{Ke99}. The exceptional locus of $\kappa_{g,n}(D)$ is contained in the union of boundary divisors of the form $\Delta_{i,I}$. Hence, it suffices to prove that $\kappa_{g,n}(D)|_{\Delta_{i,I}}$ is semiample. The proof is exactly the same as in the nef case.

\end{proof}

Note that the role of kappa divisors in this proof is twofold:  (1) they provide a nef divisor that is stable under pullback by clutching maps, and  (2) they ensure that the divisors \( \kappa_{g,n}(D) \) are F-nef. This naturally raises the question of whether one can replace \( \kappa \) with other divisors. The second property is not essential, as one may instead restrict attention to those divisors that are already F-nef. This leads to the following.

\begin{thm}
    Let \( (L_{g,n})_{2g - 2 + n > 0} \) be a family of nef line bundles such that for any clutching map $\theta:\M{g'}{n'+1} \to \M{g}{n}$, we have \( \theta^\ast L_{g,n} = L_{g',n'+1} \). Then
    \[
    L_{g,n}(D) := L_{g,n} + \sum_{\delta_{i,I} \in D} \delta_{i,I}
    \]
    is nef whenever it is F-nef. Moreover, if \( L_{g,n} \) is big and nef for every $g,n$ satisfying $2g-2+n>0$, then \( L_{g,n}(D) \) is semiample in positive characteristic if it is F-nef.
\end{thm}

For example, one may take any nef linear combination of \( \lambda \), \( 12\lambda - \delta_{\mathrm{irr}} \), and \( \kappa \).  
We omit the proof, as it proceeds identically to that of \cref{thm:kapnef}.

\begin{eg}\label{eg:kapeg}
    \begin{enumerate}
    \item If $I\ne \emptyset, [n]$, or $I=\emptyset$ and $2i>g$, or $I=[n]$ and $2i<g$, then $\kappa+\delta_{i,I}$ is a nef divisor on $\M{g}{n}$.
    \item If $(\delta_{i,I}, \delta_{j,J}, \delta_{k,K})$ is a triple such that each divisor satisfies the condition of (1), $i+j+k=g$ and $I\sqcup J\sqcup K=[n]$, then $\kappa+\delta_{i,I}+ \delta_{j,J}+\delta_{k,K}$ is nef.
    \item $\kappa+\sum \delta_{i, \emptyset}$ is nef. Moreover,  $\kappa+\sum_{i\text{ even}} \delta_{i, \emptyset}$ is also nef. However, $\kappa+\sum_{i\text{ odd}} \delta_{i, \emptyset}$ is not nef in general. Similarly, for any $d\in \Z_{>0}$, $\kappa+\sum_{d\mid i} \delta_{i, \emptyset}$ is nef.
    \item Let $I, J \subseteq [n]$ be such that $I \setminus J$, $J \setminus I$ and $I\cap J$ are nonempty. Then, for any $i, j \le g$, the pair $(\delta_{i, I}, \delta_{j, J})$ is a semigroup if $I\cup J\ne [n]$ or $i+j<g$. Hence, $\kappa + \delta_{i, I} + \delta_{j, J}$ is nef in this case.
    \end{enumerate}
\end{eg}

\begin{rmk}\label{rmk:symm}
    Any symmetric semigroup kappa divisor on $\M{0}{n}$ is of the following form: for a positive integer $d \mid n$,
    \[
    D_{d,n} := \kappa + \sum_{d \mid |I|} \delta_{0, I}.
    \]
    For $d = 2,3$, this is related to other divisors that have already been studied: type A, level 1 affine coinvariant divisors \cite{Fak12} and type A parafermion divisors \cite{Cha25}.
    
    Let $F_{1,1,i}$ be an F-curve of the form $F_6^{0,0,0,0}(I_1, I_2, I_3, I_4)$ where $|I_1| = |I_2| = 1$ and $|I_3| = i$, viewed as a class in $\mathrm{A}_1(\M{0}{n}/S_n)$. Since we are considering the symmetrized situation, the class does not depend on the exact choice of the $I_j$'s. Note that, by \cite[Corollary 2.2]{AGS14}, these classes form a basis of $\mathrm{A}_1(\M{0}{n}/S_n)$. Hence, to compare symmetric divisors, it is enough to compute their intersections with the $F_{1,1,i}$'s. Using \cref{thm:explicit}, we obtain
    \[
    D_{2,n} \cdot F_{1,1,i} =
    \begin{cases}
        0 & \text{if } 2 \mid i, \\
        4 & \text{if } 2 \nmid i,
    \end{cases}
    \]
    and
    \[
    D_{3,n} \cdot F_{1,1,i} =
    \begin{cases}
        3 & \text{if } i \equiv 2 \pmod{3}, \\
        0 & \text{otherwise}.
    \end{cases}
    \]
    By comparing these with \cite[Proposition 5.2]{Fak12} and \cite[Propositions 4.41, 4.42, 4.48, 4.49]{Cha25}, we find that $D_{2,n}$ (resp. $D_{3,n}$) coincides, up to a positive constant, with the $\mathfrak{sl}_2$ (resp. $\mathfrak{sl}_3$) level 1 affine coinvariant divisor, and with the $\mathfrak{sl}_{r}$ level $2$ (resp. level $3$) parafermion coinvariant divisor. Hence, by \cite[Lemma 2.5]{Fak12}, they are also semiample in characteristic zero.
    
    The notion underlying this is that of divisors arising from symmetric functions, as defined in \cite{Fe15}. Symmetric semigroup kappa divisors on $\M{0}{n}$, level 1 affine coinvariant divisors, and certain parafermion coinvariant divisors are all examples of divisors from symmetric functions. They lie in the same abelian group but arise from different symmetric functions, and in some simple cases, they coincide up to a positive constant.

\end{rmk}

As we will see in \cref{sec:extremality}, semigroup kappa divisors are useful for constructing nef divisors that contract certain F-curves. In view of \cref{cor:ext01}, \cref{rmk:symm}, and the semiampleness of semigroup kappa divisors in characteristic~$p$, the following questions naturally arise:

\begin{qes}\label{qes:mult}
    \begin{enumerate}
        \item In characteristic zero, which semigroup kappa divisors are semiample?
        \item In positive characteristic, is there a modular interpretation of the morphism defined by semigroup kappa divisors?
    \end{enumerate}
\end{qes}

From now on, we will list some sporadic new nef divisors, which will be used in \cref{sec:extremality}.

\begin{prop}[char $k\ne 2$]\label{prop:14nef}
    For any $0 \leq \epsilon_{i, I} \leq \frac{1}{4}$,
    \[
    \kappa + \sum_{i,I} \epsilon_{i, I} \delta_{i,I}
    \]
    is nef, where $\delta_{i, I}$'s include $\delta_{0.i}$'s.
\end{prop}

\begin{proof}
    First, we will prove that this divisor is F-nef. The intersection of $\kappa$ with any F-curve is $1$. By \cref{thm:explicit}, the intersection of $\sum_{i,I} \epsilon_{i, I} \delta_{i,I}$ with an F-curve is at least the sum of four terms of the form $-\epsilon_{i,I}$, hence at least $-1$. Therefore, the total intersection is nonnegative, and the divisor is F-nef.

    We will first prove that $\kappa - \frac{1}{4} \sum_{i \in I} \psi_i$ is nef, where $I \subseteq [n]$. Note that this is a special case of the divisor in the statement. Let $F: \M{0}{g+n} \to \M{g}{n}$ be the flag map. Then, by \cite[Lemma 1]{AC09}, the pullback of the divisor along $F$ is a divisor of the same form. Since this divisor is F-nef, by \cite[Theorem 0.3]{GKM02}, it is enough to check the genus $0$ case.

    Since the $\psi_i$ are nef, it is enough to show that $\kappa - \frac{1}{4} \psi$ is nef on $\M{0}{g}$. Let $F: \M{0}{g} \to \Mg{g}$ be the flag map, and consider the divisor $\kappa + \frac{1}{4} \delta_1$. By \cite[Lemma 1]{AC09}, we have $F^\ast(\kappa + \frac{1}{4} \delta_1) = \kappa - \frac{1}{4} \psi$, so it suffices to show that $\kappa + \frac{1}{4} \delta_1$ is nef. Since this is a special case of the divisor in the statement, it is F-nef. Using the relation $\kappa = 12\lambda - \delta$ on $\Mg{g}$, this follows from \cite[Corollary 5.4]{Gib09}.

    Now we prove the general assertion. As in the proof of \cref{thm:kapnef}, we use induction on $3g - 3 + n$. The base case follows from the F-nefness.

    Consider the induction step. Since $\kappa - \frac{1}{4} \sum_{i \in I} \psi_i$ is nef, the divisor in the statement is a sum of a nef divisor and an effective sum of boundary divisors of the form $\Delta_{i,I}$. Hence, it is enough to show that the restriction to $\Delta_{i,I}$ is nef. This follows from \cref{lem:bdry}: by \cite[Lemma 1]{AC09}, the pullback of the divisor along an attaching map $\theta$ is of the same form, and the induction hypothesis applies.
\end{proof}

% Now, assume $n\ge 1$ and define 
% \[
% B_{1, \left\{n\right\}}^0:=\left\{\delta_{j,J}\ \middle|\ I\subseteq [n-1] , \ j\leq g-1 \right\}.
% \]
% This set has a natural involution given by $\delta_{j, J}\mapsto \delta_{g-1-j,[n-1]\setminus J}$. Let $B_{1, \left\{n\right\}}$ be the quotient of $B_{1, \left\{n\right\}}^0$ by this involution and $\pi:B_{1, \left\{n\right\}}^0\to B_{1, \left\{n\right\}}$ be the quotient map.

%\begin{thm}\label{thm:14nef}
%    For any subset $D\subseteq D_{i, I}^0$, 
%    \[
%    \kappa+\delta_{i, I}+\frac{1}{8}\sum_D \left(\delta_{j,J}+\delta_{g-i-j, I^c}\right)
%    \]
%    is nef.
%\end{thm}

%\begin{proof}
%    TODO
%\end{proof}

Now, let $\delta_{k,K}$ be a boundary divisor satisfying the condition of \cref{eg:kapeg} (1)  and
\[
B_{k,K}^0 := \left\{\delta_{j,J} \ \middle|\ J\subseteq K^c,\ j\leq g-k,\  J=\emptyset \Rightarrow k<j, J=K^c\Rightarrow g-2k>j \right\}.
\]
This set has a natural involution given by $\delta_{j, J} \mapsto \delta_{g-k-j,\,[n]\setminus \left(K\cup J\right)}$. Let $B_{k, K}$ be the quotient of $B_{k,K}^0$ by this involution, and let $\pi:B_{k,K}^0 \to B_{k,K}$ be the quotient map. 

\begin{prop}[char $k\ne 2$]\label{prop:14new}
    For any subset $B\subseteq B_{k,K}$,
    \[
    D_B := \kappa+\delta_{k,K}+\frac{1}{4}\sum_{\pi(\delta_{i, I})\in B} \delta_{i, I}
    \]
    is nef.
\end{prop}

\begin{proof}
    By \cref{prop:14nef}, $D_B$ is a sum of a nef line bundle and $\delta_{k,K}$. Hence, it is enough to show that $D_B|_{\Delta_{k, K}}$ is nef. Let 
    \[
    \theta_1:\M{k}{K\cup \{p\}}\to \M{g}{n} \quad\text{and}\quad \theta_2:\M{g-k}{K^c\cup \{p\}}\to \M{g}{n}
    \]
    be the relevant attaching maps. By \cref{lem:bdry}, it is enough to show that $\theta_1^\ast D_B$ and $\theta_2^\ast D_B$ are nef. We have
    \[
    \theta_1^\ast D_B = \kappa - \psi_p = \pi_p^\ast \kappa,
    \]
    which is nef. Also,
    \[
    \theta_2^\ast D_B = \kappa - \psi_p + \frac{1}{4}\sum_{\pi(\delta_{i, I})\in B} \delta_{i, I} = \pi_p^\ast \left( \kappa+ \frac{1}{4}\sum_{\delta_{i, I}\in B} \delta_{i, I} \right),
    \]
    and by \cref{prop:14nef}, this is also nef. Therefore, $D_B$ is nef.
\end{proof}

\begin{prop}[char $k\ne 2$]\label{prop:last}
    On $\M{g}{1}$ with $g\ge 3$,
    \[
        D:=\kappa+\delta_{1, \{1\}}+\frac{1}{4}\left(\delta_{1, \emptyset}+\delta_{g-2, \emptyset} \right)
    \]
    is nef.
\end{prop}

\begin{proof}
    Most of the proof is the same as in \cref{prop:14new}. The difference is that $\theta_1^\ast D$ is not simply $\kappa - \psi_p$. In this case,
    \[
        \theta_1^\ast D=\kappa - \psi_p+\frac{1}{4}\delta_{1,\emptyset}
    \]
    which is nef on $\M{1}{2}$. Hence, $D$ is also nef.
\end{proof}

\begin{prop}\label{prop:3nnew}
    \[ D:=\kappa+\frac{2}{3}\delta_{1, \emptyset}+\frac{1}{3}\delta_{2, \emptyset} \]
    is a nef divisor on $\M{3}{n}$.
\end{prop}

\begin{proof}
    First, we will prove that this is F-nef. Since the intersection of $\kappa$ with any F-curve is $1$, by the explicit intersection formula in \cref{thm:explicit}, it is straightforward to see that the intersection of $D$ with any F-curve of types $1$–$4$ is nonnegative. By a similar argument, we can show that the intersection with a type~5 F-curve is also nonnegative, except in the case $F_5^{1,1}(\emptyset, \emptyset)$. Since the coefficient of $\delta_{2, \emptyset}$ is $\frac{1}{3}$, a direct computation gives $D \cdot F_5^{1,1}(\emptyset, \emptyset) = 0$.
    
    Now consider the type~6 case. The main factor is, among the $\delta_{g_i, I_i}$'s, how many of them are $\delta_{1, \emptyset}$ or $\delta_{2, \emptyset}$. The problematic case occurs when $2$ or $3$ among them are $\delta_{1, \emptyset}$. In these cases, we can explicitly compute the intersection number using \cref{thm:explicit}, and the result is $0$ in both cases. Hence, $D$ is F-nef.

    Since $\kappa$ is ample, it is enough to show that $D|_{\Delta_{1, \emptyset}}$ and $D|_{\Delta_{2, \emptyset}}$ are nef. Let 
    \[
    \theta_1:\M{1}{1} \to \M{3}{n}, \quad 
    \theta_2:\M{2}{n+1} \to \M{3}{n}, \quad
    \theta_3:\M{2}{1} \to \M{3}{n}, \quad
    \theta_4:\M{1}{n+1} \to \M{3}{n}
    \]
    be the relevant attaching maps. By \cref{lem:bdry}, it suffices to prove that $\theta_i^\ast D$ is nef for $i = 1, 2, 3, 4$. Since $D$ is F-nef, each $\theta_i^\ast D$ is also F-nef. In particular, $\theta_1^\ast D$ and $\theta_3^\ast D$ are nef by the known cases of the F-conjecture.

    We have
    \[
    \theta_2^\ast D = \kappa - \frac{2}{3}\psi_{n+1} + \frac{1}{3}\delta_{2, \emptyset} + \frac{1}{3}\delta_{1, \{n+1\}}
    \]
    on $\M{2}{n+1}$, by \cite[Lemma 1]{AC09}. Since $\kappa - \psi_{n+1} = \pi_{n+1}^\ast \kappa$ is nef, it is enough to show that the restriction of $\theta_2^\ast D$ to $\Delta_{2, \emptyset}$ and $\Delta_{1, \{n+1\}}$ is nef. As in the previous paragraph, there are four pullbacks to consider, but two of them are trivially nef by the known cases of the F-conjecture. The two other nontrivial cases are
    \[
    \alpha_1:\M{0}{n+2} \to \M{2}{n+1}, \quad \alpha_2:\M{1}{n+1} \to \M{2}{n+1}.
    \]
    Then we have
    \[
    \alpha_1^\ast \theta_2^\ast D = \kappa - \frac{2}{3}\psi_{n+1} - \frac{1}{3}\psi_{n+2}, \quad
    \alpha_2^\ast \theta_2^\ast D = \kappa - \frac{1}{3}\psi_{n+1}.
    \]
    Since the divisors $\kappa - \psi_i$ are nef, both of these are nef. Hence, $\theta_2^\ast D$ is nef.

    Next, we have
    \[
    \theta_4^\ast D = \kappa - \frac{2}{3}\psi_{n+1} + \frac{1}{3}\delta_{1, \emptyset}
    \]
    on $\M{1}{n+1}$. By the same argument as above, it is enough to show that $\beta^\ast \theta_4^\ast D$ is nef for 
    \[
    \beta:\M{0}{n+2} \to \M{1}{n+1}.
    \]
    Since
    \[
    \beta^\ast \theta_4^\ast D = \kappa - \frac{2}{3}\psi_{n+1} - \frac{1}{3}\psi_{n+2},
    \]
    this is also nef. Hence, $\theta_4^\ast D$ is nef, and therefore $D$ is nef.

\end{proof}

\begin{prop}[char $k\ne 2$]\label{prop:4nnef}
    \[ D:=\kappa+\frac{1}{2}\delta_{1, \emptyset}+\frac{1}{2}\delta_{2, \emptyset}  \]
    is a nef divisor on $\M{4}{n}$.
\end{prop}

\begin{proof}
    \textbf{Step 1. } On $\M{0}{n}$ with $n\ge 4$, 
    \[ D_0:=\kappa-\frac{1}{2}\sum_{i=1}^4 \psi_i+\frac{1}{2}\sum_{\{i,j\}\subset [4]}\delta_{0, \{i,j\} } \]
    is nef.

    Since the case $n=4$ is straightforward, we may assume $n\ge 5$. For any $\{i,j\}\subset [4]$, 
    \[ \kappa-\psi_i-\psi_j+\delta_{0, \{i,j\} }=\pi_i^\ast\pi_j^\ast \kappa \]
    is nef. By averaging these divisors, we find that
    \[  \kappa-\frac{1}{2}\sum_{i=1}^4 \psi_i+\frac{1}{6}\sum_{\{i,j\}\subset [4]}\delta_{0, \{i,j\} } \]
    is nef. Hence, for any $\Delta_{0,\{i, j\}}\simeq \M{0}{([n]\setminus \{i,j\})\cup p}$, it suffices to prove that $D_0|_{\Delta_{0,\{i,j\}}}$ is nef. Without loss of generality, we may assume $(i,j)=(3,4)$. On $\M{0}{([n]\setminus \{3,4\})\cup p}$, $D_0|_{\Delta_{0,\{3,4\}}}$ is 
    \[ \kappa-\frac{1}{2}\left( \psi_1+\psi_2+\psi_{p} \right)+\delta_{0, \{1,2\}}=\frac{1}{2}\pi_1^\ast \pi_2^\ast\kappa+\frac{1}{2}\pi_{p}^\ast\kappa, \]
    which is nef.

    \textbf{Step 2. } $D$ is F-nef.

    Since the intersection of $\kappa$ with any F-curve is $1$, and the coefficients of the boundary divisors are $\frac{1}{2}$, it is straightforward to verify that they intersect F-curves of types $1$ to $5$ non-negatively by \cref{thm:explicit}. For a type~6 F-curve $F_6^{g_1, g_2, g_3, g_4}(I_1, I_2, I_3, I_4)$, the result depends on the number of $\delta_{1, \emptyset}$ and $\delta_{2, \emptyset}$ among the $\delta_{g_i, I_i}$’s. If there are at most two such divisors, then the intersection is trivially non-negative. The remaining cases are $F_6^{1,1,1,1}(\emptyset,\emptyset,\emptyset,\emptyset)$ and $F_6^{0,1,1,2}([n],\emptyset,\emptyset,\emptyset)$, where a direct computation shows that the intersection is non-negative.

    \textbf{Step 3. } $D$ is nef.

    By Step~2 and \cite[Theorem 0.3]{GKM02}, it is enough to prove that $F^\ast D$ is nef. By \cite[Lemma 1]{AC09}, 
    \[ F^\ast D=\kappa-\frac{1}{2}\sum_{i=1}^4 \psi_i+\frac{1}{2}\sum_{\{i,j\}\subset [4]}\delta_{0, \{i,j\} }, \]
    which is nef by Step~1.
\end{proof}

\begin{prop}\label{prop:4n2nef}
    \[ D:=\kappa-\frac{1}{2}\psi_1+\frac{1}{2}\delta_{2, \emptyset}+\frac{1}{2}\delta_{4, \emptyset}+\frac{1}{2}\delta_{1, \{1\}} \]
    is nef on $\M{4}{n}$ for $n\ge 2$.
\end{prop}

\begin{proof}
    Since 
    \[ \kappa-\frac{1}{2}\psi_1+\frac{1}{2}\delta_{2, \emptyset}+\frac{1}{2}\delta_{4, \emptyset}=\frac{1}{2}\pi_1^\ast \kappa+\frac{1}{2}\left( \kappa+\delta_{2, \emptyset}+\delta_{4, \emptyset} \right), \]
    which is nef by \cref{eg:kapeg}~(3), it suffices to show that $D|_{\Delta_{1, \{1\}}}$ is nef. Let
    \[ \theta_1:\M{1}{\{1,p\}} \to \M{4}{n}, \quad \theta_2:\M{3}{\{1\}^c\cup p}\to \M{4}{n} \]
    be the relevant attaching maps. By \cref{lem:bdry}, it is enough to prove that $\theta_i^\ast D$ is nef for $i=1,2$. We have
    \[ \theta_1^\ast D = \kappa-\frac{1}{2}\psi_1-\frac{1}{2}\psi_p=\frac{1}{2}\left( \pi_1^\ast \kappa+\pi_p^\ast \kappa \right) \]
    and
    \[ \theta_2^\ast D=\kappa+\frac{1}{2}\delta_{2, \emptyset}-\frac{1}{2}\psi_p=\frac{1}{2}\left(\pi_p^\ast \kappa + \kappa+\delta_{2, \emptyset} \right), \]
    since $n\ge 2$. Both are nef by \cref{eg:kapeg}.
\end{proof}

For future reference, we record the following.

\begin{prop}\label{prop:charpsemi}
    The F-conjecture holds for $\M{1}{2}$, $\M{2}{1}$, $\M{2}{2}$ and $\M{3}{1}$. Moreover, if the base field has positive characteristic, then any nef line bundle on $\M{1}{2}$ and $\M{2}{1}$ is semiample.
\end{prop}

Note that the first part follows from \cite{GKM02} when the base field has characteristic not equal to $2$.

\begin{proof}
    On $\M{1}{2}$ (resp. $\M{2}{1}$), the cone of F-nef divisors is generated by $\lambda, 12\lambda - \delta_{\text{irr}}$ (resp. $\lambda, 12\lambda - \delta_{\text{irr}}, \psi_1$). These are nef regardless of the characteristic.

    Consider $\M{2}{2}$. The extremal rays of the cone of F-nef divisors are given in \cref{eg:M22}. The divisors in the first line are nef by the same argument as above. Since
    \[
    \psi_1 + \psi_2 - 2\delta_{0, \{1,2\}} = \pi_1^\ast \psi + \pi_2^\ast \psi,
    \]
    the second line is a sum of nef divisors and $\delta_{1,\{1\}}$. Hence it is enough to prove this is nef on $\Delta_{1, \{1\}}$. Since the F-conjecture holds for $\M{1}{2}$, F-nefness implies nefness on $\Delta_{1, \{1\}}$ by \cref{lem:bdry}. 
    
    For
    \[
    D = \lambda + \psi_1 + \psi_2 + \delta_{1, \emptyset},
    \]
    it is again enough, by \cref{lem:bdry}, to show that $\theta^\ast D$ is nef, where $\theta: \M{1}{3} \to \M{2}{2}$. We compute
    \[
    \theta^\ast D = \lambda + \psi_1 + \psi_2 - \psi_3 + \delta_{1, \emptyset}.
    \]
    Note that $\psi_1 = \psi_2$ on $\M{1}{2}$. Therefore,
    \[
    \psi_2 - \delta_{0, \{1,2\}} = \pi_1^\ast \psi_2 = \pi_1^\ast \psi_3 = \psi_3 - \delta_{0, \{1,3\}},
    \]
    so
    \[
    \theta^\ast D = \lambda + \psi_1 + \delta_{0, \{1,2\}} - \delta_{0, \{1,3\}} + \delta_{1, \emptyset} = \lambda + \pi_3^\ast \psi_1 + \delta_{0, \{1,2\}} + \delta_{1, \emptyset}.
    \]
    Since $\lambda + \pi_3^\ast \psi_1$ is nef and the F-conjecture holds for $\M{1}{2}$, this divisor is nef by \cref{lem:bdry}.

    Now consider $\M{3}{1}$. Note that \cite{Fab90} proves the F-conjecture for $\Mg{3}$. In particular,
    \[
    10\lambda - \delta_{\text{irr}} - 2\delta_{1, \emptyset} - 2\delta_{1, \{1\}} = \pi^\ast(10\lambda - \delta_{\text{irr}} - 2\delta_1)
    \]
    is nef. Hence, the first four divisors in \cref{eg:M31} are nef, and it remains to show that the following two divisors are also nef:
    \[
    10\lambda - \delta_{\text{irr}} + 2\psi_1 - 2\delta_{1, \emptyset}, \quad 11\lambda - \delta_{\text{irr}} + 3\psi_1 - \delta_{1, \emptyset} - 2\delta_{1, \{1\}}.
    \]
    We compute:
    \[
    10\lambda - \delta_{\text{irr}} + 2\psi_1 - 2\delta_{1, \emptyset} = (10\lambda - \delta_{\text{irr}} - 2\delta_{1, \emptyset} - 2\delta_{1, \{1\}}) + 2\psi_1 + 2\delta_{1, \{1\}},
    \]
    so it suffices to prove that its restriction to $\Delta_{1, \{1\}}$ is nef. This follows from \cref{lem:bdry} and the F-conjecture for $\M{1}{2}$ and $\M{2}{1}$.
    
    Moreover, we have:
    \[
    11\lambda - \delta_{\text{irr}} + 3\psi_1 - \delta_{1, \emptyset} - 2\delta_{1, \{1\}} = \lambda + (10\lambda - \delta_{\text{irr}} - 2\delta_{1, \emptyset} - 2\delta_{1, \{1\}}) + 3\psi_1 + \delta_{1, \emptyset}.
    \]
    Again, this follows from \cref{lem:bdry} and the F-conjecture for $\M{1}{1}$ and $\M{2}{2}$.

    We now consider the second assertion. Since $\lambda$ corresponds to the Satake compactification, it is always semiample on $\Mg{g}$ for $g \ge 2$ and on $\M{1}{1}$. Also, by \cite{SB22}, the divisor $12\lambda - \delta_{\mathrm{irr}}$ is semiample on $\Mg{g}$ for $g \ge 2$. Note that $\lambda$ is zero on genus $0$, and $12\lambda - \delta_{\mathrm{irr}}$ is zero on genus $0$ and $1$. Since $\lambda$ and $12\lambda - \delta_{\mathrm{irr}}$ on $\M{g}{n}$ are pullbacks of those on $\Mg{g}$ or $\M{1}{1}$, they are semiample in general. Therefore, the second assertion follows from this and \cite{Ke99}.
\end{proof}

%% file: sections/Knudsen.tex
The main purpose of this section is to prove the exact sequence in \cref{thm:Knuses}, which will be used in \cref{sec:Knudsen}. As before, we refer to the beginning of \cref{subsec:method} for the motivation behind this exact sequence.

The \textbf{Knudsen's construction} is the map
\[
f^\K_{g,n}:\M{g}{n}\to \M{g}{n-1}\times_{\M{g}{n-2}}\M{g}{n-1},
\]
defined in \cite{Knu83}. For any $g\ge 0$ and $n\ge 2$ such that $2g-2+n>0$, we define $F^{\K}_{g,n}$ to be
\[
\left\{F_6^{0,0,g_1,g_2}( \left\{n-1\right\}, \left\{n\right\}, I_1, I_2)\ |\ g_1+g_2=g, I_i\ne\emptyset\text{ if }g_i=0 \right\}\cup \left\{F_5^{0,0}(\left\{n-1\right\}, \left\{n\right\}) \right\}.
\]
This is exactly the set of F-curves contracted by $f^\K_{g,n}$. 

\begin{thm}\label{thm:Knuses}
    There exists an exact sequence
    \[
    0\to \Pic(\M{g}{n-2})\xrightarrow{(\pi_n^\ast, -\pi_{n-1}^\ast)} \Pic(\M{g}{n-1})\times \Pic(\M{g}{n-1})\xrightarrow{\pi_n^\ast+\pi_{n-1}^\ast}\Pic(\M{g}{n})\to \mathbb{Q}^{F^{\K}_{g,n}}\to 0.
    \]
    Moreover, the image of $\pi_n^\ast+\pi_{n-1}^\ast$ is the same as the image of $f^{\K\ast}_{g,n}$.
\end{thm}

\begin{rmk}
    In general, $f^\K_{g,n}$ is not a contraction and admits a Stein factorization
    \[
    \M{g}{n} \to X_{g,n}^{\K} \to \M{g}{n-1} \times_{\M{g}{n-2}} \M{g}{n-1},
    \]
    where $X_{g,n}^{\K}$ is the coarse moduli space of $\MM{g}{n-1} \times_{\MM{g}{n-2}} \MM{g}{n-1}$, as in \cite[Theorem 4.5]{Cho25}. As in \cite{Cho24}, \cref{thm:Knuses} characterizes the Picard group of $X_{g,n}^{\K}$ as the image of $\pi_n^\ast + \pi_{n-1}^\ast$.
\end{rmk}

The following theorem is the key step of \cref{thm:Knuses}. Indeed, a much weaker version of this theorem would suffice, but we will prove a stronger version in order to establish \cref{cor:simp}.

\begin{thm}\label{thm:Knu}(char $k\ne 2$)
    For any $C\in F^{\K}_{g,n}$, there exists a nef divisor $D_C$ such that $D_C\cdot C>0$ and $D_C\cdot C'=0$ for every $C'\in F^{\K}_{g,n}\setminus \left\{C\right\}$.
\end{thm}

\begin{proof}
    If $g=0$, then this is proved in \cite[Corollary 3.15]{Cho24}. Assume $g\ge 1$. If $C$ is the type 5 curve in $F^{\K}_{g,n}$, then
    \[
    D = \pi_{\left\{n-1,n\right\}}^\ast \left( \kappa + \sum_{1\le i \le g} \delta_{i, \emptyset} \right)
    \]
    works, where $\pi_{\left\{n-1,n\right\}} : \M{g}{n} \to \M{g}{2}$. Note that $D$ is nef by \cref{eg:kapeg}, and $D_C \cdot C = 2$.

    Now assume that $C$ is type 6. Then $C$ is the image of a curve $C_0$ in $F^{\K}_{0,g+n}$ under the flag map $F : \M{0}{g+n} \to \M{g}{n}$. By \cite[Corollary 3.15]{Cho24}, there exists a nef divisor on $\M{0}{g+n}$ which intersects $C_0$ nontrivially and contracts all other curves in $F^{\K}_{0,g+n}$. By symmetrizing this with the $S_g$-action and using \cite[Theorem 0.7]{GKM02}, there exists a nef divisor $D_0$ on $\M{g}{n}$ such that $D_0 \cdot C \ne 0$ and contracts all other type 6 curves in $F^{\K}_{g,n}$.

    Let $D_{\text{GKM}}$ be a divisor from \cite[Section 4]{GKM02}. This is a nef line bundle on $\M{g}{n}$ such that $F^\ast D_{\text{GKM}} = 0$ and it intersects nontrivially with all F-curves of types 1–5. Choose a constant $c$ such that $D_0 - c D_{\text{GKM}}$ contracts the type 5 curve in $F^{\K}_{g,n}$.

    Now, choose a sufficiently ample divisor $A$ on $\M{g}{n-1}$ and let
    \[
    D := D_0 - c D_{\text{GKM}} + \pi_{n-1}^\ast A + \pi_n^\ast A.
    \]
    Note that $D$ satisfies all the conditions of the statement except possibly nefness. Since F-curves not in $F^{\K}_{g,n}$ are not contracted by both $\pi_{n-1}$ and $\pi_n$, this is F-nef for sufficiently ample divisor $A$. Moreover, since $F^\ast D_{\text{GKM}} = 0$, we have $F^\ast D$ nef. Therefore, by \cite[Theorem 0.7]{GKM02}, $D$ is nef.
\end{proof}

\begin{cor}\label{cor:Knuind}
    $F^{\K}_{g,n}$ is linearly independent. 
\end{cor}

\begin{proof}
If $\operatorname{char} k \ne 2$, then this follows from \cref{thm:Knu}. This also implies the characteristic $2$ case, since the Picard group of the moduli space of curves does not depend on the characteristic.
\end{proof}

The following lemma is well known. See, for example, \cite{AC09}.

\begin{lem}\label{lem:picnum}
    \begin{enumerate}
        \item $\rho(\M{0}{n})=2^{n-1}-\binom{n}{2}-1$.
        \item $\rho(\M{1}{n})=2^{n}-n$.
        \item $\rho(\overline{\rm{M}}_2)=2$ and $\rho(\M{2}{n})=3\cdot 2^{n-1}$ for $n\ge 1$.
        \item If $g\ge 3$, then $\rho(\overline{\rm{M}}_g)=\lfloor\frac{g}{2}\rfloor+2$ and $\rho(\M{g}{n})=(g+1)\cdot 2^{n-1}+1$ for $n\ge 1$. 
    \end{enumerate}
\end{lem}

\begin{proof}[Proof of \cref{thm:Knuses}]
    The case of $g=0$ (resp. $g=1$) is proved in \cite[Theorem 3.17]{Cho24} (resp. \cite[Corollary 4.9]{Cho25}). Hence, we will assume that $g\ge 2$. The surjectivity of the last part follows from \cref{cor:Knuind}.

    Now, we will prove that the sequence 
    \[
    \Pic(\M{g}{n-1})\times \Pic(\M{g}{n-1})\xrightarrow{\pi_n^\ast+\pi_{n-1}^\ast}\Pic(\M{g}{n})\to \mathbb{Q}^{F^{\K}_{g,n}}
    \]
    is exact. Since the last map is surjective, it suffices to prove that there exists a set $T$ consisting of $|F^{\K}_{g,n}|$ divisors such that $T\cup \text{Im }\pi_{n-1}^\ast\cup \text{Im }\pi_{n}^\ast$ spans $\Pic(\M{g}{n})$. Note that, by \cite[Lemma 1 (i), (ii), (iii)]{AC09}, it is enough to show that every boundary divisor is contained in the span of $T\cup \text{Im }\pi_{n-1}^\ast\cup \text{Im }\pi_{n}^\ast$. Let $V_T$ be the span of $T\cup \text{Im }\pi_{n-1}^\ast\cup \text{Im }\pi_{n}^\ast$.

    First, consider the case of $n=2$. In this case, $|F^{\K}_{g,n}|=\lfloor\frac{g}{2}\rfloor+1$. Define
    \[ T:= \left\{ \delta_{i,\left\{1,2\right\}} \ |\ 0\le i\le \lfloor\frac{g}{2}\rfloor \right\}. \]
    Then, by \cite[Lemma 1 (iv)]{AC09}, it is straightforward to see that $V_T$ contains all of the boundary divisors, so $V_T=\Pic(\M{g}{n})$.

    Next, assume $n\ge 3$. In this case, $|F^{\K}_{g,n}|=(g+1)\cdot 2^{n-3}$. Define
    \[ T:= \left\{ \delta_{i,I} \ |\ I\subseteq [n-1],\ n-1,n-2\in I \right\}\cup \left\{\psi_n\right\} \]
    Then, it is straightforward to see that $|T|=(g+1)\cdot 2^{n-3}$. First, consider $\delta_{i,I}$ for $I\subseteq [n]$, such that $n-1, n\in I$. If $(i, I)=(0, \left\{n-1, n\right\})$, then $\delta_{i,I}=\psi_n-\pi_{n-1}^\ast \psi_n$, so $\delta_{i,I}\in V_T$. Assume $(i, I)\ne (0, \left\{n-1, n\right\})$. If $n-2\not\in I$, then $\delta_{i, I\setminus \left\{n-1 \right\}}\in T$, so
    \[ \delta_{i, I}=\pi_{n-1}^\ast \delta_{i, I\setminus \left\{n-1 \right\}}-\delta_{i, I\setminus \left\{n-1 \right\}}\in V_T \]
    and if $n-2\in I$, then $\delta_{i, I\setminus \left\{n \right\}}\in T$, so
    \[ \delta_{i, I}=\pi_{n}^\ast \delta_{i, I\setminus \left\{n \right\}}-\delta_{i, I\setminus \left\{n \right\}}\in V_T. \]
    Hence, $\delta_{i, I}\in V_T$ for all such $I$. Now, again by \cite[Lemma 1 (iv)]{AC09}, it is straightforward that this implies $\delta_{i, I}\in V_T$ for any $(i,I)$. Therefore, $V_T=\Pic(\M{g}{n})$. This completes the proof of the exactness of $\Pic(\M{g}{n-1})\times \Pic(\M{g}{n-1})\xrightarrow{\pi_n^\ast+\pi_{n-1}^\ast}\Pic(\M{g}{n})\to \mathbb{Q}^{F^{\K}_{g,n}}$

    Now, consider the exactness of the full sequence. It is well known that $\Pic(\M{g}{n-2})\xrightarrow{(\pi_n^\ast, -\pi_{n-1}^\ast)} \Pic(\M{g}{n-1})\times \Pic(\M{g}{n-1})$ is injective. Hence, the theorem follows from the following elementary computation:
    \[  |F^{\K}_{g,n}| = \rho(\M{g}{n-2})-1 = \rho(\M{g}{n})-\rho(\M{g}{n-1})-\rho(\M{g}{n-1})+\rho(\M{g}{n-2}) \]
    for $g\ge 3$, and
    \[ |F^{\K}_{2,n}| = \rho(\M{2}{n-2})= \rho(\M{2}{n})-\rho(\M{2}{n-1})-\rho(\M{2}{n-1})+\rho(\M{2}{n-2}) \]
    which follows directly from \cref{lem:picnum}.
\end{proof}

We observe that this establishes the extremality of F-curves in $F^{\K}_{g,n}$. The following result follows directly from \cref{thm:Knuses} and \cref{thm:Knu}.

\begin{cor}[char $k\ne 2$]\label{cor:simp}
    $\NE{f^\K_{g,n}}$ is a simplicial cone generated by $F^{\K}_{g,n}$. Moreover, any element of $F^{\K}_{g,n}$ is a regular extremal ray.
\end{cor}

%% file: sections/strategy.tex
Here, we describe the method for proving the regular extremality of F-curves, using \cref{sec:nef} and \cref{sec:Knudsen}. Let $F$ be a curve whose regular extremality we wish to prove. Our goal is to construct a set of nef divisors $T$ that contract $F$ and span a codimension $1$ subspace of the Picard group. This is precisely what was done in the proof of \cref{cor:ext01}. What we proved there is essentially the following: if we let
\[
T := \{ \pi_S^\ast \psi_i\ |\  \pi_S^\ast \psi_i\cdot F=0 \},
\]
then $T$ spans a codimension $1$ subspace. We aim to apply the same strategy for curves on $\M{g}{n}$ in general. However, there are two main obstacles to this method:

\begin{enumerate}
    \item We do not have enough nef divisors. Pullbacks of $\psi$-classes alone are insufficient.
    \item Even with enough nef divisors, it is difficult to compute the dimension of their span.
\end{enumerate}

\cref{sec:nef} addresses issue (1); we will soon see that the divisors constructed there are well-suited for contracting F-curves. Issue (2) will be resolved by an inductive argument using \cref{thm:Knuses}. Instead of proving the regular extremality of a single curve $F$, we prove the regular extremality of a family of curves on $\M{g}{n}$—with fixed $g$ and varying $n$—that is stable under the projection maps, to use the induction. Then, the regular extremality of $\pi_{\ast} F$ and \cref{thm:Knuses} reduce the problem to computing the intersection of divisors contracting $F$ with $F_{g,n}^{\text{Knu}}$. This is formalized by the following lemma.

\begin{lem}\label{lem:extind}
    Fix $g,n\in \N$. Let $C_{n+2}\subseteq \NE{\M{g}{n+2}}$ be a subset, $C_{n+1}^i:=\pi_{n+i, \ast}C_{n+2}$ for $i=1,2$, and $C_n:=\pi_{n+2, \ast}C_{n+1}^1=\pi_{n+1, \ast}C_{n+1}^2$. Let $N_{n+2}, N_{n+1}^1, N_{n+1}^2, N_n$ denote the set of nef divisors  that intersect  $C_{n+2}, C_{n+1}^1, C_{n+1}^2, C_n$ trivially. Assume the following conditions hold:
    \begin{enumerate}
        \item $I(C_{n}), I(C_{n+1}^1)$ and $I(C_{n+1}^2)$ coincide with the dimension of the subspace spanned by $C_n, C_{n+1}^1$, and $C_{n+1}^2$.
        \item The image of the intersection pairing $N_{n+2}\to \R^{F^{\K}_{g,n+2}}$ spans $\R^{F^{\K}_{g,n+2}}$.
    \end{enumerate}
    Then $I(C_{n+2})\le I(C_{n+1}^1)+I(C_{n+1}^2)-I(C_n)$. If, moreover, $I(C_{n+1}^2)=I(C_n)$, then $I(C_{n+2})=I(C_{n+1}^1)$ and this coincides with the dimension of the subspace spanned by $C_{n+2}$. 
\end{lem}

\begin{proof}
    This is a direct consequence of \cref{thm:Knuses}. We have an exact sequence
    \[
    0\to \Pic(\M{g}{n})_{\R}\to \Pic(\M{g}{n+1})_\R\times \Pic(\M{g}{n+1})_\R\to\Pic(\M{g}{n+2})_\R\to \mathbb{R}^{F^{\K}_{g,n+2}}\to 0.
    \]
    We first compute the dimension of the image of 
    \[
    (N_{n+1}^1\otimes \R)\times (N_{n+1}^2\otimes \R)\to N_{n+2}\otimes \R \subseteq \Pic(\M{g}{n+2})_\R.
    \]
    By the condition of (1), $N_n\otimes \R$ coincides with the subspace of $\R$-divisors that intersect $C_n$ trivially. Therefore, by the exact sequence above, the kernel of this map is $N_n\otimes \R$. Hence, the dimension of the image is
    \[
    \dim N_{n+1}^1 + \dim N_{n+1}^2 -\dim N_n=2\rho(\M{g}{n+1})-\rho(\M{g}{n})-(I(C_{n+1}^1)+I(C_{n+1}^2)-I(C_n)).
    \]
    By (2), there exists a $|F^{\K}_{g,n+2}|$-dimensional subspace of $N_{n+2}\otimes \R$ that is independent from the image of $(N_{n+1}^1\otimes \R)\times (N_{n+1}^2\otimes \R)$. Hence, the dimension of $N_{n+2}\otimes \R$ is at least
    \begin{align*}
         &2\rho(\M{g}{n+1})-\rho(\M{g}{n})+|F^{\K}_{g,n+2}|-(I(C_{n+1}^1)+I(C_{n+1}^2)-I(C_n))\\&=\rho(\M{g}{n+2})-(I(C_{n+1}^1)+I(C_{n+1}^2)-I(C_n)).
    \end{align*}
    Hence $I(C_{n+2})\le I(C_{n+1}^1)+I(C_{n+1}^2)-I(C_n)$. Moreover, since $\pi_{i\ast}$ on $\A_1$ is surjective, by (1), the dimension of the subspace spanned by $C_{n+2}$ is at least $I(C_{n+1}^1)$. Hence, the codimension of $N_{n+2}$ is at least $I(C_{n+1}^1)$. Therefore, $I(C_{n+1}^2)=I(C_n)$ implies $I(C_{n+2})=I(C_{n+1}^1)$. The last assertion is straightforward since the dimension of such a subspace is $\leq I(C_{n+2})$.
    
\end{proof}

Hence, while using the induction, the main part of the proof is to check the second condition. However, this remains a nontrivial task, since $|F_{g,n}^{\text{Knu}}|$ has roughly $g \cdot 2^{n-3}$ elements, which is still a large number when $n$ is large. Nevertheless, we have not yet fully exploited the inductive structure. In \cref{lem:extind}, we only used the projections $\pi_{n+1}$ and $\pi_{n+2}$, while ignoring the other $n$ projections. The second condition of \cref{lem:extind} can be made easier to verify by taking these into account. To facilitate this, we now introduce some notation.

For $n\geq 2$, we will denote
\[
V_{g,n}=
\begin{cases}
    \R^{F^{\K}_{g,n}} &\text{ if }n=2,\\
    \bigoplus_{0\leq j\leq \lfloor\frac{g}{2}\rfloor}\R e_j &\text{ if }g\text{ is odd or }n\text{ is even}\ne 2,\\
    \bigoplus_{0\leq j< \lfloor\frac{g}{2}\rfloor}\R e_j &\text{ otherwise}.
\end{cases}
\]
Also, we define $E: \R^{F^{\K}_{g,n}} \to V_{g,n}$ by the identity map when $n=2$ and
\[
E(h)=\sum_{j} \sum_{ F_6^{0,0,j,g-j}( \left\{n-1\right\}, \left\{n\right\}, I, I^c)  }(-1)^{|I|}h(F_6^{0,0,j,g-j}( \left\{n-1\right\}, \left\{n\right\}, I, I^c))  e_j,
\]
where we consider $h$ as a function from $F^{\K}_{g,n}$ to $\R$, if $n\ne 2$. Note that $E$ is defined as a direct sum of linear functions, each corresponding to $j$. The definition of $E$ may appear technical, but the underlying idea is simple: the kernel of the map
\[
E: \R^{F^{\K}_{g,n+2}} \to V_{g,n+2}
\]
is precisely the image of
\[
\bigoplus_{1\leq i\leq n} \pi_i^\ast: \bigoplus_{1\leq i\leq n} \R^{F^{\K}_{g,n+1}} \to \R^{F^{\K}_{g,n+2}}
\]
where $\pi_i^\ast$ is the dual of the pushforward $\pi_{i, \ast}:F^{\K}_{g,n+2}\to F^{\K}_{g,n+1}\cup \{0\}$. More precisely,

\[ \pi_i^\ast(h)(F)=\begin{cases}
    h(\pi_{i, \ast}(F))&\text{ if }\pi_{i, \ast}(F)\ne 0\\
    0 & \text{ otherwise}
\end{cases}. \]
Thus, the inductive strategy allows us to verify the statement after composing with $E$. This can be formalized as follows.

\begin{lem}\label{lem:second}
    Let $C_{n+2}\subseteq \NE{\M{g}{n+2}}$ be a subset, and define $C_{n+1}^i:=\pi_{i\ast}C_{n+2}$ for $1\leq i\leq n$. Let $N_{n+2}, N_{n+1}^i$ be defined as before. If 
    \begin{enumerate}
        \item For each $1\leq i\leq n$, the image of the intersection pairing $N_{n+1}^i\to \R^{F^{\K}_{g,n+1}}$ spans $\R^{F^{\K}_{g,n+1}}$.
        \item The image of the composition of the intersection pairing with $E$, 
        \[f: N_{n+2}\to  \R^{F^{\K}_{g,n+2}}\to V_{g, n+2},  \]
        spans $V_{g, n+2}$, 
    \end{enumerate}
    then the image of the intersection pairing $N_{n+2}\to \R^{F^{\K}_{g,n+2}}$ spans $\R^{F^{\K}_{g,n+2}}$.
\end{lem}

\begin{proof}
    Note that $\pi_i$ for $1\le i\le n$ maps $\R^{F^{\K}_{g,n+2}}$ to $\R^{F^{\K}_{g,n+1}}$. By (1), the span of the image of $N_{n+2}\to \R^{F^{\K}_{g,n+2}}$ contains the image of $\pi_i^\ast: \R^{F^{\K}_{g,n+1} }\to \R^{F^{\K}_{g,n+2}}$ for $1\le i\le n$. Therefore, it is enough to prove that
    \[ N_{n+2}\to  \R^{F^{\K}_{g,n+2}}\to \R^{F^{\K}_{g,n+2}}/\left(\oplus_{i=1}^n \pi_i^\ast\R^{F^{\K}_{g,n+1}}  \right)  \]
    is surjective, which reduces to (2) once we show that the kernel of $E$ is the image of $\oplus_{1\leq i\leq n}\pi_i^\ast$. 
    
    If $n=0$, there is nothing to prove. Assume $n>0$. Note that in the definition of $V_{g,n}$, as well as in the maps $E$ and $\pi_i^\ast$, every morphism respects the decomposition indexed by $j$, which parametrizes the genus occurring in type $6$ curves, and the type $5$ curve. Hence, it is enough to prove the statement within each component of the decomposition of $F^{\K}_{g,n+1}, F^{\K}_{g,n+2}$, and $V_{g,n}$.
    
    If $j\ne \frac{g}{2}$, then the statement follows from  
    \begin{align*}
        &\pi_{i\ast}^{-1}\left( F_6^{0,0,j,g-j}( \left\{n+1\right\}, \left\{n+2\right\}, I, J) \right)=\\&\left\{ F_6^{0,0,j,g-j}( \left\{n+1\right\}, \left\{n+2\right\}, I\cup \left\{i\right\}, J), F_6^{0,0,j,g-j}( \left\{n+1\right\}, \left\{n+2\right\},I, J\cup \left\{i\right\}) \right\}
    \end{align*}
    where $\pi_{i\ast}:F^{\K}_{g,n+2}\to F^{\K}_{g,n+1}$. The statement for the type 5 curve is automatic
    \[ \pi_{i\ast}^{-1}\left( F_5^{0,0}( \left\{n+1\right\}, \left\{n+2\right\}) \right)=\left\{  F_5^{0,0}( \left\{n+1\right\}, \left\{n+2\right\}) \right\}. \]

    Now assume $g = 2j$. Then, there is one additional relation:
    \[
    F_6^{0,0,j,j}( \left\{n+1\right\}, \left\{n+2\right\}, I, J) = F_6^{0,0,j,j}( \left\{n+1\right\}, \left\{n+2\right\}, J, I).
    \]
    If $n$ is even, then this relation, together with the formula for 
    \[
    \pi_{i\ast}^{-1}\left( F_6^{0,0,j,g-j}( \left\{n+1\right\}, \left\{n+2\right\}, I, J) \right),
    \]
    shows that the image of $\bigoplus_{1\leq i\leq n} \pi_i^\ast$ is surjective in this component, which completes the proof. If $n$ is odd, then this new relation lies in the image of $\bigoplus_{1\leq i\leq n} \pi_i^\ast$, so the situation is the same as in the case $g \ne 2j$. Hence, the statement also holds by the same argument.

\end{proof}

Note that the dimension of $V_{g,n}$ is $\lfloor \frac{g}{2} \rfloor + 1$, regardless of the value of $n$. This illustrates how significantly \cref{lem:second} simplifies the entire inductive process.

We will implicitly apply the following useful computation while proving theorems.

\begin{lem}\label{lem:comp}
    Let $f:\Pic(\M{g}{n})\to V_{g,n}$ be the natural extension of the map $f$ in \cref{lem:second}, i.e., the composition of the pairing map $\Pic(\M{g}{n}) \to \R^{F^{\K}_{g,n}}$ and $E:\R^{F^{\K}_{g,n+2}}\to V_{g,n}$. Then
    \[
    f(\kappa) = f(\psi_i)= -e_0, \quad\text{and}\quad f(\delta_{i, I}) = (-1)^{|I|+1} e_i \quad\text{for}\quad 0 \le i \le \left\lfloor \frac{g}{2} \right\rfloor.
    \]
\end{lem}

\begin{proof}
    If $i \ne 1$, then the coefficient of $e_i$ in $f(\kappa)$ is $0$, since the intersection of any F-curve with $\kappa$ is $1$ and
    \[
    \sum_{I\subseteq [n]} (-1)^{|I|} = 0.
    \]
    The coefficient is $-1$ for $e_0$ because there is no $F_6^{0,0,0,g}(n-1, n, \emptyset, [n-2])$.The case of \( \psi_i \) is essentially the same as that of \( \kappa \).  The last assertion follows directly from \cref{thm:explicit} by a case-by-case analysis with respect to $|I \cap \{n-1, n\}|$.
\end{proof}

Note that if $g$ is even and $n$ is odd, then there is no $e_{\lfloor\frac{g}{2}\rfloor}$. However, for convenience, we will extend the notation by including $e_{\lfloor\frac{g}{2}\rfloor}$ with the convention $e_{\lfloor\frac{g}{2}\rfloor} = 0$. \cref{lem:comp} still holds under this convention.

%% file: sections/extremality.tex
In this section, we assume that the characteristic of the base field is not equal to~2, unless stated otherwise. This section is devoted to prove \cref{thm:main2} and \cref{thm:main3} using \cref{sec:nef} to \cref{subsec:method}. Before we begin, as promised, we provide a detailed statement of \cref{thm:main2} and \cref{thm:main3}.

\begin{thm}\label{thm:main2det}
    \begin{enumerate}
        \item Types 1 and 4 F-curves on $\M{g}{n}$ span regular extremal rays of $\NE{\M{g}{n}}$, while type 2 F-curves do not.
        \item Apart from the following three exceptional cases, no type 3 $F$-curve spans an extremal ray of $\NE{\M{g}{n}}$.
        \begin{enumerate}
            \item $F_3^0([n])$ on $\M{2}{n}$ for $n\geq 1$,
            \item $F_3^1([n])$ on $\M{3}{n}$ for $n\geq 0$,
            \item $F_3^2(\emptyset)$ on $\overline{\rm{M}}_4$.
        \end{enumerate}
        In the three exceptional cases, the corresponding curve spans a regular extremal ray.
        \item Knudsen-type F-curves (cf. \cref{defn:Knudsen}) are regular extremal. Consequently, every F-curve on $\M{0}{n}$ and $\M{1}{n}$ is regular extremal. 
        \item The following type 6 F-curves are regular extremal:
        %A list of regular extremal type 6 
        \begin{enumerate}
            \item $F_6^{1,1,1,g-3}(\emptyset, \emptyset, \emptyset, [n])$,
            \item $F_6^{1,1,2,g-4}(\emptyset, \emptyset, \emptyset, [n]) \text{ for }n\ne 0$,
            \item $ F_6^{0,1,1,g-2}(i, \emptyset, \emptyset, [n]\setminus i)$ for all $i\in [n]$.
        \end{enumerate}
    \end{enumerate}
\end{thm}

\begin{thm}\label{thm:main3det}
    \begin{enumerate}
        \item Any type 5 or type 6 F-curve on $\M{2}{n}$ is regular extremal.
        \item Except $F_5^{0,1}(I, \emptyset)$, every type 5 or type 6 F-curve on $\M{3}{n}$ is regular extremal.
        \item Except for the following cases, every type 5 or type 6 F-curve on $\M{4}{n}$ is regular extremal:
        \begin{enumerate}
            \item $F_5^{i, 1}(I, \emptyset)$, $i=0,1$.
            \item $F_5^{2,1}([n],\emptyset)$ and $F_5^{0,2}([n],\emptyset)$.
        \end{enumerate}
    \end{enumerate}
\end{thm}

\subsection{Type 1,2 and 4}\label{subsec:124}

In this case, the proof of (non-)extremality is relatively straightforward and follows as a direct application of \cite{GKM02}.

\begin{thm}\label{thm:type1}
    The type $1$ F-curve on $\M{g}{n}$, $g\ge 1$, is regular extremal.
\end{thm}

\begin{proof}
    Let $L_1,\dots, L_\rho$ be a basis of $\Pic(\M{g}{n})$ consisting of nef line bundles. Consider the set of line bundles
    \[
    T = \left\{ L_i-\frac{L_i\cdot [F_1]}{12}\lambda\ \middle|\ 1\leq i\leq \rho \right\}.
    \]
    Since $\lambda\cdot F_1=12$, each line bundle $L_i-\frac{L_i\cdot [F_1]}{12}\lambda$ intersects $F_1$ trivially. Moreover, since $\lambda$ intersects all other $F$-curves trivially, it follows that $L_i-\frac{L_i\cdot [F_1]}{12}\lambda$ is $F$-nef. For the flag map $F:\M{0}{g+n}\to \M{g}{n}$, the pullback of $L_i-\frac{L_i\cdot [F_1]}{12}\lambda$ along $F$ coincides with $F^\ast L_i$, as $\lambda=0$ on $\M{0}{g+n}$. In particular, this implies nefness. Hence, by \cite[Theorem 0.3]{GKM02}, each $L_i-\frac{L_i\cdot [F_1]}{12}\lambda$ is nef. Thus, the set $T$ consists of line bundles that intersect $[F_1]$ trivially. Since $T$ together with $\lambda$ spans the Picard group, $T$ generates a codimension-$1$ subspace. Consequently, $I([F_1])=1$, implying that $[F_1]$ is regular extremal.
\end{proof}

\begin{thm}\label{thm:type2}
    The type $2$ $F$-curve on $\M{g}{n}$, for $g\geq 3$, is not extremal.
\end{thm}

\begin{proof}
    By \cite[Theorem 2.1]{GKM02},
    \[
    [F_2]=\frac{1}{2}\left( [F_3^{1}( \emptyset)]+[F_4^{1}( \emptyset)] \right).
    \]
    Since $[F_2]$ is expressed as a positive linear combination of other effective curve classes, it is not extremal.
\end{proof}

\begin{thm}\label{thm:type4}
    Any type $4$ F-curve on $\M{g}{n}$, $g\ge 2$, is regular extremal.
\end{thm}

\begin{proof}
    Let $S$ be the set of all type $4$ $F$-curves on $\M{g}{n}$. Choose any curve $F\in S$. Since $\kappa$ is an ample divisor whose intersection with any $F$-curve is $1$, it follows from \cite[Theorem 2.1]{GKM02} that there exist perturbations of $\kappa$, denoted $L_1,\dots, L_\rho$, such that:
    \begin{enumerate}
        \item $L_1,\dots, L_\rho$ form a basis of $\Pic(\M{g}{n})$ consisting of nef line bundles.
        \item For any $1\leq i\leq \rho$, we have $L_i\cdot F=\min_{F'\in S} L_i\cdot F'$.
    \end{enumerate}
    Consider the set of line bundles
    \[
    T = \left\{ L_i-\frac{L_i\cdot F}{2}(12\lambda-\delta_{\text{irr}})\ \middle|\ 1\leq i\leq \rho \right\}.
    \]
    By the same reasoning as in \cref{thm:type1}, the elements of $T$ are nef divisors that intersect $F$ trivially and span a codimension-$1$ subspace of the Picard group. Consequently, $I(F)=1$, which implies that $F$ is regular extremal.

\end{proof}

This completes the proof of \cref{thm:main2} (1).

\subsection{Type 3}\label{subsec:3}

Here, we prove \cref{thm:main2det} (2).

\begin{thm}\label{thm:type3next}
   Except for the three cases described in \cref{thm:main2det} (2), any type 3 F-curve does not form an extremal ray.
\end{thm}

\begin{proof}
    If $F_3^{i}( I)\ne F_3^{g-2}( [n])$, then by \cref{thm:explicit},
    \[
    [F_3^{i}(I)]=\frac{1}{2}\left([F_5^{1,i}(\emptyset, I)] + [F_5^{1,g-i-1}(\emptyset, I^c)] \right),
    \]
    and since $[F_5^{1,i}(\emptyset, I)]\ne [F_5^{1,g-i-1}(\emptyset, I^c)]$, it follows that $[F_3^{i}(I)]$ is not an extremal ray.
    
    If $F_3^{i}(I) = F_3^{g-2}( [n])$ with $g\geq 4$ and $(g,n)\ne (4,0)$, then
    \[
    [F_3^{g-2}([n])]=\frac{1}{2}\left([F_5^{1,2}(\emptyset, \emptyset)] + [F_5^{g-3,2}([n],\emptyset)] \right),
    \]
    and since $[F_5^{1,2}(\emptyset, \emptyset)] \ne [F_5^{g-3,2}([n],\emptyset)]$, it follows that $[F_3^{g-2}([n])]$ is not extremal.
\end{proof}

Now, we need to establish the regular extremality of type 3 F-curves for the three cases in \cref{thm:main1}~(2). The last case follows from \cref{eg:mext}, provided that the characteristic is not equal to \(2\). In the case of general characteristic, then this can be seen explicitly, since \( \lambda \), \( 12\lambda - \delta_{\text{irr}} \), and \( \kappa + \delta_{2,\emptyset} \) all contract the curve \( F_3^2(\emptyset) \) and nef by \cref{thm:kapnef}.

\begin{thm}[All characteristics]\label{thm:type3a}
    $F_3^{0}([n])$ is a regular extremal on $\NE{\M{2}{n}}$ for $n\geq 1$
\end{thm}

\begin{proof}
    We will prove that $I(F_3^{0}([n])) = 1$, which is equivalent to the given statement. For $n = 1$, this holds since both $\lambda$ and $12\lambda - \delta_{\text{irr}}$ contract $F_3^{0}([1])$. The case $n = 2$ follows from \cref{eg:M22} and \cref{prop:charpsemi}.
    
    We will first verify that $F_3^{0}([n])$ satisfies the conditions of \cref{lem:second}. Let $N_n$ be the set of nef divisors that intersect $F_3^{0}([n])$ trivially. Note that the image of $F_3^{0}([n])$ under any projection is $F_3^{0}([n-1])$. Hence, it suffices to show that the image of $f:N_n\to V_{2,n}$ spans $V_{2,n}$ for each $n$. 
    
    Define the divisor
    \[
    D_n=\kappa+\delta_{1,\emptyset}+\delta_{2,\emptyset}
    \]
    on $\M{2}{n}$ for $n\geq 2$. This divisor is nef by \cref{thm:kapnef}, see also \cref{eg:kapeg} (3). Furthermore, by \cite[Theorem 2.1]{GKM02}, we have $D_n\cdot F_3^{0}([n])=0$. For $n=2$, we compute
    \[
    f(D_2)=2\cdot F_5, \quad f(\psi_1)=F_5+F_6.
    \]
    Thus, $f$ satisfies the assertion. For $n>2$ and $n$ odd,
    \[
    f(\psi_1)=-e_0.
    \]
    For $n>2$ and $n$ even,
    \[
    f(D_n)=-2e_0-e_1, \quad f(\psi_1)=-e_0.
    \]
    Hence, in all cases, $f$ satisfies the assertion. Therefore, by \cref{lem:second}, for any $n\geq 2$, the image of $N_n\to \R^{F^{\K}_{g,n}}$ generates $\R^{F^{\K}_{g,n}}$.
    
    Now, we apply \cref{lem:extind} with $C_n=[F_3^{0}([n])]$. We have already verified condition (1) for $n=1,2$. Moreover, condition (2) holds for every $n$ by the preceding paragraph. Hence, by induction using \cref{lem:extind}, we conclude that $I(F_3^{0}([n]))=1$ for every $n\geq 1$.
\end{proof}

\begin{thm}[All characteristics]\label{thm:type3b}
    $F_3^{1}( [n])$ is regular extremal on $\NE{\M{3}{n}}$ for $n\geq 0$
\end{thm}
\begin{proof}
    The proof is almost identical to \cref{thm:type3a}, except that we are using a different divisor here. We will prove that $I(F_3^{1}( [n]))=1$ For $n=0$, this holds since $\lambda, 12\lambda-\delta_{\text{irr}}$ contracts $F_3^{1}(\emptyset)$. The case of $n=1$ follows from \cref{eg:M31} and \cref{prop:charpsemi}.
    
    We will first verify that $F_3^{1}( [n])$ satisfies the conditions of \cref{lem:second}. Let $N_n$ be the set of nef divisors that intersect $F_3^{1}( [n])$ trivially. Since the image of $F_3^{1}( [n])$ under any projection is $F_3^{1}( [n])$, it suffices to show that the image of $f:N_n\to V_{3,n}$ spans $V_{3,n}$ for each $n$. 
    
    Define the divisor
    \[
    D_n=\kappa+\delta_{2,\emptyset}
    \]
    on $\M{3}{n}$ for $n\geq 2$. This divisor is nef by \cref{thm:kapnef}, see also \cref{eg:kapeg} (1). Furthermore, by \cite[Theorem 2.1]{GKM02}, we have $D_n\cdot F_3^{1}( [n])=0$. For $n=2$, we compute
    \[
    f(D_2)=F_5, \quad f(\psi_1)=F_5+F_6.
    \]
    Thus, $f$ satisfies the assertion. For $n>2$,
    \[
    f(D_2)=-e_0+(-1)^{n+1}e_1, \quad f(\psi_1)=-e_0.
    \]
    Hence, in all cases, $f$ satisfies the assertion. Therefore, by \cref{lem:second}, for any $n\geq 2$, the image of $N_n\to \R^{F^{\K}_{g,n}}$ generates $\R^{F^{\K}_{g,n}}$.
    
    Now, we apply \cref{lem:extind} with $C_n=[F_3^{1}( [n])]$. We have already verified condition (1) for $n=0$. Moreover, condition (2) holds for every $n$ by the preceding paragraph. Hence, by induction using \cref{lem:extind}, we conclude that $I(F_3^{1}( [n]))=1$ for every $n\geq 1$.
\end{proof}

This completes the proof of \cref{thm:main2} (2). To prove \cref{thm:main1}, we need the following theorem:

\begin{thm}\label{thm:f31n}
    For $g\ge 3$ and $n\ge 1$, $I(F_3^{1}([n]))=\lfloor\frac{g}{2}\rfloor$ on $\M{g}{n}$.
\end{thm}

\begin{proof}

    Define
    \[
    C_n:=\left\{ F_3^{1}([n]), F_5^{1,i}([n], \emptyset)\ \middle|\ 1\leq i\leq g-2 \right\}.
    \]
    
    \textbf{Step 1.} If a nef divisor intersects $F_3^{1}([n])$ trivially, then it also intersects any element of $C_n$ trivially. Let
    \[
    D=a\lambda-b_{\text{irr}}\delta_{\text{irr}}-\sum b_{i,I}\delta_{i,I}
    \]
    be such a nef divisor, where we use the convention $\delta_{0, \left\{i\right\}}=-\psi_i$. Then
    \[
    D\cdot [F_3^{1}([n])]=b_{1, [n]}=0.
    \]
    Since $D$ is nef, we obtain the following inequalities:
    \begin{gather*}
        D\cdot F_5^{1,i}([n], \emptyset)=b_{i, \emptyset}+b_{1, [n]}-b_{i+1, [n]}=b_{i, \emptyset}-b_{g-i-1,\emptyset}\geq 0,\\
        D\cdot F_5^{1,g-i-1}([n], \emptyset)=b_{g-i-1, \emptyset}+b_{1, [n]}-b_{g-i, [n]}=b_{g-i-1,\emptyset}-b_{i, \emptyset}\geq 0.
    \end{gather*}
    Thus, $D\cdot F_5^{1,i}([n], \emptyset)=D\cdot F_5^{1,g-i-1}([n], \emptyset)=0$, so $D$ satisfies the assertion.
    
    Therefore, it suffices to prove that $I(C_n)=\lfloor\frac{g}{2}\rfloor$.
    
    \textbf{Step 2.} The dimension of the subspace generated by $C_n$ is $\lfloor\frac{g}{2}\rfloor$.
    
    By the explicit intersection formula shown in Step 1, identifying $\A_1$ as the dual space of $\Pic$, the mappings
    \begin{gather*}
        D\mapsto b_{1, [n]},\\
        D\mapsto b_{i, \emptyset}-b_{g-i-1,\emptyset} \quad \text{for } 1\leq i\leq \lfloor\frac{g}{2}\rfloor-1
    \end{gather*}
    form a basis of the subspace. Hence, the dimension is $\lfloor\frac{g}{2}\rfloor$.

  Therefore, $I(F_3^{1}([n]))=\lfloor\frac{g}{2}\rfloor$ implies the last condition of \cref{lem:extind} (1).

    \textbf{Step 3.} $I(F_3^{1}([n]))=\lfloor\frac{g}{2}\rfloor$ for $n=0,1$.
    
    Let $N_n$ be the set of nef divisors that intersect each element of $C_n$ trivially. To prove the case of $n=0$, by \cref{lem:picnum}, it is enough to show that $N_0$ contains two linearly independent divisors. Since $\lambda, 12\lambda-\delta_{\text{irr}}\in N_0$, this holds.

    Consider the case of $n=1$. Again, by \cref{lem:picnum}, we need to produce $g-\lfloor\frac{g}{2}\rfloor+2$ linearly independent divisors in $N_1$. The divisors
    \[
    \lambda, 12\lambda-\delta_{\text{irr}}, \kappa+\delta_{g-1, \emptyset}, \kappa+\delta_{1, \{1\}}+\frac{1}{4}(\delta_{i, \emptyset}+\delta_{g-1-i, \emptyset}) \quad \text{for } 1\leq i \leq \lfloor\frac{g-1}{2}\rfloor
    \]
    are nef by \cref{thm:kapnef}, \cref{prop:14new} and \cref{prop:last}, and they intersect $F_3^{1}([n])$ trivially by \cref{thm:explicit}. This forms a linearly independent set of
    \[
    3+\lfloor\frac{g-1}{2}\rfloor= g-\lfloor\frac{g}{2}\rfloor+2
    \]
    divisors in $N_1$. This establishes the case of $n=1$.

    \textbf{Step 4.} Using \cref{lem:second}.
    
    Note that $\pi_{i\ast} C_{n+1}=C_{n}$ for any $1\leq i\leq n+1$. Let $f:N_n\to V_{g,n}$ be the composition of $E$ and the intersection pairing. We need to show that the image of $f$ spans $V_{g,n}$. Define
    \[
    D_1:=\kappa+\delta_{1, [n]}, \quad \text{and} \quad D_i:= \kappa+\delta_{1, [n]}+\delta_{i, \left\{1\right\}} \text{ for } 2\leq i\leq \lfloor\frac{g}{2}\rfloor.
    \]
    These are nef by \cref{thm:kapnef}, and they intersect trivially with $F_3^{1}([n])$, hence with any element of $C_n$. If $n=2$,
    \begin{gather*}
        f(\psi_1)=\sum_{j=1}^{\lfloor\frac{g}{2}\rfloor}F_6^{0,0,j, g-j}, \quad f(D_1)=\sum_{j=2}^{\lfloor\frac{g}{2}\rfloor}F_6^{0,0,j, g-j},\\
        f(D_i)=\sum_{j=2}^{\lfloor\frac{g}{2}\rfloor}F_6^{0,0,j, g-j}-F_6^{0,0,i, g-i} \quad \text{for } 2\leq i\leq \lfloor\frac{g}{2}\rfloor.
    \end{gather*}
    Hence, the image of $f$ spans $V_{g,n}$. Here, $F_6^{0,0,j, g-j}$ is an abbreviation for $F_6^{0,0,j, g-j}( \left\{1\right\},\left\{2\right\}, \emptyset, \emptyset)$. 
    
    If $n>2$ and $g$ is odd or $n$ is even, 
    \[
    f(\psi_1)=-e_0, \quad f(D_1)=-e_0+(-1)^{n+1}e_1, \quad f(D_i)=-e_0+(-1)^{n+1}e_1+e_i  \quad \text{for } 2\leq i\leq \lfloor\frac{g}{2}\rfloor.
    \]
    A similar result holds when $g$ is even and $n$ is odd. Hence, in all cases, $f$ satisfies the assertion. By \cref{lem:second}, the image of $N_n\to \R^{F^{\K}_{g,n}}$ spans $\R^{F^{\K}_{g,n}}$.
    
    By Steps 3 and 4, we can apply \cref{lem:extind}, and we obtain $I(F_3^{1}([n]))=\lfloor\frac{g}{2}\rfloor$ for every $n$.

\end{proof}

\subsection{Type 5 and 6}\label{subsec:56}

In this subsection, we prove \cref{thm:main2det} (3) and (4). The situation of type 5 and 6 F-curves is more subtle. Based on an explicit computation using \cite{Choigit25}, we propose the following conjecture:

\begin{conj}\label{conj:exttype5}
    Any type 5 F-curve on $\M{g}{n}$ spans a regular extremal ray, except for the following cases:
    \begin{enumerate}
        \item $F_5^{i,j}(I, \emptyset)$, where $i + 2j < g$ and $j \ne 0$.
        \item $F_5^{i,j}([n], \emptyset)$, where $i + 2j = g$, $j \ne 0$, and $(g,n,i) \ne (2,n,0),\ (3,n,1),\ (4,0,2)$.
    \end{enumerate}
\end{conj}

As reported in \cref{thm:main3}, we have verified \cref{conj:exttype5} for $g \le 4$. This will be proved in \cref{subsec:g23}. Note that the two exceptions in \cref{conj:exttype5} are clearly not extremal, since 
\[
[F_5^{i,j}(I, \emptyset)] = \frac{1}{2}\left([F_6^{j,j,i,g-i-2j}(\emptyset, \emptyset, I, I^c)] + [F_5^{i,2j}(I, \emptyset)] \right)
\]
in the first case, and
\[
[F_5^{i,j}([n], \emptyset)] = [F_3^{i}([n])]
\]
in the second case (cf. \cref{thm:main2}).

Type 6 F-curves appear to be more complicated. For instance, \cref{rmk:type6} presents an example of a non-extremal type 6 F-curve on $\overline{\rm{M}}_7$. However, this is not the minimal genus example. On $\M{6}{1}$, we have
\[
[F_6^{0,1,2,3}(\left\{1\right\}, \emptyset, \emptyset, \emptyset)] = \frac{1}{2} \left([F_6^{0,2,2,2}(\left\{1\right\}, \emptyset, \emptyset, \emptyset)] + [F_6^{0,1,1,4}(\left\{1\right\}, \emptyset, \emptyset, \emptyset)]\right).
\]
At present, the author does not have a clear conjecture regarding which type 6 F-curves are (regular) extremal.

From now, Rather than attempting to classify all (regular) extremal rays, we instead focus on the following collection of type 5 and 6 curves.

\begin{defn}\label{defn:Knudsen}
    An F-curve $C$ on $\M{g}{n}$ is said to be of \textbf{Knudsen type} if there exists a projection map $\pi:\M{g}{n} \to \M{g}{m}$ such that $\pi_\ast [C] \in F^{\K}_{g,m}$. Equivalently, such a curve is of the form $F_5^{0,0}(I,J)$ or $F_6^{0,0,i,g-i}(I,J,K,L)$.
\end{defn}

The motivation behind this definition is as follows. Our method relies on establishing the regular extremality of certain curves for small values of $n$ as a base case. However, in contrast to the situation for larger $n$, the space of semigroup $\kappa$ divisors for small $n$ is significantly more constrained, making it difficult to prove the base case directly. For F-curves of Knudsen type, however, we can invoke \cref{cor:simp} as the base step.

The first main result of this section is the following:

\begin{thm}\label{thm:Knutype}
    Any F-curve of Knudsen type forms a regular extremal ray. 
\end{thm}

The idea of the proof is similar to that of \cref{subsec:3}, but the argument is more involved due to two main reasons: (1) They are not stable under permutation of the marked points, and (2) a curve may become an element of $F^{\K}_{g,n}$ after forgetting a marked point. Hence, extra care is needed when applying \cref{lem:extind} and \cref{lem:second}. Note that $F^{\K}_{g,n}$ implicitly depends on the choice of two marked points $i$ and $j$. We will use the original slightly ambiguous notation for simplicity.

Roughly speaking, the proof breaks into two steps. First, we consider the case of \( F_5^{0,0}(p, \{q,r\}) \) or \( F_6^{0,0,i,g-i}(p, \{q,r\}, I, J) \) as in \cref{Prop:extn}, which corresponds to the situation where there exists a point \( x \) such that \( \pi_{x,\ast}C \in F_{g,n-1}^{\K} \) for some choice of \( x \) and \( F_{g,n-1}^{\K} \). This part is more ad hoc, but manageable, since the setting is restricted. We will use induction on $n$, where \cref{Prop:ext3} serves as the base step. After this, we may then assume that \( \pi_{x,\ast}C \notin F_{g,n-1}^{\K} \) for any \( x \) and \( F_{g,n-1}^{\K} \). This is the more general and simpler case.

We introduce the following notation. For any subset $S \subseteq \NE{\M{g}{n}}$, define $N_S$ to be the $\Q$-linear span of the set of nef divisors that intersect every element of $S$ trivially.

\begin{prop}\label{Prop:extlem2}
    If $C$ is a curve on $\M{g}{n}$ whose curve class is not contained in the linear span of $F^{\K}_{g,n}$ and spans a regular extremal ray, then $C$ satisfies the conclusion of \cref{lem:second}, i.e., 
    \[
    N_C \twoheadrightarrow \Q^{F^{\K}_{g,n}}.
    \]
\end{prop}

\begin{proof}
    Since $C$ is a regular extremal ray, $N_C$ is a codimension $1$ subspace of the $\Q$-Picard group consisting of divisors that intersect trivially with $C$. If the map $N_C \twoheadrightarrow \Q^{F^{\K}_{g,n}}$ is not surjective, then there exists a class $C'$ in $\text{A}_1(\M{g}{n})$, contained in the linear span of $F^{\K}_{g,n}$, such that every divisor in $N_C$ intersects trivially with $C'$. However, this contradicts the codimension of $N_C$, since $C$ and $C'$ are linearly independent.
\end{proof}

\begin{prop}\label{Prop:ext3}
    If $C$ is $F_5^{0,0}(p, \left\{q,r\right\})$ or $F_6^{0,0, i, g-i}(p, \left\{q,r\right\}, \emptyset, \emptyset)$ on $\M{g}{3}$ , then 
    \begin{enumerate}
        \item $N_C \twoheadrightarrow \Q^{F^{\K}_{g,n}}$ for any choice of $F^{\K}_{g,n}$, and
        \item $C$ spans a regular extremal ray.
    \end{enumerate}
\end{prop}

\begin{proof}
    \textbf{Case 1.} $F^{\K}_{g,n}$ corresponds to $\{q, r\}$.
    
    The space \( N_C \) contains \( \pi_p^\ast \text{Pic}\left( \M{g}{2} \right) \), so we apply \cref{lem:second}. By \cref{thm:kapnef}, and in particular \cref{eg:kapeg}~(2), we have
    \[
    D_j = \kappa + \delta_{0,\{q,r\}} + \delta_{j,p} + \delta_{g-j,\emptyset}
    \]
    which are nef for \( 0 < j < \frac{g}{2} \), and contract both \( F_5^{0,0}(p,\{q,r\}) \) and \( F_6^{0,0,i,g-i}(p,\{q,r\},\emptyset,\emptyset) \).
    
    For the map \( f: N_3 \to V_{g,3} \) in \cref{lem:second}, we have
    \[
    f(\psi_r) = -e_0 \quad
    f(D_j) = -2e_0-2e_j
    \]
    
    Hence, we obtain a surjection \( N_C \twoheadrightarrow \Q^{F^{\K}_{g,n}} \).

    \textbf{Case 2.} $F^{\K}_{g,n}$ corresponds to $\{p,q\}$ or $\{p,r\}$, and $C=F_5^{0,0}(p,\{q,r\})$.
    
    Without loss of generality, assume that $F^{\K}_{g,n}$ corresponds to $\{p,q\}$. We have $\pi_r^\ast N_{F_5^{0,0}(p,q)} \subseteq N_C$. Although we have \cref{thm:Knuses}, this is not enough, since (1) they cannot distinguish between $F_6^{0,0,j,g-j}(p,q,r,\emptyset)$ and $F_6^{0,0,j,g-j}(p,q,\emptyset,r)$, and (2) they always vanish on $F_5^{0,0}(p,q)$. Issue (1) is resolved by using exactly the same divisor as in Case 1. Issue (2) is resolved because $\kappa + \delta_{0,\{q,r\}}$ is a nef divisor that contracts $C=F_5^{0,0}(p,\{q,r\})$ but does not contract $C=F_5^{0,0}(p,q)$.
    
    \textbf{Case 3.} $F^{\K}_{g,n}$ corresponds to $\{p,q\}$ or $\{p,r\}$, and $C=F_6^{0,0,i,g-i}(p,\{q,r\},\emptyset,\emptyset)$.
    
    We follow the approach in Case 2. We also have $\pi_r^\ast N_{F_6^{0,0,i,g-i}(p,q,\emptyset,\emptyset)} \subseteq N_C$. This set of divisors (1) cannot distinguish between $F_6^{0,0,j,g-j}(p,q,r,\emptyset)$ and $F_6^{0,0,j,g-j}(p,q,\emptyset,r)$, and (2) always vanishes on $F_6^{0,0,i,g-i}(p,q,r,\emptyset)$ and $F_6^{0,0,i,g-i}(p,q,\emptyset,r)$. This issue is resolved because $\kappa + \delta_{i,r}$, $\kappa + \delta_{g-i,r}$, and the divisors in Step 1 are all contained in $N_C$.
    
    (2) This follows from (1), \cref{cor:simp}, and \cref{lem:extind} for $\{p,q\}$.

\end{proof}

\begin{prop}\label{Prop:extn}
    If $C$ is $F_5^{0,0}(p, \left\{q,r\right\})$ or $F_6^{0,0, i, g-i}(p, \left\{q,r\right\}, I,J)$ on $\M{g}{n}$, then 
    \begin{enumerate}
        \item $C$ spans a regular extremal ray, and
        \item $N_C \twoheadrightarrow \Q^{F^{\K}_{g,n}}$ for any choice of $F^{\K}_{g,n}$ does not contain \([C]\).
    \end{enumerate}
\end{prop}

\begin{proof}
We may exclude the cases \( |I| = 1,\, i = 0 \) and \( |J| = 1,\, i = g \), since these cases follow from \cref{cor:simp} and \cref{Prop:extlem2}.

We will use induction on $n$ and prove (1) and (2) at once. The $n=3$ case is proved in \cref{Prop:ext3}. Consider the case $n\ge 4$. Fix a marked point $s\not \in \left\{p,q,r\right\}$. 

First, we will prove that $N_C \twoheadrightarrow \Q^{F^{\K}_{g,n}}$, where this $F^{\K}_{g,n}$ corresponds to $\left\{r,s\right\}$. By \cref{cor:simp}, \cref{Prop:extlem2}, and the induction hypothesis for (2), (1) of \cref{lem:second} holds. Hence, it is enough to verify (2) of \cref{lem:second}. Let $f:N_C\to V_{g,n}$ be the map in \cref{lem:second} (2). Note that
\[
D_0=\kappa+\delta_{0,\left\{q,r\right\}}, \quad D_j=\kappa+\delta_{0,\left\{q,r\right\}}+\delta_{j,p}+\delta_{g-j,\left\{p,q,r\right\}^c}
\]
are nef line bundles by \cref{thm:kapnef}, and contained in $N_C$ for $1\le j\le \frac{g}{2}$. Then by \cref{lem:comp},
\[
f(D_0)=-2e_0,\quad f(D_j)=-2e_0+2e_j.
\]
Thus, (2) of \cref{lem:second} is satisfied, so $N_C \twoheadrightarrow \Q^{F^{\K}_{g,n}}$.

Now, (1) follows from \cref{lem:extind} with $\left\{r,s\right\}$ and the induction hypothesis for (1). (2) follows from (1) and \cref{Prop:extlem2}.

\end{proof}

\begin{proof}[Proof of \cref{thm:Knutype}]
    We will use induction on $n$ to prove the following statement: If $C$ is of Knudsen type, i.e., of the form $F_5^{0,0}(I,J)$ or $F_6^{0,0,i,g-i}(I,J,K,L)$, then
    \begin{enumerate}
        \item $C$ spans a regular extremal ray, and
        \item $N_C \twoheadrightarrow \Q^{F^{\K}_{g,n}}$ for any choice of $F^{\K}_{g,n}$ such that $[C]\not\in F^{\K}_{g,n}$. 
    \end{enumerate}
   For simplicity of the argument, assume that if $C=F_6^{0,0,i,g-i}(I,J,K,L)$ and $i=0$ (resp. $i=g$), then $|I|, |J| \le |K|$ (resp. $|I|, |J| \le |L|$), and $|I| \le |J|$.

    The case $n\le 3$ is covered by \cref{cor:simp} and \cref{Prop:extn}. Hence, we may assume that $n\ge 4$. Since (2) follows from (1) and \cref{Prop:extlem2}, it suffices to prove (1). If $|I|+|J|\le 3$, then this is again proved in \cref{cor:simp} and \cref{Prop:extn}. Therefore, we may assume that $|I|+|J|\ge 4$. 
    
    Note that this condition implies $\pi_{x,\ast}C\not\in F^{\K}_{g,n-1}$ for any choice of $F^{\K}_{g,n-1}$ and $x$. Hence, by the induction hypothesis for (2), condition (1) of \cref{lem:second} is satisfied. Choose $p\in I$ and $q,r\in J$ (this is possible since $|I|+|J|\ge 3$ and $|I|\le |J|$), and let $F^{\K}_{g,n}$ be the set of curves corresponding to $\left\{p,q\right\}$. Define
    \[
    D_0=\psi_r, \quad D_j=\kappa+\delta_{0,J}+\delta_{j,\{p,r\}}
    \]
    for $1\le j\le \frac{g}{2}$, which are nef line bundles contracting $C$ by \cref{thm:kapnef}. Then by \cref{lem:comp},
    \[
    f(D_0)=-e_0, \quad f(D_j)=((-1)^{|J|+1}-1)e_0-e_j,
    \]
    Thus, condition (2) of \cref{lem:second} is satisfied. Therefore, $N_C \twoheadrightarrow \Q^{F^{\K}_{g,n}}$
    
    Note that condition (1) of \cref{lem:extind} for \( \{p, q\} \) also holds by the induction hypothesis. If \( I = \{p\} \), then \( I(\pi_{p,\ast}C) = I(\pi_{p,\ast}\pi_{q,\ast}C) = 0 \). Otherwise, by the induction hypothesis, we have
\[
I(\pi_{p,\ast}C) = I(\pi_{q,\ast}C) = I(\pi_{p,\ast}\pi_{q,\ast}C) = 1.
\]
Hence, \( I(C) = 1 \), so \( C \) spans a regular extremal ray. This also implies condition (2).

\end{proof}

Next, we prove the regular extremality of certain type~6 F-curves presented in \cref{thm:main2}~(4). The method used here is entirely different from that of \cref{subsec:method}. Instead of that approach, we employ \cref{thm:01push} together with Hassett's moduli spaces of weighted pointed stable curves~\cite{Has03}. By \cref{thm:01push}, the pushforward of any regular extremal ray under the morphism \(F : \M{0}{g+n} / S_g \to \M{g}{n} \) remains a regular extremal ray. By analyzing the F-curves contracted by a specific Hassett space, we deduce the regular extremality of certain F-curves on \( \M{0}{g+n} / S_g \). To address the quotient by \( S_g \), we require the following lemma.

\begin{lem}\label{lem:quot}
    Let $V,W$ be finite-dimensional $\Q$-vector spaces with an action of a finite group $G$, and $f:V\to W$ be a $G$-equivariant map. Let $V^G$ (resp. $V_G$) be the set of $G$-invariants (resp. $G$-coinvariants), and note that $(V^\ast)^G$ can be naturally identified with $V_G^\ast$. Under this identification, the orthogonal complement of the image of $f^\ast:(W^\ast)^G\to (V^\ast)^G$ is $(\ker f)_G\subseteq V_G$. 
\end{lem}

\begin{proof}
    By restriction, we have a natural map $(V^\ast)^G\to V_G^\ast$. It is straightforward to see that this is an isomorphism by breaking $V$ into simple $G$-representations, since the base field is of characteristic $0$.

    Now we will consider the second assertion. We have an exact sequence
    \[ 0\to \ker f\to V\to W. \]
    By taking dual and $G$-invariants, we have
    \[ (W^\ast)^G\to (V^\ast)^G\to (\ker f^\ast)^G\to 0. \]
    By the identification above, the last map is equal to $(V_G)^\ast \to (\ker f)_G^\ast$ induced from the inclusion $(\ker f)_G\to V_G$. Hence the statement follows from taking duals.
\end{proof}

\begin{thm}\label{thm:type6}
    The following type 6 F-curves are regular extremal:
        \begin{enumerate}
            \item $F_6^{1,1,1,g-3}(\emptyset, \emptyset, \emptyset, [n])$,
            \item $F_6^{1,1,2,g-4}(\emptyset, \emptyset, \emptyset, [n]) \text{ for }n\ne 0$,
            \item $ F_6^{0,1,1,g-2}(i, \emptyset, \emptyset, [n]\setminus i)$ for all $i\in [n]$.
        \end{enumerate}
\end{thm}

\begin{proof}
    We will denote the index set for points on $\M{0}{g+n}$ by $[g]\cup [n]$, where $[g]$ is the set of symmetric points (i.e. where $S_g$ is acting on) and $[n]$ is the set of asymmetric points. Moreover, the $i$th symmetric (resp. asymmetric) point will be denoted by $i_s$ (resp. $i_a$).

    Define a sequence $\mathcal{A}_1:=(a_j)_{j\in [g]\cup [n]}$ by
    \[ a_j=\begin{cases}
        \frac{1}{3} & \text{if } j\in [g],\\
         1 & \text{if } j\in [n].
    \end{cases} \]
    Consider the natural contraction $f:\M{0}{g+n}\to \M{0}{\mathcal{A}_1}$. By \cite[Proposition 4.5]{Has03}, the exceptional locus of $f$ is the union of $\M{0}{I+1}\times \M{0}{I^c+1}$ where $I\subseteq [g]$ and $|I|=3$. Moreover, $f$ contracts this boundary divisor into $\M{0}{I^c+1}$. Therefore, the set of F-curves contracted by $f$ is exactly $F_6^{0,0,0,0}(p_s,q_s,r_s,I^c)$ for $I=\{p_s,q_s,r_s\}$. Since $f$ is a smooth blowdown corresponding to images of F-curves (see also \cite[Lemma 4.6, Proposition 4.7]{Fak12}), if we let $U$ be the subspace of $A_1(\M{0}{g+n})$ generated by these $F$-curves, then we have a natural exact sequence
    \[ 0\to \Pic(\M{0}{\mathcal{A}_1})\to \Pic(\M{0}{g+n})\to U^\ast  \]
    and hence, by taking duals, we have 
    \[ \ker f_{\ast}=\text{Im }U:=\text{image of }U\text{ in }A_1(\M{0}{g+n}) \]
    where $f_\ast:A_1(\M{0}{g+n})\to A_1(\M{0}{\mathcal{A}_1})$. Since $\mathcal{A}_1$ is $S_g$-invariant, $f_\ast$ is also $S_g$-equivariant, so we can apply \cref{lem:quot}, which gives
    \[ (f^\ast \Pic(\M{0}{\mathcal{A}_1})^{S_g})^\perp = \text{Im }U/S_g\subseteq A_1(\M{0}{g+n})/S_g=A_1(\M{0}{g+n}/S_g).  \]
    Note that $f^\ast \Pic(\M{0}{\mathcal{A}_1})^{S_g}$ is generated by nef divisors, and $\text{Im }U/S_g$ is $1$-dimensional since the set of $F_6^{0,0,0,0}(p,q,r,I^c)$ is transitive under the $S_g$-action. Therefore, the image of $F_6^{0,0,0,0}(p_s,q_s,r_s,I^c)$ in $A_1(\M{0}{g+n}/S_g)$ is a regular extremal ray of $\NE{\M{0}{g+n}/S_g}$. By \cref{thm:01push}, applying the flag map yields that $F_6^{1,1,1,g-3}(\emptyset, \emptyset, \emptyset, [n])$ is a regular extremal ray. 

    For $F_6^{0,1,1,g-2}(i, \emptyset, \emptyset, [n]\setminus i)$, we use essentially the same argument, with $\mathcal{A}_2:=(a_j)_{j\in [g]\cup [n]}$ defined by
    \[ a_j=\begin{cases}
        \frac{1}{2}-\epsilon & \text{if } j\in [g],\\
         \epsilon & \text{if } j=i_a,\\
         1 & \text{otherwise}.
    \end{cases} \]
    for a sufficiently small positive number $\epsilon$. Then the F-curves contracted by $f:\M{0}{n}\to \M{0}{\mathcal{A}_2}$ are $F_6^{0,0,0,0}(i_a,p_s,q_s,\{p_s,q_s,i_a\}^c)$ where $p_s,q_s\in [g]$. The proof is identical, so we omit the details.

    For the last case, $F_6^{1,1,2,g-4}(\emptyset, \emptyset, \emptyset, [n])$, we use $\mathcal{A}_3:=(a_j)_{j\in [g]\cup [n]}$ defined by
    \[
    a_j=\begin{cases}
        \frac{1}{4} & \text{if } j\in [g],\\
         1 & \text{if } j\in [n].
    \end{cases}
    \]
    By \cite[Proposition 4.5]{Has03}, the exceptional locus of $f:\M{0}{n}\to \M{0}{\mathcal{A}_3}$ is the union of $\M{0}{I+1}\times \M{0}{I^c+1}$ where $I\subseteq [g]$ and $|I|=3$ or $4$, and $f$ contracts this divisor to $\M{0}{I^c+1}$. Therefore, the set of $F$-curves contracted by $f$ is $F_6^{0,0,0,0}(p_s,q_s,r_s,I^c)$ and $F_6^{0,0,0,0}(p_s,q_s,\{r_s, r_s'\},I^c)$, where $I=\{p_s,q_s,r_s\}$ and $I=\{p_s,q_s,r_s, r_s'\}$, respectively. Note that the set of such $F$-curves has two $S_g$-orbits, corresponding to $|I|=3$ and $|I|=4$. We already know that the $|I|=3$ case is regular extremal in $\NE{\M{0}{g+n}/S_g}$, and we will prove that the $|I|=4$ case is also regular extremal. By applying the flag map and \cref{thm:01push}, this finishes the proof.
        
    By the same argument as in the first case, we have
    \[
    (f^\ast \Pic(\M{0}{\mathcal{A}_3})^{S_g})^\perp = \text{Im }U/S_g\subseteq A_1(\M{0}{g+n})/S_g=A_1(\M{0}{g+n}/S_g).
    \]
    However, in this case, $U/S_g$ is $2$-dimensional, generated by the images of $F_6^{0,0,0,0}(p_s,q_s,r_s,I^c)$ and $F_6^{0,0,0,0}(p_s,q_s,\{r_s, r_s'\},I^c)$. Hence, it is enough to construct an $S_g$-invariant nef line bundle which contracts $F_6^{0,0,0,0}(p_s,q_s,\{r_s, r_s'\},I^c)$ and does not contract $F_6^{0,0,0,0}(p_s,q_s,r_s,I^c)$. 
    
    For this, we use nef divisors on $\M{0}{g+n}$ coming from GIT quotients \cite{AS11, GG12}. We follow the notation of \cite[Section 2]{GG12}. Define a sequence $(x_i)_{i\in [g]\cup [n]}$ with $x_i=\frac{1}{2}$ for $i\in [g]$ and $x_i=0$ or $\frac{1}{2}$ for $i\in [n]$ so that the sum of all $x_i$'s is an integer $d+1$ for $d\ge 1$. Then, by \cite[Theorem 2.1]{GG12}, we have a nef divisor $D:=\varphi_{d, \vec{x}}^\ast \mathcal{O}(1)$ satisfying
    \[
    D\cdot F_6^{0,0,0,0}(p_s,q_s,r_s,I^c)\ne 0,\quad D\cdot F_6^{0,0,0,0}(p_s,q_s,\{r_s, r_s'\},I^c)=0.
    \]
    This proves the $|I|=4$ case, and hence the theorem.
\end{proof}

Note that, by \cref{rmk:type6}, the conclusion of \cref{thm:type6} (2) does not hold when $n=0$.

\subsection{Small genus}\label{subsec:g23}

In this subsection, we prove \cref{thm:main3}, which amounts to classifying regular extremal F-curves of types~5 and~6. Note that dealing with type~5 F-curves amounts to proving \cref{conj:exttype5} in this case. The non-extremal type~5 F-curves in \cref{thm:main3} are exactly the exceptions listed in \cref{conj:exttype5}, so it suffices to show the extremality of the remaining type~5 F-curves.

We begin with the case of $\M{2}{n}$. By \cref{thm:Knutype}, it suffices to show that the curves $F_5^{0,1}(I, J)$ are regular extremal. 

\begin{thm}\label{thm:genus2}
    $F_5^{0,1}(I, J)$ are regular extremal on $\M{2}{n}$.
\end{thm}

\begin{proof} 
We proceed by induction on $n$. For $n \le 5$, the F-conjecture is known to hold, so the F-cone coincides with the nef cone and is polyhedral. In this range, we can verify the extremality of these F-curves directly by brute-force computation. This verification is carried out using a Python script available at \cite{Choigit25}.

Now assume $n \ge 6$ and let $F = F_5^{0,1}(I, J)$ (in fact, $n \ge 4$ suffices). Choose three indices $p, q, r \in [n]$ such that $I \ne \{r\}$ and $I \not\subseteq \{p, q\}$. Then $\pi_{p, \ast}F$, $\pi_{q, \ast}F$, and $\pi_{q, \ast} \pi_{p, \ast}F$ are F-curves of the same type. By the induction hypothesis, \cref{lem:extind} (1) holds with index of extremality 1. Therefore, to apply \cref{lem:extind}, it remains to verify that conditions in \cref{lem:second} holds. Since each $\pi_{i, \ast}F$ is either an F-curve of the same type or zero, \cref{lem:second} (1) follows from the induction hypothesis and \cref{Prop:extlem2}. Thus, it suffices to check \cref{lem:second} (2).

\textbf{Case 1. } $J \ne \emptyset$

By \cref{thm:kapnef}, the divisor $D = \kappa + \delta_{1, J}$ is nef and contracts $F$. Moreover, for the map $f : N_n \to V_{2,n}$, by \cref{lem:comp}, we have
\[
f(\psi_r) = -e_0, \quad f(D) = -e_0+(-1)^{|J|+1}e_1.
\]
Hence, the condition in \cref{lem:second} (2) is satisfied.

\textbf{Case 2. } $J = \emptyset$

By \cref{thm:kapnef}, the divisor $D = \kappa + \delta_{1, \emptyset} + \delta_{2, \emptyset}$ is nef and contracts $F$. Moreover, for the map $f : N_n \to V_{2,n}$, by \cref{lem:comp}, we have
\[
f(\psi_r) = -e_0, \quad f(D) = - 2e_0 - e_1.
\]
Therefore, the condition in \cref{lem:second} (2) is satisfied.
\end{proof}

Now consider the case of $\M{3}{n}$. Apart from the cases covered by \cref{thm:Knutype}, we need to consider the following four cases:
\begin{itemize}
    \item $F_5^{0,1}(I, J)$ with $J \ne \emptyset$,
    \item $F_5^{1,1}(I, J)$,
    \item $F_5^{0,2}(I, J)$,
    \item $F_6^{0,1,1,1}(I, J, K, L)$.
\end{itemize}

\begin{thm}\label{thm:genus3}
    All of the above four types of F-curves are regular extremal on $\M{3}{n}$.
\end{thm}

\begin{proof}
The proof follows the strategy of \cref{thm:genus2}. As before, we proceed by induction and verify the statement directly for $n \le 4$, using \cite{Choigit25}. Now assume $n \ge 4$. For each F-curve $F$ of above types, choose different indicies $p,q,r$ so that $\pi_{p,\ast}\pi_{q, \ast}F$ is nonzero and $\psi_r\cdot F=0$. This is always possible since $n\ge 4$. Again, following the same reasoning as in \cref{thm:genus2}, it suffices to check that the condition in \cref{lem:second} (2) holds with respect to $p,q$.

\textbf{Case 1. } $F = F_5^{0,1}(I, J)$ with $J \ne \emptyset$

By \cref{thm:kapnef}, the divisor $D = \kappa + \delta_{1, J}$ is nef and contracts $F$. Moreover, for the map $f : N_n \to V_{3,n}$, we have
\[
f(\psi_r) = -e_0, \quad f(D) = -e_0 + (-1)^{|J| + 1}e_1.
\]
Hence, the condition in \cref{lem:second} (2) is satisfied.

\textbf{Case 2. } $F = F_5^{1,1}(I, J)$

Assume that not both $I$ and $J$ are empty. Without loss of generality, we may assume $J \ne \emptyset$. Then we can take $D = \kappa + \delta_{1, J}$, and the proof proceeds as in Case 1.

If $I = J = \emptyset$, let
\[
D = \kappa + \frac{2}{3}\delta_{1, \emptyset} + \frac{1}{3}\delta_{2, \emptyset}.
\]
Then $D$ is nef by \cref{prop:3nnew} and contracts $F$. Moreover, for the map $f : N_n \to V_{3,n}$, we have
\[
f(\psi_1) = -e_0, \quad f(D) = -e_0+\left(-\frac{1}{3} + (-1)^{n-1}\frac{2}{3}\right)e_1.
\]
Hence, the condition in \cref{lem:second} (2) is satisfied.

\textbf{Case 3. } $F = F_5^{0,2}(I, J)$

By \cref{thm:kapnef}, the divisor $D = \kappa + \delta_{2, J}$ is nef and contracts $F$. Moreover, for the map $f : N_n \to V_{3,n}$, we have
\[
f(\psi_r) = -e_0, \quad f(D) = -e_0 + (-1)^{|J^c| + 1}e_1.
\]
Hence, the condition in \cref{lem:second} (2) is satisfied.

\textbf{Case 4. } $F = F_6^{0,1,1,1}(I, J, K, L)$

This is a special case of \cref{thm:type6} (1).

\end{proof}

Now we consider the case of $\M{4}{n}$. Except for the exceptions in \cref{thm:main3} (3) and the cases covered by \cref{thm:Knutype}, we need to prove the regular extremality of the following F-curves:

\begin{itemize}
    \item $F_6^{0,1,1,2}(I,J,K,L)$, $F_6^{1,1,1,1}(I,J,K,L)$, and $F_5^{0,3}(I,J)$.
    \item $F_5^{0,1}(I,J)$ for $J\ne \emptyset$.
    \item $F_5^{1,1}(I,J)$ for $I, J\ne \emptyset$.
    \item $F_5^{0,2}(I,J)$ and $F_5^{2,1}(I,J)$ for $(I,J)\ne ([n], \emptyset)$.  
\end{itemize}

\begin{thm}\label{thm:genus4}
    All of the above F-curves are regular extremal on $\M{4}{n}$.
\end{thm}

\begin{proof}
    The overall strategy is the same. As before, if $n\le 3$, the claim follows from the known cases of the F-conjecture and \cite{Choigit25}, so we assume $n\ge 4$. The following computation will be used multiple times (cf. \cref{lem:comp}):
    \begin{gather*}
        f(\psi_r)=-e_0,\
        f(\kappa+\delta_{1, J})=-e_0+(-1)^{|J|+1}e_1,\
        f(\kappa+\delta_{2, J})=-e_0+(-1)^{|J|+1}e_2,\\
        f\!\left(\kappa+\sum_{i=1}^{4}\delta_{i, \emptyset}\right)=\big((-1)^{n+1}-1\big)e_0+\big((-1)^{n+1}-1\big)e_1-e_2,\\
        f(\kappa+\delta_{2, \emptyset}+\delta_{4, \emptyset})=\big((-1)^{n+1}-1\big)e_0-e_2.
    \end{gather*}
    Note that all of these divisors are nef by \cref{thm:kapnef} and \cref{eg:kapeg}. Moreover, in the proof, we need to verify that the listed divisors contract certain F-curves; however, we omit this step since the computation is evident from \cref{thm:explicit}.

    \textbf{Case 1.} $F_6^{0,1,1,2}(I,J,K,L)$.
    
    Choose $p,q,r$ so that $I \not \subseteq \{p,q\}$ and $I\ne \{r\}$, which is possible since $n\ge 4$. Then $\pi_{p, \ast}F$, $\pi_{q, \ast}F$, and $\pi_{q, \ast}\pi_{p, \ast}F$ are all F-curves of the same type. Hence, to apply \cref{lem:extind}, it remains to check that the condition in \cref{lem:second} holds, which reduces to verifying \cref{lem:second} (2) with respect to $p$ and $q$.
    
    Note that 
    \begin{align*}
        f\!\left(\kappa+\delta_{0, I}+\delta_{1, L\cup i}\right)&=\big((-1)^{|I|+1}-1\big)e_0+(-1)^{|L|}e_1,\\
        f\!\left(\kappa+\delta_{2, J}+\delta_{1, I\cup j}\right)&=-e_0+(-1)^{|I|}e_1+(-1)^{|J|+1}e_2.
    \end{align*}
    These divisors are nef under certain condition, by \cref{eg:kapeg} (4).
    
    To conclude, we select in each situation three of the above divisors whose $f$-images generate $V_{4,n}$ (and, if $n$ is odd, we may ignore $e_2$). When $J=K=L=\emptyset$, this is covered by \cref{thm:type6} (2). If $L=\emptyset$ then at least one of $J$ or $K$ is nonempty (say $J$), so we can take 
    \[
    \psi_r,\ \kappa+\delta_{1, J},\ \kappa+\sum_{i=1}^{4}\delta_{i, \emptyset}.
    \]
    If $L\ne\emptyset$ and at least one of $J$ or $K$ is nonempty (again say $J$), we may choose 
    \[
    \psi_r,\ \kappa+\delta_{1, J},\ \kappa+\delta_{2, L}.
    \]
    When $J=K=\emptyset$ and $L\ne\emptyset$, either $|I|>1$ or $|L|>1$, since $n\ge 4$. In the first case, for any $i\in I$ we can use 
    \[
    \psi_r,\ \kappa+\delta_{0, I}+\delta_{1, L\cup i},\ \kappa+\delta_{2, L}.
    \]
    In the second, for any $l\in L$ we can take 
    \[
    \psi_r,\ \kappa+\delta_{2, L}+\delta_{1, I\cup l},\ \kappa+\delta_{2, L}.
    \]
    In all subcases, these choices generate $V_{4,n}$, and thus $F_6^{0,1,1,2}(I,J,K,L)$ is extremal.

    \textbf{Case 2.} $F_6^{1,1,1,1}(I,J,K,L)$
    
    We proceed as in Case~1: choose any $p,q,r$ and verify \cref{lem:second} (2) for $p$ and $q$.
    
    If at least three of $I,J,K,L$ are empty, then the situation falls under \cref{thm:type6} (1). When exactly two sets are nonempty, say $I$ and $J$, a suitable choice is 
    \[
    \psi_r,\ \kappa+\delta_{1, I},\ \kappa+\sum_{i=1}^{4}\delta_{i, \emptyset}.
    \]
    If at least three sets are nonempty, say $I$, $J$, and $K$, we may take 
    \[
    \psi_r,\ \kappa+\delta_{1, I},\ \kappa+\delta_{1, I}+\delta_{1, J}+\delta_{2, I\cup J}.
    \]
    In each case, the image of these divisors generate the required space, completing the argument.

    \textbf{Case 3.} $F_5^{0,3}(I,J)$.  
    
    Choose $p, q, r$ so that $I \not\subseteq \{p, q\}$ and $I \ne \{r\}$. We will check \cref{lem:second} (2) for $p, q$. Note that  
    \[
    f(\psi_r) = -e_0, \quad f(\kappa + \delta_{3, J}) = -e_0 + (-1)^{|J^c|+1} e_1
    \]
    and that $\psi_r$ and $\kappa + \delta_{3, J}$ are nef divisors. Since $I \ne \emptyset$, we can apply \cref{eg:kapeg} (1). Hence, it suffices to find a nef divisor $D_2$ that contracts $F_5^{0,3}(I,J)$ and whose $e_2$ coefficient in $f(D_2)$ is nonzero when $n$ is even.  
    
    If $J = \emptyset$, we can use $D_2 = \kappa + \sum_{i=1}^{4} \delta_{i, \emptyset}$, for which  
    \[
    f(D_2) = -2e_0 - 2e_1 - e_2
    \]
    when $n$ is even.  
    
    If $|J| \ge 2$, choose $j \in J$ and use  
    \[
    D_2 = \kappa + \delta_{1, j} + \delta_{2, J \setminus j} + \delta_{3, J},
    \]
    where  
    \[
    f(D_2) = -e_0 + \big(1 + (-1)^{|J^c|+1}\big) e_1 + (-1)^{|J|} e_2.
    \]
    
    Hence, from now on, assume $J = \{j\}$. If $|I| \ge 2$, choose $i \in I$ and use $D_2 = \kappa + \delta_{0, I} + \delta_{2, \{i, j\}}$, where  
    \[
    f(D_2) = \left((-1)^{|I|+1}-1\right) e_0 - e_2.
    \]
    Thus, we may assume $I = \{i\}$. Finally, let $D_2 = \kappa + \delta_{1, J^c} + \delta_{2, I^c}$, where  
    \[
    f(D_2) = -e_0 + (-1)^{|J^c|+1} e_1 + (-1)^{|I|+1} e_2.
    \]
    Note that, since $n\ge 4$, this is a semigroup Kappa divisor by \cref{eg:kapeg} (4).
    
    \textbf{Case 4. }$F_5^{0,1}(I,J)$ for $J\ne \emptyset$.

    Choose $p, q, r$ so that $I, J \not\subseteq \{p, q\}$ and $I \ne \{r\}$; this is possible since $n \ge 4$. Using $\psi_r$ and $\kappa+\delta_{1, J}$, we can generate the space spanned by $e_0$ and $e_1$, so it remains to find a nef divisor $D$ contracting $F_5^{0,1}(I,J)$ whose $e_2$-coefficient in $f(D)$ is nonzero. If $|J|\ge 2$, choose $j\in J$; then 
    \[
    \kappa + \delta_{1, J}+\delta_{2, j\cup I}
    \]
    works by \cref{eg:kapeg} (4). If $|I|\ge 2$, choose $i\in I$ and 
    \[
    \kappa+\delta_{0, I}+\delta_{2, i\cup J}
    \]
    suffices. If $I=\{i\}$ and $J=\{j\}$, then 
    \[
    \kappa+\delta_{1, i}+\delta_{1, j}+\delta_{2, i\cup j}
    \]
    is the desired divisor. Therefore, $F_5^{0,1}(I,J)$ is regular extremal.

    \textbf{Case 5. }$F_5^{1,1}(I,J)$ for $I, J\ne \emptyset$.

     Choose $p, q, r$ so that $I, J \not\subseteq \{p, q\}$. This is possible since $n\ge 4$. Using $\psi_r$ and $\kappa+\delta_{1, I}$, we can generate the space spanned by $e_0$ and $e_1$, so it remains to find a nef divisor $D$ contracting $F_5^{1,1}(I,J)$ whose $e_2$-coefficient in $f(D)$ is nonzero. If $I\cup J\ne [n]$, then 
    \[
    \kappa+\delta_{1, I}+\delta_{1, J}+\delta_{2, I\cup J}
    \]
    works. If $I\cup J=[n]$, then $|I|>1$ or $|J|>1$. Without loss of generality, assume $|J|>1$ and choose $j\in J$. Then
    \[
    \kappa+\delta_{1, J}+\delta_{2, I\cup j}
    \]
    suffices by \cref{eg:kapeg} (4). Therefore, $F_5^{1,1}(I,J)$ is regular extremal.

    \textbf{Case 6. }$F_5^{0,2}(I,J)$ for $(I,J)\ne ([n], \emptyset)$.  

    Choose $p,q,r$ so that $I\ne \{r\}$ and the image of $F_5^{0,2}(I,J)$ under $\pi_p, \pi_q, \pi_{\{p,q\}}$ is an F-curve of the same type. If $|I|\ge 3$, choose $p,q$, and $r$ from $I$. If $|I|=2$, choose $p,r\in I$ and $q$ from the complement. If $|I|=1$, choose all of the elements from $I^c$. This is possible since $n\ge 4$.  

    We divide this case into subcases, and for each subcase, we choose three divisors whose images under $f$ generate $V_{4,n}$, spanned by $e_0, e_1$, and $e_2$. We can always take $\psi_r$ for $e_0$. If $J\ne \emptyset$, we may take $\kappa+\delta_{2, J}$ for $e_2$, and we then need one more divisor for $e_1$. If $|J|\ge 2$, choose $j\in J$ and take
    \[
        \kappa+\delta_{2,J}+\delta_{1, j\cup I}.
    \]
    This is nef by \cref{eg:kapeg}~(4), contracts $F_5^{0,2}(I,J)$, and
    \[
        f\left(\kappa+\delta_{2,J}+\delta_{1, j\cup I}\right)=-e_0+(-1)^{|I|}e_1+(-1)^{|J|+1}e_2,
    \]
    so it works. If $|I|\ge 2$, choose $i\in I$ and let 
    \[
        \kappa+\delta_{0,I}+\delta_{1, i\cup J}.
    \]
    This is sufficient for the same reasoning. If $I=\{i\}$ and $J=\{j\}$, then, since $n\ge 3$, $(\delta_{2,j}, \delta_{3,i})$ forms a semigroup. Hence,
    \[
        \kappa+\delta_{2, j}+\delta_{3,i}
    \]
    is nef, contracts $F_5^{0,2}(I,J)$, and 
    \[
        f\left( \kappa+\delta_{2, j}+\delta_{3,i} \right)=-e_0+(-1)^n e_1+e_2.
    \]
    Therefore, the case $J\ne \emptyset$ is settled.
    
    Now consider the case $J=\emptyset$. In this case, $\kappa+\delta_{2, \emptyset}+\delta_{4,\emptyset}$ accounts for the $e_2$ term, so together with $\psi_r$, we still need one more divisor for the $e_1$ term. If $|I|\ge 2$, choose $i\in I$ and $k\in I^c$. Since $(I,J)\ne ([n], \emptyset)$, such $k$ exists. Then
    \[
        \kappa+\delta_{0, I}+\delta_{1, \{i,k\}}
    \]
    works by a similar argument. If $I=\{i\}$, then by \cref{prop:4n2nef},
    \[
        D:=\kappa-\frac{1}{2}\psi_i+\frac{1}{2}\delta_{2, \emptyset}+\frac{1}{2}\delta_{4, \emptyset}+\frac{1}{2}\delta_{1, \{i\}}
    \]
    is nef, contracts $F_5^{0,2}(I,J)$, and 
    \[
        f(D)=\left(\frac{1}{2}(-1)^{n+1}-\frac{1}{2}  \right)e_0 +\frac{1}{2}e_1-\frac{1}{2}e_2.
    \]
    This completes the argument.

    \textbf{Case 7. }$F_5^{2,1}(I,J)$ for $(I,J)\ne ([n], \emptyset)$. 

    Choose $p,q,r$ such that the image under $\pi_p, \pi_q, \pi_{\{p,q\}}$ is an F-curve of the same type. The process is similar to Case~6, so we omit it.

    Note that 
    \begin{align*}
        f\!\left(\kappa+\frac{1}{2}\delta_{1, \emptyset}+\frac{1}{2}\delta_{2, \emptyset} \right)&=-e_0-\frac{1}{2}e_1-\frac{1}{2}e_2.
    \end{align*}
    This divisor is nef by \cref{prop:4nnef}.
    
    To conclude, we select in each situation three of the above divisors whose $f$-images generate $V_{4,n}$. If $I, J \ne \emptyset$, we can take 
    \[
    \psi_r,\ \kappa+\delta_{1, J},\ \kappa+\delta_{2, I}.
    \]
    If $J\ne \emptyset$ and $I=\emptyset$, then 
    \[
    \psi_r,\ \kappa+\delta_{1, J},\ \kappa+\delta_{2, \emptyset}+\delta_{4, \emptyset}
    \]
    works. If $I=J=\emptyset$, then 
    \[
    \psi_r,\ \kappa+\tfrac{1}{2}\delta_{1, \emptyset}+\tfrac{1}{2}\delta_{2, \emptyset},\ \kappa+\delta_{2, \emptyset}+\delta_{4, \emptyset}
    \]
    suffices. Finally, if $J=\emptyset$ and $I\ne \emptyset$, since we already have $\psi_r$ and $\kappa+\delta_{2, I}$, it remains to find one more nef divisor $D$ contracting our F-curve whose $e_1$-coefficient in $f(D)$ is nonzero. In this case, $I\ne \emptyset$ and $I^c\ne \emptyset$, so choose $i\in I$ and $j\in I^c$. Then $(\delta_{2, I}, \delta_{1, \{i,j\}})$ is a semigroup. Hence 
    \[
        D:=\kappa+\delta_{2,I}+ \delta_{1, \{i,j\}}
    \]
    is nef, contracts $F_5^{2,1}(I,J)$, and 
    \[
        f(D)=-e_0-e_1+(-1)^{|I|+1}e_2,
    \]
    so this divisor works. All in all, $F_5^{2,1}(I,J)$ with $(I,J)\ne ([n], \emptyset)$ spans an extremal ray.

\end{proof}

%% file: sections/mainproof.tex
\begin{proof}[Proof of \cref{thm:main1}]
    We begin with the case $g=2$. The case $n=2$ is already covered in \cref{sec:examples}, so assume $n\ge 3$. Define
    \[
        F=\left\{ D\in \mathrm{Nef}(\M{2}{n})\ \middle|\ D\cdot F_3^{0}([n])=0 \right\},
    \]
    and let $E$ be the linear subspace of $F$ in which the coefficient of $-\delta_{0, [n-1]}$ is equal to the coefficient of $\psi_n$. By \cref{thm:type3a}, $E$ is a face of $\mathrm{Nef}(\M{2}{n})$ of codimension~$1$. Since $\psi_n\in F$, $E$ is a proper subspace of $F$. Let $f:\M{2}{2}\to \M{2}{n}$ be the map attaching a rational stable curve with $n$ marked points. Then, for any $D\in F\setminus E$, the divisor $f^\ast D$ has the property that the coefficient of $\psi_1$ is not equal to the coefficient of $\psi_2$. Moreover, since $D\cdot F_3^{0}([n])=0$, we have $f^\ast D\cdot F_3^{0}([2])=b_{0, \{1,2\}}=0$, so $f^\ast D$ satisfies the condition of \cref{thm:semin}. Hence, $f^\ast D$ is not semiample, and therefore $D$ is also not semiample.

   Now consider the case $g\ge 3$. The proof is almost the same as above. Define 
    \[
        F=\left\{ D\in \mathrm{Nef}(\M{g}{n})\ \middle|\ D\cdot F_3^{1}([n])=0 \right\},
    \]
    which is a face of codimension $\lfloor\frac{g}{2}\rfloor$ by \cref{thm:f31n}, and let $E$ be the linear subspace in which the coefficient of $-\delta_{1, [n-1]}$ is equal to the coefficient of $\psi_n$. Since $\psi_n\in F$, $E$ is a proper subspace of $F$. Let $f:\M{g-1}{2}\to \M{g}{n}$ be the map attaching a genus~$1$ stable curve with $n$ marked points. Then, by a similar argument, we can verify that if $D\in F\setminus E$, we may apply \cref{thm:semin} to $f^\ast D$ and deduce that $D$ is not semiample.
\end{proof}

\begin{proof}[Proof of \cref{cor:main1cor}]
    As above, let 
    \[
        F=\left\{ D\in \mathrm{Nef}(\M{g}{n})\ \middle|\ D\cdot F_3^{1}([n])=0 \right\},
    \]
    and let $E$ be the linear subspace in which the coefficient of $-\delta_{1, [n-1]}$ is equal to the coefficient of $\psi_n$. Moreover, define 
    \[
        C_n:=\left\{ F_3^{1}([n]),\ F_5^{1,i}([n], \emptyset)\ \middle|\ 1\leq i\leq g-2 \right\},
    \]
    as in the proof of \cref{thm:f31n}. Let $C$ be the intersection of $\NE{\M{g}{n}}$ with the subspace generated by $C_n$. By the proof of \cref{thm:f31n}, since $C_n$ generates a $\lfloor\frac{g}{2}\rfloor$-dimensional subspace and $I(C_n)=\lfloor\frac{g}{2}\rfloor$, $C$ is a $\lfloor\frac{g}{2}\rfloor$-dimensional face of $\NE{\M{g}{n}}$.
        
    Assume there exists a projective contraction $f:\M{g}{n}\to X$ whose relative closed cone of curves is $C$. Let $D=f^\ast \mathcal{O}_X(1)$. Then $D\in F$ and $D$ is semiample, so by the proof of \cref{thm:main1}, $D\in E$. Write
    \[
        D=a\lambda-b_{\mathrm{irr}}\delta_{\mathrm{irr}}+\sum_{i=1}^n b_{0,i}\psi_i-\sum b_{i, I}\delta_{i,I}.
    \]
    Since $D\in E$,
    \[
        b_{1, [n]}=0, \quad b_{0,n}=b_{1, [n-1]},
    \]
    hence
    \[
        D\cdot F_5^{g-1,0}(\emptyset, n)
        =b_{g-1, \emptyset}+b_{0,n}-b_{g-1, n}
        =b_{1, [n]}+b_{0,n}-b_{1, [n-1]}=0.
    \]
    Therefore, $f$ also contracts $F_5^{g-1,0}(\emptyset, n)$. It is straightforward to see that $F_5^{g-1,0}(\emptyset, n)\not \in C$, so such $f$ does not exist. 
\end{proof}

Now we will prove that, in positive characteristic, $F_3^{1}([n])$ on $\M{3}{n}$ is contractible, i.e., there exists a projective contraction $f:\M{3}{n}\to X$ whose relative closed cone of curves is exactly the extremal ray spanned by $F_3^{1}([n])$. 

\begin{thm}\label{thm:M3ncont}
    Assume that the base field has positive characteristic. Then there exists a divisorial contraction $f:\M{3}{n}\to X$ of relative Picard number $1$, whose relative closed cone of curves is precisely the extremal ray spanned by $F_3^{1}([n])$. More precisely, $f$ is an isomorphism outside $\Delta_{2,\emptyset}$ and restricts to the projection
    \[
    \Delta_{2,\emptyset}\simeq \M{2}{1}\times \M{1}{n+1}\to\Mg{2}\times \M{1}{n+1}
    \]
    on $\Delta_{2,\emptyset}$.

\end{thm}

\begin{proof}
    First, we prove all assertions except the last. We proceed by induction. Consider first the case $n=1$. By the proof of \cref{thm:cone31}, the face of the nef cone contracting $F_3^{1}([n])$ is generated by
    \[
        \lambda,\quad 12\lambda-\delta_{\mathrm{irr}},\quad \psi_1,\quad 10\lambda-\delta_{\mathrm{irr}}+2\psi_1-2\delta_{1, \emptyset}.
    \]
    The first three divisors are semiample by the proof of \cref{prop:charpsemi}. Hence, to prove the existence of such a contraction, it is enough to show that there is another divisor in this cone that is also semiample. Let $D$ be any divisor in the interior of the cone. Then $D$ intersects positively with every F-curve except $F_3^{1}([1])$, so by a known case of F-conjecture, $D-\epsilon\,\delta_{1, [1]}$ is ample for sufficiently small $\epsilon>0$. In particular, the exceptional locus of $D$ is contained in $\Delta_{1, [1]}$, so by \cite{Ke99}, it is enough to prove that $D|_{\Delta_{1, [1]}}$ is semiample. This follows from \cref{lem:bdry} and \cref{prop:charpsemi}.

    Now consider the induction step. We need to produce a codimension-$1$ subcone of the nef cone, intersecting trivially with $F_3^{1}([n])$, consisting entirely of semiample divisors. The proof is exactly the same as the proof of \cref{thm:type3b}, using the fact that $\psi_i$ and semigroup kappa classes are semiample (cf.~\cref{thm:kapnef}).

Now we prove the last assertion, which implies that $f$ is divisorial. It suffices to show that an integral subcurve $C \subseteq \M{3}{n}$ is contained in a fiber of
\[
\M{2}{1}\times \M{1}{n+1}\longrightarrow \Mg{2}\times \M{1}{n+1}
\]
if and only if its class in $\mathrm{A}_1(\M{3}{n})$ is proportional to $F_3^{1}([n])$.

Let $C$ be an integral curve contained in a fiber. Then $C = C_0 \times \{p\}$ for some integral subcurve $C_0 \subseteq \M{2}{1}$ contracted by $\pi:\M{2}{1}\to \Mg{2}$. Since the Picard rank of $\M{2}{1}$ (resp.~$\Mg{2}$) is $3$ (resp.~$2$), the classes of such curves are proportional. Noting that $F_3^0(\{1\})$ is contracted by $\pi$, we see that $[C_0]$ is proportional to $[F_3^0(\{1\})]$. As $C$ is the image of $C_0$ under the clutching map attaching the curve corresponding to $p$, it follows that $[C]$ is proportional to $[F_3^{1}([n])]$.

Conversely, let $C$ be an integral curve on $\M{3}{n}$ whose class is proportional to $[F_3^{1}([n])]$. Since
\[
\delta_{2,\emptyset}\cdot F_3^{1}([n])=-1,
\]
we have
\[
\delta_{2,\emptyset}\cdot C<0,
\]
so $C$ is contained in $\Delta_{2,\emptyset}\simeq \M{2}{1}\times \M{1}{n+1}$.

Let $\pi_1,\pi_2$ be the projections from $\M{2}{1}\times \M{1}{n+1}$. We claim that $\pi_2(C)$ is a point, hence $C$ is of the form $C_0\times \{p\}$. This follows from
\[
D\cdot \pi_2(C)=\pi_2^\ast D\cdot C=0
\]
for every $D\in \Pic(\M{1}{n+1})$.

Let $i:\M{2}{1}\times \M{1}{n+1}\to \M{3}{n}$ be the clutching map. Define
\[
S=\{\delta_{0,I}\mid |I|\ge 2,\ I\subseteq [n]\}\ \cup\ \{\delta_{1,I}\mid I\subsetneq [n],\ I\ne \emptyset\}\ \cup\ \{\psi_i\mid 1\le i\le n\}\ \cup\ \{\lambda\}.
\]
Note that each $D\in S$ can be regarded both as a divisor on $\M{3}{n}$ and on $\M{1}{n+1}$. To avoid confusion, view $S$ as divisors on $\M{3}{n}$ and let $S'$ be the corresponding set on $\M{1}{n+1}$; for $D\in S$ denote by $D'\in S'$ the corresponding divisor. For $D\in S$, we have $D\cdot F_3^{1}([n])=0$, hence $D\cdot C=0$. Moreover, for $D\in S\setminus\{\lambda\}$,
\[
i^\ast D=\pi_2^\ast D'
\]
on $\M{2}{1}\times \M{1}{n+1}$. Therefore
\[
D'\cdot \pi_2(C)=\pi_2^\ast D'\cdot C=i^\ast D\cdot C=D\cdot C=0.
\]
For $\lambda$,
\[
0=\lambda\cdot C=i^\ast \lambda\cdot C=(\pi_1^\ast \lambda+\pi_2^\ast \lambda)\cdot C.
\]
Since both $\pi_1^\ast \lambda$ and $\pi_2^\ast \lambda$ are nef, we obtain
\[
\lambda\cdot \pi_2(C)=\pi_2^\ast \lambda\cdot C=0.
\]
Hence, for every $D'\in S'$,
\[
\pi_2^\ast D'\cdot C=0.
\]

To finish the claim, it remains to show that $S'$ spans $\Pic(\M{1}{n+1})$. By \cite[Theorem~4(c)]{AC09}, the boundary divisors span $\Pic(\M{1}{n+1})$. Among the boundary divisors, all except $\delta_{0,[n+1]}$ lie in $S'$ since $12\lambda=\delta_{\text{irr}}$. Moreover, $\psi_1\in S'$, so \cite[Theorem~4(c)]{AC09} yields
\[
\delta_{\mathrm{irr}}+12\,\delta_{0,[n+1]}\in \operatorname{Span} S'.
\]
Thus $S'$ spans $\Pic(\M{1}{n+1})$, proving the claim.

Consequently, $C=C_0\times \{p\}$ for some integral curve $C_0\subseteq \M{2}{1}$. Since $\lambda$ and $12\lambda-\delta_{\mathrm{irr}}$ intersect trivially with $F_3^{1}([n])$, their pullbacks to $\M{2}{1}$ intersect trivially with $C_0$. But these two classes span $\pi^\ast\Pic(\Mg{2})\subseteq \Pic(\M{2}{1})$. Therefore, $C_0$ is contained in a fiber of $\pi:\M{2}{1}\to \Mg{2}$.

\end{proof}

%% file: sections/discussion.tex
In this section, we revisit some questions previously presented in the body of the paper, along with several new ones introduced here. They are organized into four themes.

\subsection{Non-semiample nef divisors on \texorpdfstring{$\M{g}{n}$}{\texttwoinferior}}

In this subsection, we assume that the characteristic of the base field is $0$.

\cref{thm:main1} shows that a large portion of the nef cone of $\M{g}{n}$ is non-semiample, and that if $g=2$ or $3$, the subset of non-semiample nef divisors attains the smallest possible codimension. Hence, a natural question arising from \cref{thm:main1} is \cref{qes:nonsemi}, which asks whether the same is also true for higher genus. Unfortunately, the non-semiampleness criterion \cref{thm:semin} does not appear to be sufficient for answering the question when $g\ge 4$, since it requires the divisor to contract a type~3 F-curve. Therefore, to address this question, one would needs to find other examples of non-semiample nef divisors that can be utilized in this context, such as the divisor on $C\times C$ used in the proof of \cref{thm:semin}.

The following question about semiample divisors on $\M{g}{n}$ is natural, as the answer is known in other cases.

\begin{qes}\label{qes:semi01g}
    Is every nef divisor on $\M{0}{n}$, $\M{1}{n}$, and $\Mg{g}$ semiample?
\end{qes}

\cref{qes:semi01g} for $\M{0}{n}$ is considered in \cite{Fe15, MS19}, and it is known to hold for $n\le 6$, since $\M{0}{n}$ is log Fano for this range, as well as for symmetric divisors when $n\le 19$ by \cite{MS19}. The case of \cref{qes:semi01g} for $\Mg{g}$ was posed in \cite{Far09}. The author is not aware of other sources that have posed \cref{qes:semi01g} for $\M{1}{n}$. The motivation for this question is to find a non-semiample divisor that is unrelated to $\psi$-classes. Note that \cref{thm:semin} originated from the proof of the non-semiampleness of $\psi$-classes for higher genus in \cite{Ke99}, whereas $\psi$-classes are semiample in genus~$1$ by \cref{thm:psisemi}.

\subsection{Semiample divisors in positive characteristic}\label{subsec:pos}

Thanks to \cite{Ke99}, it is easier to prove the semiampleness of certain divisors in positive characteristic, and indeed, there are more semiample divisors in positive characteristic. This naturally leads to the following question.

\begin{conj}\label{qes:semip}
    Over a field of positive characteristic, every nef divisor on $\M{g}{n}$ is semiample.
\end{conj}

If true, this conjecture would reveal a drastic difference in the nature of $\M{g}{n}$ between characteristic~$0$ and positive characteristic. We note the following fact regarding \cref{qes:semip}.

\begin{prop}[char $k\ne 2$]
    \cref{qes:semip} for genus $0$ implies \cref{qes:semip}.
\end{prop}

\begin{proof}
    We proceed by induction on $\dim \M{g}{n}$. Let $D$ be a nef divisor on $\M{g}{n}$. Since the conjecture is known for $g=0$, we may assume $g\ge 1$. If $g\ge 2$, then by \cite[Theorem 0.9]{GKM02}, either $D=\pi_i^\ast D'$ for some projection map $\pi_i$, or $D$ is big and the exceptional divisor is contained in the boundary of $\M{g}{n}$. In the first case, since $D'$ is also nef, the induction hypothesis implies that $D$ is semiample. In the second case, by \cite[Theorem 0.2]{Ke99}, it suffices to prove that the restriction of $D$ to boundary divisors are semiample. By \cref{lem:bdry}, it is enough to show that the pullbacks of $D$ along the attaching maps $\theta:\M{g'}{n'}\to \M{g}{n}$ and $\xi:\M{g-1}{n+2}\to \M{g}{n}$ are semiample. This also follows from the induction hypothesis. 
    
    The case $g=1$ is almost the same, except that in the first case we have \( D=\pi_1^\ast D_1+ \pi_2^\ast D_2\), where $\pi_1:\M{1}{n}\to \M{1}{S}$ and $\pi_2:\M{1}{n}\to \M{1}{S^c}$, and $D_1$ (resp.~$D_2$) is a nef divisor on $\M{1}{S}$ (resp.~$\M{1}{S^c}$).
\end{proof}

\subsection{Modular interpretation of morphisms defined by semiample divisors}\label{subsec:morph}

In \cref{rmk:cont}, we motivated \cref{qes:psisemi} by relating it to contracting an F-curve on $\M{1}{n}$. There is another motivation for seeking a modular interpretation of other morphisms associated to semiample divisors. As mentioned, many divisors are semiample only in positive characteristic; for example, $\psi_i$ is such a case by \cite{Ke99}. However, to the author’s knowledge, the modular interpretation of morphisms corresponding to $\psi_i$ is not known, except in genus~$0$, where the corresponding map is given by Kapranov’s construction \cite{Ka92, Ka93}. If we could obtain a modular interpretation of such morphisms defined only in positive characteristic, we might uncover the origin of this difference. 

A good analogue in the moduli of abelian varieties is the complete subvariety problem \cite{GMMT25}. In that setting, there are more complete subvarieties of $A_g$ in positive characteristic, given by the locus of abelian varieties with $p$-rank $0$. Here, the difference arises from the $p$-rank, which is meaningful only in positive characteristic.

Regarding this, we ask for the modular interpretation of the morphism corresponding to semigroup Kappa divisors (cf.~\cref{qes:mult}). Moreover, it would be interesting to find a modular interpretation of the morphism on $\M{3}{n}$ that contracts only $F_3^1([n])$.

Finally, we note that in \cite{SJ25}, the authors observed that the effectiveness of a divisor class depends on the characteristic of the base field.

\subsection{Extremality of boundary strata of higher codimension}

There is a body of literature devoted to extremal cycles of $\M{g}{n}$. However, relatively little is known in the case of higher codimension (see, e.g., \cite{Bla22}). One of the obstacles is that, unlike in low codimension, there are many relations between higher-dimensional boundary strata, which makes the investigation significantly more complicated.

In this paper, we focus on the special case of $1$-dimensional boundary strata, namely, F-curves. We study this case using the dual cone of nef divisors of $\M{g}{n}$. However, even in this setting, the case of type~5 F-curves (cf.~\cref{conj:exttype5}) remains unresolved, and for the type~6 case, we do not even have a conjectural description of the extremal curves. Hence, the following question is worth investigating.

\begin{qes}\label{qes:exthigh}
    Which of the boundary strata of $\M{g}{n}$ are extremal?
\end{qes}